\documentclass[11pt]{amsart}
\usepackage[margin=0.96in]{geometry}

\usepackage{comment,todonotes}
\usepackage{amsmath, amsthm, amssymb,amsfonts}
\usepackage{mathtools}
\usepackage{hyperref,pdfsync}
\usepackage{bbm, dsfont}
\usepackage{cleveref}

\newtheorem{theorem}{Theorem}[section]
\newtheorem{corollary}[theorem]{Corollary}
\newtheorem{lemma}[theorem]{Lemma}
\newtheorem{proposition}[theorem]{Proposition}
\newtheorem*{definition}{Definition}

\newtheorem{question}{Question}

\numberwithin{equation}{section}

\begin{document}

\title[Space curves with totally positive torsion]{A new proof of the description of the convex hull of space curves with totally positive torsion}
\author[J. de Dios Pont, P.~Ivanisvili, and J.~Madrid]{J.~de Dios Pont, P.~Ivanisvili, and J.~Madrid}
\thanks{J.D. was partially supported by the NSF grant DMS-1764034; P.I. was partially supported by the NSF grants DMS-2152346 and CAREER-DMS-2152401}

\address{Department of Mathematics, University of California, Los Angeles, CA}
\email{jdedios@math.ucla.edu \textrm{(J.\ de Dios Pont)}}

\address{Department of Mathematics, University of California, Irvine, CA}
\email{pivanisv@uci.edu \textrm{(P.\ Ivanisvili)}}

\address{Department of Mathematics, University of California, Los Angeles, CA}
\email{jmadrid@math.ucla.edu \textrm{(J.\ Madrid)}}

\makeatletter
\@namedef{subjclassname@2010}{
  \textup{2010} Mathematics Subject Classification}
\makeatother
\begin{abstract} 
We give new proofs of the description convex hulls of  space curves $\gamma : [a,b] \mapsto \mathbb{R}^{d}$ having totally positive torsion. These are curves such that all the leading principal minors of $d\times d$  matrix $(\gamma', \gamma'', \ldots, \gamma^{(d)})$ are positive. 
 In particular, we recover  parametric representation of the boundary of the convex hull,  different formulas for its surface area and the volume of the convex hull,   and the  solution to a general moment problem corresponding to $\gamma$. 
\end{abstract}
\maketitle

\section{Introduction and a summary of main results}

Convex hull of a set $K \subset \mathbb{R}^{d}$ is defined as 
\begin{align*}
    \mathrm{conv}(K) = \left\{ \sum_{j=1}^{m} \lambda_{j}x_{j}, \, x_{j} \in K, \, \sum_{j=1}^{m} \lambda_{j}=1, \, \lambda_{j} \geq 0,\,  j=1, \ldots, m\; \text{for all}\; m\geq 1\right\}.
\end{align*}
Describing  the convex hull of a given set $K$ is a basic problem in mathematics.  By imposing  additional geometric structures on $K$ one may hope to give a {\em simpler} description of $\mathrm{conv}(K)$. Perhaps a good starting point is when $K$ is a space curve which is the topic of our paper.

Let $[a,b]$ be an interval in $\mathbb{R}$, and let $\gamma_{1}(t), \ldots, \gamma_{n+1}(t)$ be real valued  functions on $[a,b]$.  We start with two main questions which are ultimately related to each other.
\begin{question}\label{que2} Describe the boundary of the  convex hull of $\gamma([a,b])$, where 
$$
\gamma(t)=(\gamma_{1}(t), \ldots, \gamma_{n+1}(t)), \quad t \in [a,b].
$$
\end{question}

The next question, known as the {\em general moment problem} \cite{Kem1968, Karlin1, krn},  is a certain probabilistic reformulation of Question~\ref{que2}.
\begin{question}\label{que1} Find
\begin{align}
 M^{\mathrm{sup}}(x_{1}, \ldots, x_{n}) &\stackrel{\mathrm{def}}{=} \sup \;  \{ \mathbb{E}\gamma_{n+1}(Y) \, :\, \mathbb{E}\gamma_{1}(Y)=x_{1}, \ldots , \mathbb{E} \gamma_{n}(Y)=x_{n}\}, \label{sup1}\\
 M^{\mathrm{inf}}(x_{1}, \ldots, x_{n}) &\stackrel{\mathrm{def}}{=} \inf \;  \{ \mathbb{E}\gamma_{n+1}(Y) \, :\, \mathbb{E}\gamma_{1}(Y)=x_{1}, \ldots , \mathbb{E} \gamma_{n}(Y)=x_{n}\}, \label{inf1}
\end{align}
where supremum or infimum is taken over all  random variables $Y$ with values in $[a,b]$   such that $\gamma_{j}(Y)$ are measurable for all $j$, $1\leq j \leq n+1$. 
\end{question}

The answers to both of these questions are given in terms of {\em lower and upper principal representations} in two remarkable monographs \cite{krn, Karlin1} (see also a brief survey \cite{pin01})  under the assumption (A1) which says  that the sequences $(1, \gamma_{1}(t), \ldots, \gamma_{n}(t))$ and $(1, \gamma_{1}(t), \ldots, \gamma_{n+1}(t))$ are $T_{+}$-systems on $[a,b]$, we refer the reader to Subsection~\ref{markovs} for more details.

In this paper we give a new self-contained geometric approach to both of these questions for a subclass of (A1),  curves  with so called  {\em totally positive torsion}. 
%As a corollary,  we obtain 1) formulas for the volume of $\mathrm{conv}(\gamma([a,b]))$; 2) {\em Carath\'eodory number} of  $\mathrm{conv}(\gamma([a,b]))$; 3) any single affine hyperplane intersects $\gamma([a,b])$ in at most $n+1$ points
\begin{definition}
A curve $\gamma \in C^{n+1}((a,b), \mathbb{R}^{n+1}) \cap C([a,b], \mathbb{R}^{n+1})$ is said to have totally positive torsion if all the leading  principal minors of the matrix 
\begin{align}\label{mm22}
(\gamma'(t), \gamma''(t), \ldots, \gamma^{(n+1)}(t))
\end{align}
are positive for all $t \in (a,b)$. 
\end{definition}

Perhaps an instructive example to keep in mind is $\gamma(t)=(t, t^{2}, \ldots, t^{n}, \gamma_{n+1}(t))$ where the total positivity of the torsion on $(a,b)$ is the same as  $\gamma_{n+1}^{(n+1)}(t)>0$ on $(a,b)$.

In fact the only property that will be needed from the principal minors  of the matrix (\ref{mm22}) is that they are non-vanishing.  Indeed, we can consider an invertible linear  image of $\gamma$, namely a new curve $t \mapsto (\varepsilon_{1}\gamma_{1}(t), \ldots, \varepsilon_{n+1} \gamma_{n+1}(t))$ with an appropriate choice of signs $\varepsilon_{j} = \pm 1$ and reduce the study of the convex hulls  to the curves with totally positive torsion (an invertible linear transformation $T$ maps convex hull of a set $K$   to the convex hull of the image $T(K)$).

%We also provide applications  refining the recent result due to Carlen, Frank, Lieb ~\cite{Car1} on sharpening the classical triangle inequalities in $L^{p}$ spaces for many functions.  

In Section~\ref{istoria} we provide an overview of the literature on results related to Questions 1 and 2.  Section~\ref{glavnaya} is devoted to the statements of  main results of the paper, and Section ~\ref{damtkiceba} contains the proofs. Here we give a short summary of the theorems  that we recover in this paper and that were previously known in \cite{krn, Karlin1}. The results we state hold in $\mathbb{R}^{n+1}$ for all $n\geq 1$, and all  space curves $\gamma : [a,b] \to \mathbb{R}^{n+1}$ with totally positive torsion.  Set $\bar{\gamma}(t) \stackrel{\mathrm{def}}{=} (\gamma_{1}(t), \ldots, \gamma_{n}(t))$, and let us denote by  $\mathrm{conv}(\gamma([a,b]))$  the convex hull of the image of $[a,b]$ under the map $\gamma$.
\subsection*{Summary of the results:}
\begin{itemize}
    \item[(1)] Boundary of the convex hull of $\gamma([a,b])$ will be given in a parametric form.
    \item[(2)] Explicit diffeomorphism will be constructed between the interior of simplicies  and the interior of the convex hull of $\gamma([a,b])$
    \item[(3)]   Formulas for the surface area of the boundary of the convex hull of $\gamma([a,b])$ will be obtained, Corollary~\ref{area1}, and  different formulas for the volume of the convex hull will be presented, Corollary~\ref{provolume}.
    \item[(4)] Any single affine hyperplane intersects the space curve  $\gamma :[a,b] \to \mathbb{R}^{n+1}$  in at most $n+1$ points.  Minimal number $k$  points required to represent any point $x \in \mathrm{conv}(\gamma([a,b]))$ as a convex combination of $k$ points of $\gamma([a,b])$ is at most $\lfloor \frac{n+3}{2}\rfloor$. Moreover, $k = \lfloor \frac{n+3}{2}\rfloor$  for any interior point of $\mathrm{conv}(\gamma([a,b]))$.
    \item[(5)]
    Parametric representations will be given for  functions $M^{\sup}$ and $M^{\inf}$. The obtained parametric forms change depending on whether $n$ is even or odd. 
    
    \textup{(i)} If $n$ is even then 
\begin{align*}
&M^{\sup}\left(\lambda_{0} \bar{\gamma}(b)+\sum_{j=1}^{\frac{n}{2}} \lambda_{j} \bar{\gamma}(x_{j}) \right) = \lambda_{0} \gamma_{n+1}(b)+\sum_{j=1}^{\frac{n}{2}}\lambda_{j} \gamma_{n+1}(x_{j}),\\
&M^{\inf}\left(\lambda_{0} \bar{\gamma}(a)+\sum_{j=1}^{\frac{n}{2}} \lambda_{j} \bar{\gamma}(x_{j}) \right) = \lambda_{0} \gamma_{n+1}(a)+\sum_{j=1}^{\frac{n}{2}}\lambda_{j} \gamma_{n+1}(y_{j}),
\end{align*}
for all $\lambda_{0}, \lambda_{j} \in [0,1], x_{j} \in [a,b]$, $j=1, \ldots, \frac{n}{2}$ with $\sum_{0\leq k \leq \frac{n}{2}} \lambda_{k}=1$.

\textup{(ii)} If $n$ is odd then 
\begin{align*}
&M^{\sup}\left(\lambda_{0}\bar{\gamma}(a)+\lambda_{1}\bar{\gamma}(b)+\sum_{j=2}^{\frac{n+1}{2}} \lambda_{j} \bar{\gamma}(x_{j}) \right) = \lambda_{0} \gamma_{n+1}(a)+\lambda_{1} \gamma_{n+1}(b)+\sum_{j=2}^{\frac{n+1}{2}}\lambda_{j} \gamma_{n+1}(x_{j}), \\
&M^{\inf}\left(\sum_{j=1}^{\frac{n+1}{2}} \beta_{j} \bar{\gamma}(x_{j}) \right) = \sum_{j=1}^{\frac{n+1}{2}}\beta_{j} \gamma_{n+1}(x_{j}),
\end{align*}
for all $\lambda_{0}, \lambda_{j}, \beta_{j} \in [0,1], x_{j} \in [a,b]$, $j=1, \ldots, \frac{n+1}{2}$ with $\sum_{0\leq j \leq \frac{n+1}{2}} \lambda_{j}=\sum_{1\leq j \leq \frac{n+1}{2}}\beta_{j}=1$.
   
 \item[(6)] Explicit random variables $Y$ will be constructed which attain supremum and  infimum  correspondingly in (\ref{sup1}) and (\ref{inf1}) for each given $x = (x_{1}, \ldots, x_{n})$ from the domain of definition of $M^{\sup}$ and $M^{\inf}$. 
\end{itemize}
We will also see that 
 \begin{align*}
 \partial\,  \mathrm{conv}(\gamma([a,b]))=\{(x,M^{\mathrm{sup}}(x)), x \in \mathrm{conv}(\bar{\gamma}([a,b]))\} \cup \{(x,M^{\mathrm{inf}}(x)), x \in \mathrm{conv}(\bar{\gamma}([a,b]))\}, 
\end{align*}
i.e., the {\em upper hull} of $\mathrm{conv}(\gamma([a,b]))$ coincides with the graph of $M^{\sup}$, and the lower hull with the graph of $M^{\inf}$. Besides this summary, we also recover several results previously known to Karlin--Sharpley \cite{Karlin2} for {\em moment curves} using our techniques (see  Corollary~\ref{nobel2}).  In Proposition~\ref{sensitive}, we also show that the results obtained in this paper are sensitive to the assumption on a curve having totally positive torsion.

%Then the optimization problem in Question~\ref{que1} takes a compact form 
%\begin{align}
%M^{\mathrm{sup}}(x) &\stackrel{\mathrm{def}}{=}  \sup_{a\leq Y\leq b} \{ \mathbb{E}\gamma_{n+1}(Y) \, :\, \mathbb{E}\bar{\gamma}(Y)=x\};\label{sup1}\\
%M^{\mathrm{inf}}(x) &\stackrel{\mathrm{def}}{=}  \inf_{a\leq Y\leq b} \{ \mathbb{E}\gamma_{n+1}(Y) \, :\, \mathbb{E}\bar{\gamma}(Y)=x\}.\label{inf1}
%\end{align}

\subsection{What is known about Questions 1 and 2? } \label{istoria}
In what follows we set  $x \stackrel{\mathrm{def}}{=}(x_{1}, \ldots, x_{n}) \in \mathbb{R}^{n}$, and $\mathbb{E} \bar{\gamma}(Y) \stackrel{\mathrm{def}}{=} (\mathbb{E}\gamma_{1}(Y), \ldots, \mathbb{E}\gamma_{n}(Y))$. 
We remark that both $M^{\mathrm{\sup}}$ and $M^{\mathrm{\inf}}$ depend on $n \geq 1$,  $x \in \mathbb{R}^{n}, [a,b] \subset \mathbb{R}$, and $\gamma$. We shall remind the basic fact that  the convex hull of a compact set is compact. For simplicity we shall use the symbol $M$ for $M^{\mathrm{sup}}(x)$. 

There are series of results describing $M$ for some particular $\gamma$. A common goal is to have a parametric representation for it. However, as soon as $n$ is large it becomes difficult to find parametric representation for $M$ in such generality.  

\subsubsection{Convex envelopes and Carath\'eodory number}
 Under some mild assumptions on $\gamma$, say $\gamma$ is continuous on $[a,b]$ is sufficient  (see \cite{Kem1968, Rog1958}), $M$ is defined on $\mathrm{conv}(\bar{\gamma}([a,b]))$. Moreover, for any $x \in \mathrm{conv}(\bar{\gamma}([a,b]))$, $M(x)$ is the solution of the {\em dual problem} 
\begin{align}\label{dual}
M(x) = \inf_{d_{0} \in \mathbb{R}, d \in \mathbb{R}^{n}} \{ d_{0}+ \langle d,x \rangle \;\; \text{such that}\;\; d_{0}+ \langle d, \bar{\gamma}(t) \rangle   \geq \gamma_{n+1}(t)\;  \text{for all} \;  t \in [a,b]\},
\end{align}
where $\langle a,b\rangle$ denotes the dot product in $\mathbb{R}^{n}$. Thus $M$ is the minimal concave function defined on $ \mathrm{conv}(\bar{\gamma}([a,b]))$ with the obstacle condition $M(\bar{\gamma}(t)) \geq \gamma_{n+1}(t)$ for all $t \in [a,b]$. So the graph $(x,M(x))$, $x \in  \mathrm{conv}(\bar{\gamma}([a,b]))$ belongs to the boundary of $\mathrm{conv} (\gamma([a,b]))$. Carath\'eodory's theorem says that $(x,M(x))$ is convex combination of at most $n+2$ points  from $\gamma([a,b])$. However, due to the fact $(x, M(x)) \in \partial\,  \mathrm{conv} (\gamma([a,b]))$, one can see that $n+1$ points  suffice by considering any affine hyperplane $H$ supporting $\mathrm{conv} (\gamma([a,b]))$ at $(x,M(x))$. Since $\gamma([a,b])$ lies on one side of $H$, it follows that the points, whose convex combination is $(x,M(x))$, must lie in $H$, and we can apply Carath\'eodory's theorem to  $H \cap \gamma([a,b])$ in $n+1$ dimensional space $H$. 
This leads us to another representation 
\begin{align}\label{carath}
M(x) = \sup_{\sum_{j=1}^{n+1}c_{j} \bar{\gamma}(t_{j})=x} \left\{ \sum_{j=1}^{n+1} c_{j} \gamma_{n+1}(t_{j})\; :\; \sum_{j=1}^{n+1}c_{j}=1, \; c_{\ell} \geq 0,\; t_{\ell} \in [a,b], \;1\leq \ell \leq n+1\right\}.
\end{align}

Probabilistic way of looking at (\ref{carath}) is that the supremum and infimum  in (\ref{sup1}) and  (\ref{inf1}) is attained on random variables $Y$ whose density is the sum of delta masses on at most $n+1$ points in $[a,b]$, i.e., $\sum_{j=1}^{n+1}c_{j}\delta_{t_{j}}$, with $t_{j} \in [a,b]$ for all $j=1, \ldots, n+1$.  

A direction of research focuses on understanding for which curves $\gamma$, the number $n+1$ appearing in $\sum_{j=1}^{n+1}c_{j}\delta_{t_{j}}$ can be made smaller. As we just described this is related to the following question: {\em 
given a curve $\gamma :[a,b] \to \mathbb{R}^{n+1}$, and a point $y \in \partial\, \mathrm{conv}(\gamma([a,b]))$, find the smallest number of points $b(y)$  on $\gamma([a,b])$ whose convex combination coincides with $x$.}
The integer $b(y)$ is called  Carath\'eodory number for $y$, and it is defined for all $y \in \mathrm{conv}(\gamma([a,b]))$.  Carath\'eodory number $b(\gamma)$ of a set $\gamma([a,b])$ is defined as 
\begin{align}\label{karate1}
b(\gamma) \stackrel{\mathrm{def}}{=}\sup_{x \in \mathrm{conv}(\gamma([a,b]))}b(x).
\end{align} 

By Carath\'eodory's theorem $b(\gamma) \leq n+2$ for curves in $\mathbb{R}^{n+1}$. For certain curves $\gamma$, the number $b(\gamma)$ can be strictly smaller than $n+2$. Fenchel's theorem \cite{Fenchel, Hanner} asserts that if the compact set $\gamma([a,b])$ cannot be separated by a hyperplane into two non-empty disjoint sets then $b(\gamma)\leq n+1$. In particular, for continuous curves $\gamma$ over closed intervals $[a,b]$ the Carath\'eodory's number is at most $n+1$ giving one more justification of (\ref{carath}) for continuous maps $\gamma$.   See \cite{Baran} where Carath\'eodory number and an extension of Fenchel's theorem  is studied for certain type of sets in $\mathbb{R}^{n+1}$.

\subsubsection{A Convex Optimization Approach}
Another direction of research reduces (\ref{dual}) to what is called {\em positive semidefinite optimization problem} under the assumption  
$$
\gamma(t) = (t,t^{2}, \ldots, t^{n}, \mathbbm{1}_{I}(t)),
$$
 where $I$ is an interval in $\mathbb{R}$. 
 %and $[a,b]$ is replaced by $\mathbb{R}$, or $\mathbb{R}_{+}=[0, \infty)$.
 Finding upper or lower bounds on $\mathbb{E} \mathbbm{1}_{I}(Y) = \mathbb{P}(Y \in I)$ given the first $n$ moments of $Y$ is of important interests as it would  refine the classical Chebyshev and Markov inequalities. To give a feeling how the corresponding positive semidefinite optimization problem looks like we cite Theorem~11 in \cite{Berts}:  the tight upper bound on $\mathbb{P}(Y \geq 1)$ over all nonnegative random variables $Y$ given the first $n$ moments $\mathbb{E}Y^{j}=x_{j}$, $1\leq j \leq n$ coincides with 
 
 \begin{align*}
 M^{\mathrm{sup}}(x)\,  =\,  \min_{d_{0}, \ldots, d_{n} \in \mathbb{R}} \quad  d_{0}+\sum_{j=1}^{n}d_{j} x_{j} 
 \end{align*}
 Subject to
 \begin{align*}
 \quad &0 = \sum_{i,j\, :\, i+j=2\ell-1} t_{ij}, \qquad \qquad  \quad \; \,\ell=1, \ldots, n,\\
 &(d_{0}-1)+\sum_{j=\ell}^{n} d_{j} \binom{j}{\ell}=t_{00},\\
 &\sum_{j=\ell}^{n} d_{j} \binom{j}{\ell} = \sum_{i,j\, :\, i+j=2\ell}t_{ij}, \quad \quad \; \,\ell=1, \ldots, n,\\
 &0 = \sum_{i,j\, :\, i+j=2\ell-1} z_{ij}, \qquad  \qquad \quad \; \, \ell = 1, \ldots, n,\\
 &\sum_{j=0}^{\ell} d_{j} \binom{n-j}{\ell-j} = \sum_{i,j\, :\, i+j=2\ell} z_{ij} \quad \ell=0, \ldots, n,\\
 &T, Z \geq 0,
 \end{align*}
 where $T, Z \geq 0$ means that the matrices $T=\{t_{ij}\}_{i,j=0}^{n}, Z = \{z_{ij}\}_{i,j=0}^{n}$ are positive semidefinite. 
 %We remark that there is a typo in the statement of part a), Theorem~11 in \cite{Berts}, namely, the summation   in $(d_{0}-1)+\sum_{j=\ell}^{n} d_{j} \binom{j}{\ell}=t_{00}$ is over the range $j=\ell, \ldots, n$ instead of $j=1, \ldots, n$. 
 
 The advantage of having such a semidefinite optimization problem is that it can be solved in a {\em polynomial time}. However, it is not clear to us how practical are these results if one wants to verify bounds $M(x) \leq R(x)$ for a given function $R$ and all $x$ in $\mathrm{conv}(\overline{\gamma}([0,1]))$. In \cite{Berts} the authors provide explicit formulas for the tight upper bound on $\mathbb{P}(Y>\lambda)$ for  $n=3$ over all nonnegative random variables with given first 3 moments.

\subsubsection{Tchebysheff systems, convex curves,  and Markov moment problem}\label{markovs}

The system of continuous functions $(\gamma_{0}(t), \ldots, \gamma_{n}(t))$,  on an interval $[a,b]$ is called Tchebysheff system (or $T$-system) if any nontrivial linear combination $\sum_{j=0}^{n} a_{j} \gamma_{j}(t)$ has at moat $n$ roots on $[a,b]$.  As the monographs \cite{krn, Karlin1} deal with general Markov moment problem with arbitrary Borel measures, and in this paper we consider only probability measures,  in what follows we will be assuming that $\gamma_{0}(t)=1$  to make the presentation consistent with \cite{krn, Karlin1}. Under such an assumption the corresponding curve $t\mapsto (\gamma_{1}(t), \ldots, \gamma_{n}(t))$ is called {\em convex curve}. 

The sequence $(\gamma_{0}(t), \ldots, \gamma_{n}(t))$ is called $T_{+}$-system if 
\begin{align}\label{nudel}
\mathrm{det}(\{ \gamma_{i}(t_{j})\}_{i,j=0}^{n})>0
\end{align}
on the simplex $\Sigma = \{ a\leq t_{0}<\ldots, <t_{n} \leq b\}$. Notice that any $T$-system can be made into $T_{+}$-system just by flipping the sign in front of $\gamma_{n}$ if necessary. 
If $(\gamma_{0}(t), \ldots, \gamma_{k}(t))$ is $T_{+}$-system on $[a,b]$ for any $k=0,\ldots, n$ then the sequence $(\gamma_{0}(t), \ldots, \gamma_{n}(t))$ is called $M_{+}$-system on $[a,b]$.  Checking the positivity of the determinant (\ref{nudel}) seems a bit unpractical  as one needs to verify the inequality on the simplex
 $\Sigma$. The following proposition gives a simple sufficient criteria for the system to be $M_{+}$ system. 
 
\begin{theorem}[Chapter VIII, \cite{Karlin1}]\label{man1}
Let $\gamma_{0}(t), \ldots, \gamma_{n}(t)$ be in $C([a,b])\cap C^{n}((a,b))$.  Then for the sequence $(\gamma_{0}(t), \ldots, \gamma_{n}(t))$ to be $M_{+}$-system on $[a,b]$ it is necessary\footnote{Here $\gamma_{j}^{(0)}(t)=\gamma_{j}(t)$} that $\mathrm{det}(\{ \gamma_{i}^{(j)}(t)\}_{i,j=0}^{k})\geq 0$  on $(a,b)$ for all $k=0,\ldots, n$, and it is sufficient that $\mathrm{det}(\{ \gamma_{i}^{(j)}(t)\}_{i,j=0}^{k})> 0$ on $(a,b)$ for all $k=1,\ldots, n$. 
\end{theorem}

We say that  $(\gamma_{1}(t), \ldots, \gamma_{n+1}(t))$ satisfies $(A1)$ condition if  $\gamma_{1}(t), \ldots, \gamma_{n+1}(t)$ are in $C([a,b])\cap C^{n+1}((a,b))$  such that
\begin{align*}
(1,\gamma_{1}(t),  \ldots, \gamma_{n}(t)) \quad \text{and} \quad (1,\gamma_{1}(t), \ldots, \gamma_{n+1}(t)) \quad \text{are} \quad T_{+}-\text{systems on} \quad [a,b] \quad (A1)
\end{align*}
Clearly if $\gamma(t) = (\gamma_{1}(t), \ldots, \gamma_{n+1}(t))$ has totally positive torsion on $(a,b)$ then the condition $(A1)$ holds by Theorem~\ref{man1}. On the other hand if the sequence $(\gamma_{0}(t), \ldots, \gamma_{n+1})$ satisfies only the assumption  (A1) then the probability distribution of a random variable  $X$  achieving supremum or infimum in Question~\ref{que1} is  given in terms of {\em upper and lower principal representations}, see Chapter III and IV in \cite{krn}, and also Proposition 2 in  a brief survey \cite{pin01}. In particular, Carath\'eodory number is at most $\lfloor \frac{n+3}{2}\rfloor$ for the curves $t\mapsto (\gamma_{1}(t), \ldots, \gamma_{n+1}(t))$ in $\mathbb{R}^{n+1}$ satisfying the assumption (A1).

A typical example of the convex curve is the moment curve 

 $$
 \gamma(t) = (t, \ldots, t^{n+1}) \in \mathbb{R}^{n+1},
 $$
 %Having a small Carath\'eodory number for $\gamma$ significantly reduces {\em complexity} of computing $M$ though it does not give any explicit descriptions for it, i.e., one still has to solve the optimization problem (\ref{carath}).  
Assume $[a,b]=[0,1]$. In~\cite{Karlin2} the authors show that if $x=(x_{1}, \ldots, x_{n})$ belongs to the interior of $\mathrm{conv}(\bar{\gamma}([0,1]))$ then $M^{\mathrm{sup}}(x)$   and  $M^{\mathrm{inf}}(x)$ are  the unique solutions $x_{n+1}$ of the linear equations 
 \begin{align}\label{nobel}
 K_{n+1}=0 \quad \text{and} \quad   S_{n+1}=0,
 \end{align}
 correspondingly, where $K_{k}, S_{k}$ are defined as 
 \begin{align}\label{Sharp1}
 S_{2k} = \det 
 \begin{pmatrix}1 & x_{1} & \ldots & x_{k}\\
 \vdots & & & \\
 x_{k} & x_{k+1} & \ldots & x_{2k}\end{pmatrix}, \quad  S_{2k+1} = \det 
 \begin{pmatrix}x_{1} & x_{2} & \ldots & x_{k+1}\\
 \vdots & & & \\
 x_{k+1} & x_{k+2} & \ldots & x_{2k+1}\end{pmatrix},
 \end{align}
 and 
 \begin{align}\label{kar1}
  K_{2k} = \det 
 \begin{pmatrix}x_{1}-x_{2} & x_{2}-x_{3} & \ldots & x_{k}-x_{k+1}\\
 \vdots & & & \\
 x_{k}-x_{k+1} & x_{k+1}-x_{k+2} & \ldots & x_{2k-1}-x_{2k}\end{pmatrix},\\
   K_{2k+1} = \det 
 \begin{pmatrix}1-x_{1} & x_{1}-x_{2} & \ldots & x_{k}-x_{k+1}\\
 \vdots & & & \\
 x_{k}-x_{k+1} & x_{k+1}-x_{k+2} & \ldots & x_{2k}-x_{2k+1}\end{pmatrix}. \nonumber
 \end{align}
%Here $x_{n+1}=B^{\sup / \inf}(x_{1}, \ldots, x_{n})$. 
 %Notice that even though $x_{n+1}$ appears in (\ref{f1}) and (\ref{f2}), the values $x_{n+1}+\frac{U_{n+1}}{U_{n-1}}$ and $x_{n+1}-\frac{L_{n+1}}{L_{n-1}}$ do not depend on $x_{n+1}$, i.e., one can set $x_{n+1}=0$. We remark that $U_{n-1}, L_{n-1}$ do not vanish due to the fact that $x$ belongs to the interior of $\mathrm{conv}(\bar{\gamma}([0,1]))$, see \cite{Karlin2}.
 
 %If $x$ belongs to the boundary of $\mathrm{conv}(\bar{\gamma}([0,1]))$ then there exists a unique representation of a point $x$ as a convex combination of points on $\bar{\gamma}([0,1])$, $x = \sum_{j=1}^{n} c_{j} \bar{\gamma}(t_{j})$ with $\sum_{j=1}^{n}c_{j}=1$, $c_{j}\geq 0$, $t_{j} \in [0,1]$ for all $j=1, \ldots, n$. In this case 
 %\begin{align*}
 %B^{\mathrm{sup}}(x) = B^{\mathrm{inf}}(x) = \sum_{j=1}^{m} c_{j} t_{j}^{n+1}. 
 %\end{align*}
 
 %Perhaps it is the case of the moment curve when we obtain the most satisfactory description for  $M^{\mathrm{sup}/ \mathrm{inf}}$ thanks to the main results in \cite{Karlin2}. It also follows from \cite{Karlin2}, that the boundary of the convex hull of $\gamma([0,1])$ is the union of the graphs of $M^{\mathrm{sup}}$ and $M^{\mathrm{inf}}$, i.e., 
 %$$
 %\partial\,  \mathrm{conv}(\gamma([0,1]))=\{(x,M^{\mathrm{sup}}(x)), x \in \mathrm{conv}(\bar{\gamma}([0,1]))\} \cup \{(x,M^{\mathrm{inf}}(x)), x \in \mathrm{conv}(\bar{\gamma}([0,1]))\}.
 %$$
 An important contribution of \cite{Karlin2}  is that the authors give complete description of $\partial \, \mathrm{conv}(\gamma([0,1]))$ which allowed them to obtain a geometric point of view on the classical orthogonal polynomials. For example, knowing the width  in $x_{n+1}$ direction of the set $\mathrm{conv}(\gamma([0,1]))$ one can recover the classical fact that among all polynomials of degree $n+1$ on $[0,1]$ with the leading coefficient  $1$ the Tchebyshev polynomials minimize the maximum of the absolute value on $[0,1]$ (Theorem 25.2 in ~\cite{Karlin2}).
 
 Karlin--Sharpley did announce an intend to settle the case when $[a,b]$ is replaced by $[-1,1]$, $\mathbb{R}^{+}$ or $\mathbb{R}$. After looking into a literature, to the best of our knowledge the corresponding results  appeared in the monograph of Karlin--Studden~\cite{Karlin1}. 
 
%Let us also mention a classical Hausdorff moment problem and its relevance to Question~\ref{que1}. Given an infinity sequence of nonnegative numbers $x_{1}, x_{2}, \ldots$ the Hausdorff moment problem asks to decide if there exists a random variable $Y$ with values in $[0,1]$ such that $\mathbb{E}Y^{j}=x_{j}$ for all $j\geq 1$.  It is well--known (see Theorem 14.4 in \cite{Karlin2}) such $Y$ exists if and only if the sequence $x_{0}=1, x_{1}, x_{2}, \ldots$ is {\em completely monotonic}, i.e., 
% \begin{align}\label{hmp}
%\sum_{j=0}^{r} (-1)^{r+j}\binom{r}{j} x_{r+s-j} \geq 0
 %\end{align}
% for all integers $r,s\geq 0$. We would like to point out that as (\ref{hmp}) involves infinity sequence, the equations in (\ref{nobel}) on the variable $x_{n+1}$ given $x_{1}, \ldots, x_{n}$  do not seem to immediately follow from (\ref{hmp}). 

In ~\cite{Sch0101} Schoenberg obtained a formula for the volume of a smooth closed\footnote{Here closed curve means $\nu(0)=\nu(2\pi)$} convex curve $\nu : [0, 2\pi] \mapsto \mathbb{R}^{n}$ in even-dimensional Euclidean space
\begin{align*}
    \mathrm{Vol}(\mathrm{conv}(\nu([0, 2\pi]))) = \pm  \frac{1}{n!(n/2)!}\int_{[0,2\pi]^{\frac{n}{2}}}\det (\nu(t_{1}), \ldots, \nu(t_{n/2}), \nu'(t_{1}), \ldots, \nu'(t_{n/2}))dt_{1}\ldots dt_{n/2},
\end{align*}
and as a corollary, using Fourier series,  he derived an isoperimetric inequality 
$$
(\mathrm{length}(\nu))^{n}\geq (\pi n)^{n/2}(n/2)! n! \mathrm{Vol}(\mathrm{conv}(\nu([0, 2\pi]))), 
$$
where $\mathrm{length}(\nu)$ denotes the Euclidean length of $\nu$,  and $\mathrm{Vol}(\cdot)$ denotes the Euclidean volume. The volumes of the convex hull of $\gamma([a,b])$, such that $\gamma(0)=0$ and the sequence  $(1, \gamma_{1}(t), \ldots, \gamma_{n}(t))$ forms the $T$-system were obtained both in odd and even dimensions in \cite{krn, Karlin1}, see for example, Theorem 6.1,  Ch. IV in \cite{Karlin1}.

\subsubsection{Other results for systems different from $T$-system} In \cite{sed1, sed2} Sedykh describes  possible {\em singularities} of the boundary of  convex hulls of a curve in $\mathbb{R}^{3}$. In~\cite{Krist1}, using tools from algebraic geometry, namely, {\em De Jonqui\`eres'  formula}, the authors compute number of {\em complex tritangent planes} of the {\em algebraic boundary} of the convex hull of an algebraic space curve  in $\mathbb{R}^{3}$ in terms of its genus and degree of the curve. Moreover, in \cite{Krist1} the authors also find an algebraic elimination method for computing  {\em tritangent planes} and {\em edge surfaces} of the boundary of the convex hulls of algebraic space curves in $\mathbb{R}^{3}$. {\em Algebraic boundary} of the convex hull of an algebraic variety was studied \cite{KRBS1}, where the authors extended several results from  \cite{Krist1} to higher dimensions. In \cite{Freed}, using topological results it is shown that the number of tritangent planes to a smooth {\em generic} curve in $\mathbb{R}^{3}$ with nonvanishing torsion is even.

Convex hulls of space curves have appeared implicitly or explicitly in other works in relation to problems not directly related to them. We do not intend to provide  the full list of references, however, let us mention some of the  examples.  Finding sharp constants in such classical estimates as John--Nirenberg inequality is  related to finding convex hulls {\em in non-convex domains} of certain space curves. In particular,  in \cite{Iv1, Iv2}, an algorithm is presented which finds the convex hull of a space curve $\gamma(t) = (t, t^{2}, f(t))$ defined on $\mathbb{R}$, under the assumption that $f'''(t)$ changes sign finitely many times (notice that the sign of $f'''$ coincides with the sign of the torsion of $\gamma(t)$). As the number of sign changes of $f'''$ increase the ``complexity'' of computing the convex hull of $\gamma(t)$ increases too. The method obtained in \cite{Iv1,Iv2} is  illustrated  on a particular example in \cite{Vasyunin} for the family of space curves  $\gamma_{\alpha}(t)=(t,t^{2}, g_{\alpha}(t))$ where $g_{\alpha}(t)$ is a parametric family of functions defined for all $\alpha>0$ as follows 
\begin{align*}
    g_{\alpha}(t) =     \begin{cases}
    -\cos(t), & |t|\leq \alpha \\
    \frac{1}{2}(t^{2}-\alpha^{2})\cos \alpha+(\sin \alpha-\alpha \cos \alpha)(|t|-\alpha)-\cos \alpha, & |t|\geq \alpha. 
    \end{cases}
\end{align*}
Notice that the quadratic part for $|t|\geq \alpha$ is chosen in such a way that $g_{\alpha} \in C^{2}(\mathbb{R})$. Clearly $g'''_{\alpha}(t)=-\sin(t)$ for $|t|\leq \alpha$, and $g'''_{\alpha}(t)=0$ for $|t|\geq \alpha$. We see that as $\alpha$ increases the number of sign changes of $g'''_{\alpha}(t)$ increases too. In \cite{Vasyunin} the upper boundary of the convex hull of the space curve $\gamma_{\alpha}(t)$, $t \in \mathbb{R}$, is found in the  non-convex parametric domain\footnote{By convex hull of $\gamma_{\alpha}$ in $\Omega_{\varepsilon}$ we mean all possible  convex  combinations of those points on $\gamma_{\alpha}$ such that the projection of the resulting convex hull of these points onto $\mathbb{R}^{2}$ lies inside $\Omega_{\varepsilon}$}.   
$\Omega_{\varepsilon} = \{ (x,y) \in \mathbb{R}^{2}\, :\, x^{2} \leq y \leq x^{2}+\varepsilon^{2}\}.$
In the limiting case $\varepsilon \to \infty$ one recovers the upper boundary of the convex hull of the space curve $\gamma_{\alpha}(t)$.

In  sharpening the triangle inequality in $L^{p}$ spaces, for each $p \in \mathbb{R} \setminus\{0\}$ the paper  \cite{IM} finds the boundary of the convex hull of a space curve $\gamma(t) = (t, \sqrt{1-t^{2}}, ((1-t)^{1/p}+(1+t)^{1/p})^{p}), t \in [-1,1]$. In~\cite{IVZ} the boundary of the convex hull of a closed space curve is described which is the union of the following three curves 
\begin{align*}
&\left(\frac{1}{t^{p}+(1-t)^{p}+1}, \frac{t^{p}}{t^{p}+(1-t)^{p}+1}, \frac{(1+t)^{p}}{t^{p}+(1-t)^{p}+1}\right), \quad t \in [0,1];\\
&\left(\frac{(1-t)^{p}}{t^{p}+(1-t)^{p}+1}, \frac{1}{t^{p}+(1-t)^{p}+1}, \frac{(2-t)^{p}}{t^{p}+(1-t)^{p}+1}\right), \quad t \in [0,1];\\
&\left(\frac{t^{p}}{t^{p}+(1-t)^{p}+1}, \frac{(1-t)^{p}}{t^{p}+(1-t)^{p}+1}, \frac{|1-2t|^{p}}{t^{p}+(1-t)^{p}+1}\right), \quad t \in [0,1].
\end{align*}

\subsection*{Acknowledgments}
We are grateful to Pavel Zatitskiy for pointing our attention to the reference~\cite{krn}. 
The authors would like to thank V.~Sedykh for providing references  on topological results on the convex hulls of space curves.

%\section{Informal description of the main results}
%[Jaume do you want to add some pictures illustrating the main %results of the paper?]. 

\section{Statements of main results}\label{glavnaya}
\label{sec:Statements}
For any $v=(v_{1}, \ldots, v_{d}) \in \mathbb{R}^{d}$ we set $\overline{v}=(v_{1}, \ldots, v_{d-1})$ to be the projection onto the first $d-1$ coordinates, and we set $v^{z}=v_{d}$ to be the projection onto the last coordinate. For any $a<b$ define the following sets
\begin{align*}
&\Delta^{k}_{c} := \{ (r_{1}, \ldots, r_{k}) \in \mathbb{R}^{k}\, :\,  r_{j} \geq 0, j=1, \ldots, k, \, r_{1}+\ldots+r_{k}\leq 1\},\\
&\Delta_{*}^{k} := \{ (y_{1}, \ldots, y_{k}) \in \mathbb{R}^{k}\, :\, a\leq y_{1}\leq y_{2} \leq  \ldots\leq y_{k}\leq b\}. 
\end{align*}
Let $n\geq 1$. If $n=2\ell$ we define 
\begin{align*}
&U_{n} :   \Delta_{c}^{\ell} \times \Delta_{*}^{\ell} \ni (\lambda_{1}, \ldots, \lambda_{\ell}, x_{1}, \ldots, x_{\ell}) \mapsto \sum_{j=1}^{\ell} \lambda_{j} \gamma(x_{j})  + (1-\sum_{j=1}^{\ell}\lambda_{j}) \gamma(b);\\
&L_{n} :\Delta_{c}^{\ell} \times \Delta_{*}^{\ell} \ni  (\lambda_{1}, \ldots, \lambda_{\ell}, x_{1}, \ldots, x_{\ell}) \mapsto (1-\sum_{j=1}^{\ell}\lambda_{j})\gamma(a)+\sum_{j=1}^{\ell} \lambda_{j} \gamma(x_{j}),
\end{align*}
and if $n=2\ell-1$ we define 
\begin{align*}
   &U_{n} :\Delta_{c}^{\ell} \times \Delta_{*}^{\ell-1} \ni  (\beta_{1}, \ldots, \beta_{\ell}, x_{2},\ldots, x_{\ell}) \mapsto (1-\sum_{j=1}^{\ell}\beta_{j})\gamma(a) +\sum_{j=2}^{\ell} \beta_{j} \gamma(x_{j})+\beta_{1} \gamma(b);\\
   &L_{n} :\Delta_{c}^{\ell-1} \times \Delta_{*}^{\ell} \ni (\beta_{2}, \ldots, \beta_{\ell}, x_{1},\ldots, x_{\ell}) \mapsto  (1-\sum_{j=2}^{\ell} \beta_{j})\gamma(x_{1})+\sum_{j=2}^{\ell} \beta_{j} \gamma(x_{j}).
\end{align*}
If $n=1$ we set $U_{1} : [0,1]=:\Delta_{c}^{1}\times \Delta_{*}^{0} \mapsto (1-\beta_{1})
\gamma(a)+\beta_{1}\gamma(b)$, and $L_{1} : [a,b]=:\Delta_{c}^{0}\times \Delta_{*}^{1} \mapsto \gamma(x_{1})$. 

Together with maps $U_{n}$ and $L_{n}$ we  define  functions $B^{\sup}$ (and $B^{\inf}$) on the image of $\overline{U}$ (or $\overline{L}$) such that 
\begin{align}
    &B^{\sup}(\overline{U}_{n})=U^{z}_{n}, \label{vog}\\
    &B^{\inf}(\overline{L}_{n}) = L^{z}_{n} \label{vyp}.
\end{align}
 We remark that at this moment $B^{\sup}$  (and $B^{\inf}$) is not {\em well defined}, i.e., it could be that there are points $s_{1}, s_{2}$, $s_{1} \neq s_{2}$ such that $\overline{U}_{n}(s_{1})=\overline{U}_{n}(s_{2})$ and at the same time $U^{z}_{n}(s_{1})\neq U^{z}_{n}(s_{2})$. However, we will see that the next theorem, in particular, claims that  both functions $B^{\sup}, B^{\inf}$ are well defined. 
\begin{theorem}\label{mth010}
Let $\gamma : [a,b] \mapsto \mathbb{R}^{n+1}$  be in $C([a,b])\cap C^{n+1}((a,b))$ with totally positive torsion. 

If $n =2\ell$, $\ell \geq 1$,  we have 
\begin{align}
   & \overline{U}_{2\ell}(\partial\,  (\Delta_{c}^{\ell} \times \Delta_{*}^{\ell})) =\overline{L}_{2\ell}(\partial\,  (\Delta_{c}^{\ell} \times \Delta_{*}^{\ell}))= \partial\,  \mathrm{conv}(\overline{\gamma}([a,b])), \label{b2l}\\
   &\overline{U}_{2\ell} : \mathrm{int} (\Delta_{c}^{\ell} \times \Delta_{*}^{\ell}) \mapsto \mathrm{int}(\mathrm{conv}(\overline{\gamma}([a,b]))) \quad \text{is diffeomorphism}, \label{diff2lu}\\
   &\overline{L}_{2\ell} : \mathrm{int} (\Delta_{c}^{\ell} \times \Delta_{*}^{\ell}) \mapsto \mathrm{int}(\mathrm{conv}(\overline{\gamma}([a,b]))) \quad \text{is diffeomorphism}. \label{diff2ll}
\end{align}
If $n=2\ell-1$ we have 
\begin{align}
   & \overline{U}_{2\ell-1}(\partial\,  (\Delta_{c}^{\ell} \times \Delta_{*}^{\ell-1})) =\overline{L}_{2\ell-1}(\partial\,  (\Delta_{c}^{\ell-1} \times \Delta_{*}^{\ell}))= \partial\,  \mathrm{conv}(\overline{\gamma}([a,b])), \label{b2l-1}\\
   &\overline{U}_{2\ell-1} : \mathrm{int} (\Delta_{c}^{\ell} \times \Delta_{*}^{\ell-1}) \mapsto \mathrm{int}(\mathrm{conv}(\overline{\gamma}([a,b]))) \quad \text{is diffeomorphism}, \label{diff2l-1u}\\
   &\overline{L}_{2\ell-1} : \mathrm{int} (\Delta_{c}^{\ell-1} \times \Delta_{*}^{\ell}) \mapsto \mathrm{int}(\mathrm{conv}(\overline{\gamma}([a,b]))) \quad \text{is diffeomorphism}. \label{diff2l-1l}
\end{align}
For all $n\geq 1$, 
\begin{align}\label{welld}
B^{\sup}, B^{\inf} \quad \text{are well defined}, \quad B^{\sup}, B^{\inf} \in C(\mathrm{conv}(\overline{\gamma}([a,b]))) \cap C^{1}(\mathrm{int}(\mathrm{conv}(\overline{\gamma}([a,b])))).
\end{align}
Next, for all $n\geq 1$ we have\footnote{When $n=1$ the equality $B^{\sup}(\overline{\gamma})=\gamma_{2}$ should be replaced by $B^{\sup}(\overline{\gamma})\geq \gamma_{2}.$}
\begin{align}
    &B^{\sup} \quad \text{is minimal concave  on} \quad \mathrm{conv}(\overline{\gamma}([a,b])) \quad \text{with} \quad   \, B^{\sup}(\overline{\gamma})=\gamma_{n+1}; \label{mincon1}\\
    &B^{\inf} \quad \text{is maximal convex on} \quad \mathrm{conv}(\overline{\gamma}([a,b])) \quad \text{with}\quad   \, B^{\inf}(\overline{\gamma})=\gamma_{n+1}; \label{maxcon2}
\end{align}
Moreover, 
\begin{align}
&B^{\inf}(y)=B^{\sup}(y) \quad \text{if and only if} \quad y \in \partial\,  \mathrm{conv}(\overline{\gamma}([a,b])), \label{giff}\\
& \partial\,  \mathrm{conv}(\gamma([a,b]))=\{(x,B^{\mathrm{sup}}(x)), x \in \mathrm{conv}(\bar{\gamma}([a,b]))\} \cup \{(x,B^{\mathrm{inf}}(x)), x \in \mathrm{conv}(\bar{\gamma}([a,b]))\}. \label{union}
\end{align}
\end{theorem}

The statement of  Theorem~\ref{mth010} may seem a bit technical, however, we think that the  intuition behind the construction of the convex hulls is natural. We refer the reader to schematic pictures in Fig.~\ref{fig:sketches} for better understanding of the claims made in the theorem. In Fig.~\ref{fig:4d} the domain $\mathrm{conv}(\overline{\gamma}([a,b]))$  of  $B^{\sup}$ in $\mathbb{R}^{3}$ is foliated by triangles where $B^{\sup}$ is linear on each such triangle. 

\newcommand{\showsketch}[1]{%
\begin{minipage}{.33\textwidth}%
\begin{center}
\includegraphics[trim = 1 1 1 1 , clip, width=\textwidth]{d=#1.pdf}
$n+1 = #1$
\end{center}
\end{minipage}%
}

\begin{figure}[t]
    \centering
    \showsketch{1}\showsketch{3}\showsketch{5}
    \showsketch{2}\showsketch{4}\showsketch{6}
    \caption{These schematic pictures  clarify how the convex hull of the space $\gamma$ with totally positive torsion is parametrized. If $n$ is even then  the {\em upper hull} is described by convex combination of $\frac{n}{2}+1$ points  of $\gamma$, where among these points, $\frac{n}{2}$ are {\em free}, i.e., they are chosen in an arbitrary way on the space curve, and the last point $\gamma(b)$ is always fixed. For the {\em lower hull} $\gamma(a)$ is fixed instead of $\gamma(b)$. If $n$ is odd the picture is asymmetric. In this case the {\em upper hull} fixes $2$ endpoints $\gamma(a)$ and $\gamma(b)$ and has $\frac{n-1}{2}$ free points. The lower hull has $\frac{n+1}{2}$ free points, and no fixed points. The case $n=0$ (the convex hull of an interval), not mentioned in Theorem~\ref{mth010}, was helpful to guess the construction in higher dimensions, it has two fixed points $\gamma(a)$ and $\gamma(b)$. Compare with the exact pictures for the cases $n+1= 2,3,4$ shown in Figures \ref{fig:2d},\ref{fig:3d} and \ref{fig:4d}. }
    \label{fig:sketches}
\end{figure}

\begin{figure}[t]
    \centering
    \includegraphics[width=.49\textwidth]{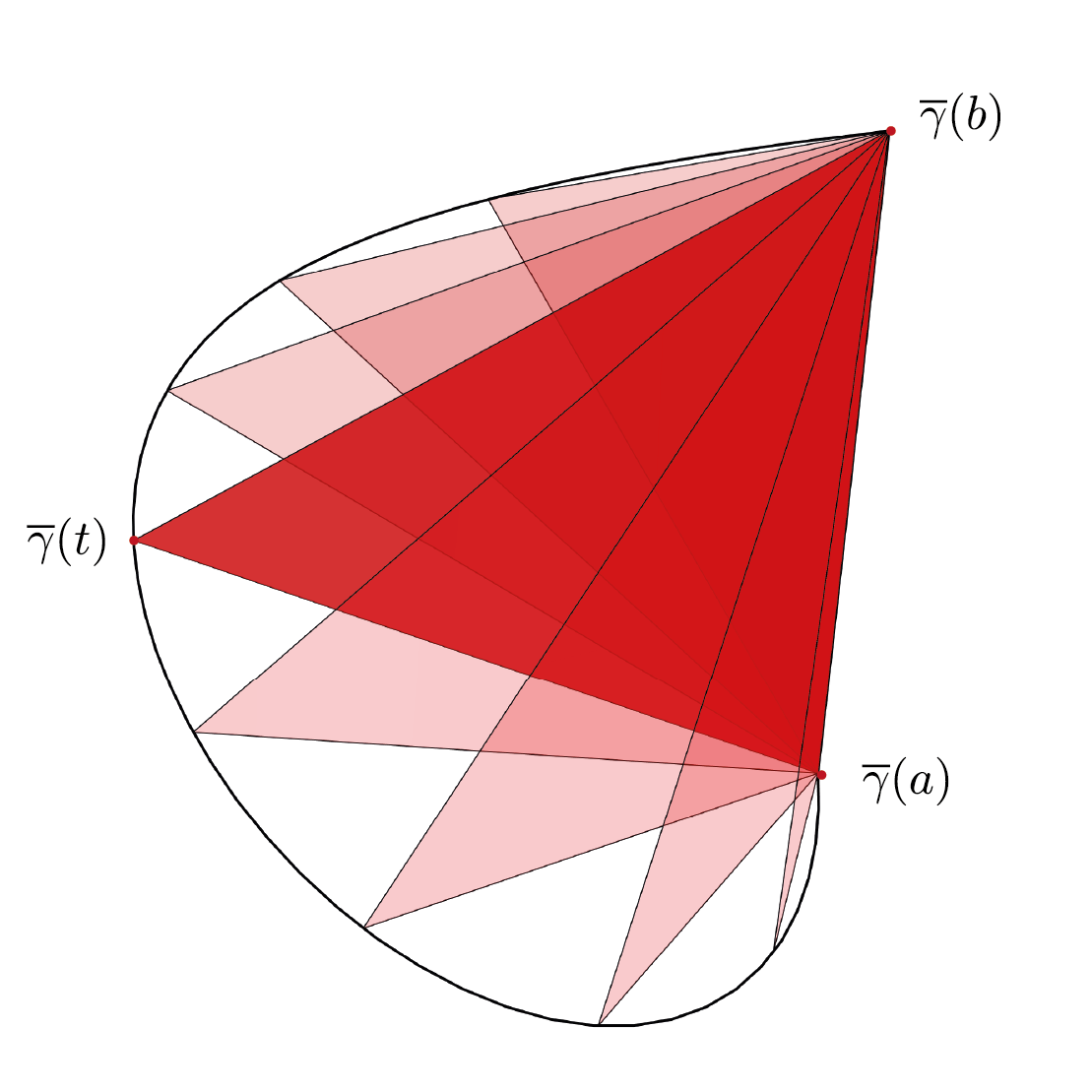}
    \caption{For $n+1=3+1$ the set $\mathrm{conv}(\overline{\gamma}([a,b]))$ is foliated by triangles (simplices) with vertices $\overline{\gamma}(a), \overline{\gamma}(b)$ and $\overline{\gamma}(t)$ for each $t \in (a,b)$.  The function $B^{\sup}$ is linear on each such triangle and $B^{\sup}(\overline{\gamma})=\gamma_{4}$. Also $B^{\sup}=B^{\inf}$ on edges of each triangle.}
    \label{fig:4d}
\end{figure}

\vskip0.3cm 

Perhaps it may seem that the total positivity of the torsion, i.e., the fact that the leading principal minors of $(\gamma', \ldots, \gamma^{(n+1)})$ have positive signs on $(a,b)$, is a redundant assumption for Theorem~\ref{mth010} to hold true. However, the next proposition shows that the total positivity is a sensitive assumption.
\begin{proposition}\label{sensitive}
There exists a curve $\gamma : [-1,1] \to \mathbb{R}^{2+1}$ in $C^{\infty}([-1,1])$ such that the leading principal minors of $(\gamma', \gamma'', \gamma''')$ are positive on $[-1,1]$ except the $2\times2$ and $3\times 3$ principal minors vanish at $t=0$, and the map $B^{\sup}$ defined by (\ref{vog}) is not concave on $\mathrm{conv}(\overline{\gamma}([-1,1]))$. 
\end{proposition}

The next theorem answers Question~\ref{que1}, and also provides us with optimizers, i.e., the random variables $Y$ which attain supremum  (infimum) in  Question~\ref{que1}. 

\begin{theorem}\label{mth1} Let $\gamma : [a,b] \to \mathbb{R}^{n+1}$,  $\gamma \in C([a,b]) \cap C^{n+1}((a,b))$ be such that all the leading principal minors of the  $(n+1)\times (n+1)$ matrix $(\gamma'(t), \ldots, \gamma^{(n+1)}(t))$  are positive for all $t \in (a,b)$. Then 
\begin{align}
\sup_{a\leq Y\leq b} \{ \mathbb{E}\gamma_{n+1}(Y) \, :\, \mathbb{E}\bar{\gamma}(Y)=x\} &=B^{\mathrm{sup}}(x),\label{extr01}\\
\inf_{a\leq Y\leq b} \{ \mathbb{E}\gamma_{n+1}(Y) \, :\, \mathbb{E}\bar{\gamma}(Y)=x\} &=B^{\mathrm{inf}}(x),\label{extr02}
\end{align}
hold for all $x \in \mathrm{conv}(\overline{\gamma}([a,b]))$, where $B^{\sup}$ and $B^{\inf}$ are given by (\ref{vog}) and (\ref{vyp}). Moreover, given $x \in \mathrm{conv}(\overline{\gamma}([a,b]))$ the supremum in (\ref{extr01}) (infimum in (\ref{extr02})) is attained  by the random variable $\zeta(x)$ (the random variable $\xi(x)$) defined as follows: 

Case 1: $n=2\ell-1$. Then by (\ref{b2l-1}) and (\ref{diff2l-1u}),  $x =(1-\sum_{j=1}^{\ell}\beta_{j})\overline{\gamma}(a)+\sum_{j=2}^{\ell}\beta_{j} \overline{\gamma}(x_{j})+\beta_{1}\overline{\gamma}(b)$ for some $(\beta_{1}, \ldots, \beta_{\ell}, x_{2}, \ldots, x_{\ell})\in \Delta_{c}^{\ell}\times\Delta_{*}^{\ell-1}$. Set $\mathbb{P}(\zeta(x)=a)=1-\sum_{j=1}^{\ell}\beta_{j}$, $\mathbb{P}(\zeta(x)=b)=\beta_{1}$, and $\mathbb{P}(\zeta(x)=x_{j})=\beta_{j}$ for $j=2,\ldots, \ell$. Also, by  (\ref{b2l-1}) and (\ref{diff2l-1l}),  $x =(1-\sum_{j=2}^{\ell}\lambda_{j})\overline{\gamma}(y_{1})+\sum_{j=2}^{\ell}\lambda_{j} \overline{\gamma}(y_{j})$ for some $(\lambda_{2}, \ldots, \lambda_{\ell}, y_{1}, \ldots, y_{\ell})\in \Delta_{c}^{\ell-1}\times\Delta_{*}^{\ell}$. Set $\mathbb{P}(\xi(x)=y_{1})=1-\sum_{j=2}^{\ell}\lambda_{j}$,  and $\mathbb{P}(\xi(x)=y_{j})=\lambda_{j}$ for $j=2,\ldots, \ell$.

Case 2: $n=2\ell$. Then by (\ref{b2l}) and (\ref{diff2lu}),  $x =\sum_{j=1}^{\ell}\beta_{j} \overline{\gamma}(x_{j})+(1-\sum_{j=1}^{\ell}\beta_{j})\overline{\gamma}(b)$ for some $(\beta_{1}, \ldots, \beta_{\ell}, x_{1}, \ldots, x_{\ell})\in \Delta_{c}^{\ell}\times\Delta_{*}^{\ell}$. Set  $\mathbb{P}(\zeta(x)=b)=1-\sum_{j=1}^{\ell}\beta_{j}$, and $\mathbb{P}(\zeta(x)=x_{j})=\beta_{j}$ for $j=1,\ldots, \ell$. Also, by  (\ref{b2l}) and (\ref{diff2ll}),  $x =(1-\sum_{j=1}^{\ell}\lambda_{j})\overline{\gamma}(a)+\sum_{j=1}^{\ell}\lambda_{j} \overline{\gamma}(y_{j})$ for some $(\lambda_{1}, \ldots, \lambda_{\ell}, y_{1}, \ldots, y_{\ell})\in \Delta_{c}^{\ell}\times\Delta_{*}^{\ell}$. Set $\mathbb{P}(\xi(x)=a)=1-\sum_{j=1}^{\ell}\lambda_{j}$,  and $\mathbb{P}(\xi(x)=y_{j})=\lambda_{j}$ for $j=1,\ldots, \ell$.
\end{theorem}

The next corollary recovers the result of Karlin--Sharpley \cite{Karlin2}, i.e., the equations  (\ref{nobel}) in case of the moment curve.
\begin{corollary}\label{nobel2}
Let $\gamma(t) = (t, \ldots, t^{n}, t^{n+1}) : [0,1] \to \mathbb{R}^{n+1}$. If $x=(x_{1}, \ldots, x_{n}) \in \mathrm{int}(\mathrm{conv} (\overline{\gamma}([0,1])))$ then $B^{\sup}(x)$ and $B^{\inf}(x)$ are the unique solutions $x_{n+1}$ of the equations $K_{n+1}=0$ and $S_{n+1}=0$ correspondingly, where $K_{n+1}$ and $S_{n+1}$ are defined by (\ref{Sharp1}) and (\ref{kar1}).
\end{corollary}

In the next corollary we give a sufficient local description of convex curves. Recall that a curve $\gamma : [a,b] \to \mathbb{R}^{n}$ is called {\em convex} if no $n+1$ its different points lie in a single affine hyperplane. 

\begin{corollary}\label{karatecor}
Let $\gamma : [a,b] \to \mathbb{R}^{n}$,  $\gamma \in C([a,b]) \cap C^{n}((a,b))$ be such that all the leading principal minors of the  $n\times n$ matrix $(\gamma'(t), \ldots, \gamma^{(n)}(t))$  are positive for all $t \in (a,b)$. Then $\gamma$ is convex. In particular, for any integer $k$, $1\leq k \leq n$,  the equation $c_{0}+c_{1}\gamma_{1}(t)+\ldots+c_{k}\gamma_{k}(t)=0$ has at most $k$ roots on $[a,b]$ provided that $(c_{0}, \ldots, c_{k})\neq (0, \ldots, 0)$. 
\end{corollary}

Recall the definition of Carath\'eodory number $b(\gamma)$ of a curve $\gamma : [a,b] \to \mathbb{R}^{n}$, i.e., the smallest integer $k$ such that any point of $\mathrm{conv}(\gamma([a,b]))$ can be represented as convex combination of at most $k$ points of $\gamma([a,b])$, see (\ref{karate1}). The next corollary directly follows  from Theorem~\ref{mth010} (parts (\ref{b2l-1}), (\ref{diff2l-1l}), (\ref{b2l}), and  (\ref{diff2ll})).
\begin{corollary}\label{karatekid}
    Let $\gamma : [a,b] \to \mathbb{R}^{n}$,  $\gamma \in C([a,b]) \cap C^{n}((a,b))$ be a curve with totally positive torsion. Then its Carath\'eodory number equals to $\lfloor \frac{n+2}{2}\rfloor$.  
\end{corollary}

In the next corollary we obtain formulas for the volumes of the convex hulls of a space curve having totally positive torsion both in even and odd dimensions. 
\begin{corollary}\label{provolume}
Let $\gamma : [a,b] \to \mathbb{R}^{n}$,  $\gamma \in C([a,b]) \cap C^{n}((a,b))$ be a curve with totally positive torsion. If $n=2 \ell$ then 
\begin{align*}
    \mathrm{Vol}(\mathrm{conv}(\gamma([a,b])))& \\
    =\frac{(-1)^{\frac{\ell(\ell-1)}{2}}}{(2\ell)!}& \int_{a\leq x_{1}\leq \ldots \leq x_{\ell} \leq b} \mathrm{det}(\gamma(x_{1})-\gamma(a), \ldots, \gamma(x_{\ell})-\gamma(a), \gamma'(x_{1}), \ldots, \gamma'(x_{\ell})) \, dx \\
    =\frac{(-1)^{\frac{\ell(\ell-1)}{2}}}{(2\ell)!}& \int_{a\leq x_{1}\leq \ldots \leq x_{\ell} \leq b} \mathrm{det}(\gamma(x_{1})-\gamma(b), \ldots, \gamma(x_{\ell})-\gamma(b), \gamma'(x_{1}), \ldots, \gamma'(x_{\ell})) \, dx. 
\end{align*}
If $n=2\ell-1$ then 
\begin{align*}
    &\mathrm{Vol}(\mathrm{conv}(\gamma([a,b]))) \\
    &=\frac{(-1)^{\frac{(\ell-1)(\ell-2)}{2}}}{(2\ell-1)!} \int_{a\leq x_{2}\leq \ldots \leq x_{\ell} \leq b} \mathrm{det}(\gamma(b)-\gamma(a), \gamma(x_{2})-\gamma(a), \ldots, \gamma(x_{\ell})-\gamma(a), \gamma'(x_{2}), \ldots, \gamma'(x_{\ell})) \, dx \\
    &=\frac{(-1)^{\frac{\ell(\ell-1)}{2}}}{(2\ell-1)!} \int_{a\leq x_{1}\leq \ldots \leq x_{\ell} \leq b} \mathrm{det}(\gamma(x_{2})-\gamma(x_{1}), \ldots, \gamma(x_{\ell})-\gamma(x_{1}), \gamma'(x_{1}), \ldots, \gamma'(x_{\ell}))\,  dx. 
\end{align*}
\end{corollary}

Let $\mathrm{Area}$ denote $n$ dimensional Lebesgue measure in $\mathbb{R}^{n+1}$, and let $A^{\mathrm{Tr}}$ be the transpose of a matrix $A$.  
\begin{corollary}\label{area1}
Let $\gamma : [a,b] \to \mathbb{R}^{n+1}$,  $\gamma \in C^{1}([a,b]) \cap C^{n+1}((a,b))$ be a curve with totally positive torsion. If $n=2 \ell$ then 
\begin{align*}
\mathrm{Area}(\partial \; \mathrm{conv}(\gamma([a,b]))) = \frac{1}{n!} \int_{a\leq x_{1}\leq \ldots \leq x_{\ell}\leq b} \left( \sqrt{\det S_{a}^{\mathrm{Tr}}S_{a}} +\sqrt{\det S_{b}^{\mathrm{Tr}}S_{b}} \right) dx, 
\end{align*}
where $S_{r} = (\gamma(x_{1})-\gamma(r), \ldots, \gamma(x_{\ell})-\gamma(r), \gamma'(x_{1}), \ldots, \gamma'(x_{\ell}))$ is $(2\ell+1)\times 2\ell$ matrix, and $dx$ is $\ell$ dimensional Lebesgue measure. 

If $n=2\ell-1$ then 
\begin{align*}
\mathrm{Area}(\partial \; \mathrm{conv}(\gamma([a,b]))) = \frac{1}{n!} \int_{a\leq x_{2}\leq \ldots \leq x_{\ell}\leq b}  \sqrt{\det \Psi^{\mathrm{Tr}}\Psi} \, d\tilde{x} +\frac{1}{n!} \int_{a\leq x_{1}\leq \ldots \leq x_{\ell}\leq b}  \sqrt{\det \Phi^{\mathrm{Tr}}\Phi} \, dx, 
\end{align*}
where $\Psi = (\gamma(b)-\gamma(a), \gamma(x_{2})-\gamma(a), \ldots, \gamma(x_{\ell})-\gamma(a), \gamma'(x_{2}), \ldots, \gamma'(x_{\ell}))$,  $\Phi = (\gamma(x_{2})-\gamma(x_{1}), \ldots, \gamma(x_{\ell})-\gamma(x_{1}), \gamma'(x_{1}), \ldots, \gamma'(x_{\ell}))$ are $2\ell \times (2\ell-1)$ size matrices, and $d\tilde{x}$ denotes $\ell-1$ dimensional Lebesgue measure. 
\end{corollary}

\section{The proof of main results}\label{damtkiceba}
Sometimes we will omit the index $n$ and simply write $U, L$ instead of $U_{n}, L_{n}$, and it will be clear from the context what is the corresponding number $n$. Before we start proving Theorem~\ref{mth010}, first let us state several lemmas that will be helpful throughout the rest of the paper.  The next lemma illustrates {\em local to global} principle. 
\begin{lemma}\label{klasika}
If the torsion of $\gamma$ is totally positive on $(a,b)$ then 
\begin{align}\label{dplane}
\det(\gamma'(x_{1}), \gamma'(x_{2}), \ldots \gamma'(x_{n+1}))>0
\end{align}
for all $a<x_{1}<\ldots<\ldots<x_{n+1}<b$.
\end{lemma}

\begin{proof}
Without loss of generality assume $[a,b]=[0,1]$. The lemma can be derived from the identity (9) obtained in \cite{DW}. As the lemma is an important step in the proofs of the main results stated in this paper, for the readers convenience we decided to include the  proof of the lemma without invoking the identity from \cite{DW}. 

We have 
\begin{align}
   &\det  \begin{pmatrix}
    \gamma'_{1}(x_{1}) & \gamma'_{1}(x_{2}) & \dots & \gamma'_{1}(x_{n+1})\\
    \gamma'_{2}(x_{1}) & \gamma'_{2}(x_{2}) & \ldots & \gamma'_{2}(x_{n+1})\\
    \vdots & \vdots &\ddots & \vdots \\
    \gamma'_{n+1}(x_{1}) & \gamma'_{n+1}(x_{2}) & \dots & \gamma'_{n+1}(x_{n+1})
    \end{pmatrix} = \nonumber\\
    &\det \begin{pmatrix}
    1 & 1& \dots & 1\\
    \frac{\gamma'_{2}(x_{1})}{\gamma'_{1}(x_{1})} & \frac{\gamma'_{2}(x_{2})}{\gamma'_{1}(x_{2})} & \ldots & \frac{\gamma'_{2}(x_{n+1})}{\gamma'_{1}(x_{n+1})}\\
    \vdots & \vdots &\ddots & \vdots \\
    \frac{\gamma'_{n+1}(x_{1})}{\gamma'_{1}(x_{1})} & \frac{\gamma'_{n+1}(x_{2})}{\gamma'_{1}(x_{2})} & \dots & \frac{\gamma'_{n+1}(x_{n+1})}{\gamma'_{1}(x_{n+1})}
    \end{pmatrix}\, \prod_{j=1}^{n+1} \gamma'_{1}(x_{j})   = \nonumber\\
     &\det \begin{pmatrix}
    1 & 0& \dots & 0\\
    \frac{\gamma'_{2}(x_{1})}{\gamma'_{1}(x_{1})} & \frac{\gamma'_{2}(x_{2})}{\gamma'_{1}(x_{2})}-  \frac{\gamma'_{2}(x_{1})}{\gamma'_{1}(x_{1})}  & \ldots & \frac{\gamma'_{2}(x_{n+1})}{\gamma'_{1}(x_{n+1})}-  \frac{\gamma'_{2}(x_{1})}{\gamma'_{1}(x_{1})} \\
    \vdots & \vdots &\ddots & \vdots \\
    \frac{\gamma'_{n+1}(x_{1})}{\gamma'_{1}(x_{1})} & \frac{\gamma'_{n+1}(x_{2})}{\gamma'_{1}(x_{2})} -\frac{\gamma'_{n+1}(x_{1})}{\gamma'_{1}(x_{1})}  & \dots & \frac{\gamma'_{n+1}(x_{n+1})}{\gamma'_{1}(x_{n+1})} -\frac{\gamma'_{n+1}(x_{1})}{\gamma'_{1}(x_{1})} 
    \end{pmatrix}\, \prod_{j=1}^{n+1} \gamma'_{1}(x_{j})  = \nonumber\\
    &\det \begin{pmatrix}
    \frac{\gamma'_{2}(x_{2})}{\gamma'_{1}(x_{2})}-  \frac{\gamma'_{2}(x_{1})}{\gamma'_{1}(x_{1})}  & \ldots & \frac{\gamma'_{2}(x_{n+1})}{\gamma'_{1}(x_{n+1})}-  \frac{\gamma'_{2}(x_{1})}{\gamma'_{1}(x_{1})} \\
     \vdots &\ddots & \vdots \\
     \frac{\gamma'_{n+1}(x_{2})}{\gamma'_{1}(x_{2})} -\frac{\gamma'_{n+1}(x_{1})}{\gamma'_{1}(x_{1})}  & \dots & \frac{\gamma'_{n+1}(x_{n+1})}{\gamma'_{1}(x_{n+1})} -\frac{\gamma'_{n+1}(x_{1})}{\gamma'_{1}(x_{1})} 
    \end{pmatrix}\, \prod_{j=1}^{n+1} \gamma'_{1}(x_{j})  \stackrel{(*)}{=} \nonumber\\
    &\det \begin{pmatrix}
    \frac{\gamma'_{2}(x_{2})}{\gamma'_{1}(x_{2})}-  \frac{\gamma'_{2}(x_{1})}{\gamma'_{1}(x_{1})}  & \ldots & \frac{\gamma'_{2}(x_{n+1})}{\gamma'_{1}(x_{n+1})}-  \frac{\gamma'_{2}(x_{n})}{\gamma'_{1}(x_{n})} \\
     \vdots &\ddots & \vdots \\
     \frac{\gamma'_{n+1}(x_{2})}{\gamma'_{1}(x_{2})} -\frac{\gamma'_{n+1}(x_{1})}{\gamma'_{1}(x_{1})}  & \dots & \frac{\gamma'_{n+1}(x_{n+1})}{\gamma'_{1}(x_{n+1})} -\frac{\gamma'_{n+1}(x_{n})}{\gamma'_{1}(x_{n})} 
    \end{pmatrix}\, \prod_{j=1}^{n+1} \gamma'_{1}(x_{j})  = \nonumber\\
    &\int_{x_{1}}^{x_{2}} \int_{x_{2}}^{x_{3}} \dots \int_{x_{n}}^{x_{n+1}} \det 
    \begin{pmatrix}
    \left(\frac{\gamma'_{2}(y_{1})}{\gamma'_{1}(y_{1})}\right)'  & \ldots & \left(\frac{\gamma'_{2}(y_{n})}{\gamma'_{1}(y_{n})}\right)' \\
     \vdots &\ddots & \vdots \\
    \left(\frac{\gamma'_{n+1}(y_{1})}{\gamma'_{1}(y_{1})}\right)' & \dots & \left(\frac{\gamma'_{n+1}(y_{n})}{\gamma'_{1}(y_{n})}\right)' 
    \end{pmatrix} dy_{1} dy_{2}\dots dy_{n}\, \prod_{j=1}^{n+1} \gamma'_{1}(x_{j}), \nonumber
\end{align}
where in the equality $(*)$ we used the property of the determinant that if $v_{1}, \ldots, v_{k}$ are column vectors in $\mathbb{R}^{k}$ then $\det(v_{2}-v_{1}, v_{3}-v_{1}, \ldots, v_{k}-v_{1})) = \det(v_{2}-v_{1}, v_{3}-v_{2}, \ldots, v_{k}-v_{k-1})$ by subtracting the columns from each other.

The leading principal minors of the matrix $(\gamma', \gamma'', \ldots, \gamma^{(n+1)})$ are positive. In particular $\gamma'_{1}$ is positive on $(0,1)$, and hence the factor $\prod_{j=1}^{n+1}\gamma'(x_{j})>0$. To verify (\ref{dplane}) it suffices to show 
\begin{align}\label{dplane2}
\det 
    \begin{pmatrix}
    \left(\frac{\gamma'_{2}(y_{1})}{\gamma'_{1}(y_{1})}\right)'  & \ldots & \left(\frac{\gamma'_{2}(y_{n})}{\gamma'_{1}(y_{n})}\right)' \\
     \vdots &\ddots & \vdots \\
    \left(\frac{\gamma'_{n+1}(y_{1})}{\gamma'_{1}(y_{1})}\right)' & \dots & \left(\frac{\gamma'_{n+1}(y_{n})}{\gamma'_{1}(y_{n})}\right)' 
    \end{pmatrix} >0 \quad \text{for all} \quad 0<y_{1}<y_{2}<\ldots<y_{n}<1. 
\end{align}
We will repeat the same computation as before but now for the determinant in (\ref{dplane2}), and, eventually, we will see that the proof of the lemma will be  just $n$  times the application of the previous computation together with an identity for determinants that we have not described yet. 

Before we proceed let us make couple of observations. We started with  the determinant of $(n+1)\times (n+1)$ matrix.  Next, we divided the columns by the entries in the first row which consist of $\gamma'_{1}>0$, and after the Gaussian elimination and the fundamental theorem of calculus we ended up with the integral of the determinant of $n \times n$, and we also acquired the factor $\prod_{j=1}^{n+1} \gamma'_{1}(x_{j}) >0$. To repeat the same computation for the determinant in (\ref{dplane2}) and the ones that we obtain in a similar manner we should verify that the entries in the first row of all such new matrices (of smaller sizes) are positive. Such entries are changed as follows 
\begin{align}\label{iteracia}
\gamma'_{1} \stackrel{\mathrm{step \,1}}{\mapsto} \left(\frac{\gamma'_{2}}{\gamma'_{1}}\right)' \stackrel{\mathrm{step\, 2}}{\mapsto} \left( \frac{\left(\frac{\gamma'_{3}}{\gamma'_{1}}\right)'}{\left(\frac{\gamma'_{2}}{\gamma'_{1}}\right)'} \right)' \stackrel{\mathrm{step\, 3}}{\mapsto}  \left(\frac{\left( \frac{\left(\frac{\gamma'_{4}}{\gamma'_{1}}\right)'}{\left(\frac{\gamma'_{2}}{\gamma'_{1}}\right)'} \right)'}{\left( \frac{\left(\frac{\gamma'_{3}}{\gamma'_{1}}\right)'}{\left(\frac{\gamma'_{2}}{\gamma'_{1}}\right)'} \right)'}\right)' \stackrel{\mathrm{step\, 4}}{\mapsto} \ldots\, .
\end{align}
We claim that after $k$'th step, $1 \leq k \leq n$,  the obtained entry is of the form $\frac{\Delta_{k+1}  \Delta_{k-1}}{\Delta^{2}_{k}}$, 
where $\Delta_{\ell }$ denotes the leading $\ell \times \ell $ principal minor of the matrix $(\gamma',\gamma'', \ldots, \gamma^{(n+1)})$ (by definition we set $\Delta_{0}:=1$).  Assuming the claim,  Lemma~\ref{klasika} follows immediately  because of the condition $\Delta_{\ell}>0$ on $(0,1)$ for all $0\leq \ell \leq n+1$. 

To verify the claim we set $T = (\gamma', \gamma'', \ldots, \gamma^{(n+1)})$. Given subsets $I, J \subset \{1, \ldots, n+1\}$ we define $T_{I\times J}$ to be the determinant of the submatrix of $T$ formed by choosing the rows of the index set $I$ and the columns of index set $J$. 
We have 
\begin{align}
 \left(\frac{\gamma'_{2}}{\gamma'_{1}}\right)' &= \frac{\gamma''_{2}\gamma'_{1}-\gamma''_{1}\gamma'_{2}}{\gamma'_{1}}=\frac{T_{\{1,2\}\times\{1,2\}}}{T^{2}_{\{1\}\times \{1\}}}, \nonumber\\
 \left(\frac{\gamma'_{\ell}}{\gamma'_{1}}\right)' &= \frac{T_{\{1,\ell\}\times\{1,2\}}}{T^{2}_{\{1\}\times \{1\}}}, \quad \text{for all} \quad  \ell \geq 2; \nonumber\\
  \left( \frac{\left(\frac{\gamma'_{\ell}}{\gamma'_{1}}\right)'}{\left(\frac{\gamma'_{2}}{\gamma'_{1}}\right)'} \right)' &= \left(\frac{T_{\{1,\ell\}\times\{1,2\}}}{T_{\{1,2\}\times\{1,2\}}}\right)' \stackrel{(*)}{=} \frac{T_{\{1,\ell\}\times\{1,3\}} T_{\{1,2\}\times \{1,2\}} - T_{\{1,\ell\}\times\{1,2\}} T_{\{1,2\}\times \{1,3\}}}{T^{2}_{\{1,2\}\times\{1,2\}}} \nonumber\\
 &\stackrel{(**)}{=}\frac{T_{\{1,2,\ell\}\times\{1,2,3\}}\, T_{\{1\}\times\{1\}}}{T^{2}_{\{1,2\}\times\{1,2\}}}, \quad \text{for all} \quad \ell \geq 3, \label{ind2} 
\end{align}
where $(*)$ follows from the identity $(T_{I\times \{1,2,\ldots, k-1, k\}})'=T_{I\times \{1,2,\ldots, k-1, k+1\}}$, and $(**)$ follows from the following general identity for determinants:
\begin{align}\label{tozhd1}
T_{\{[k-2], \ell\}\times \{[k-2], k\}} T_{[k-1]\times [k-1]}-T_{\{[k-2], \ell\}\times [k-1]} T_{[k-1]\times\{[k-2], k\}}=T_{\{[k-1],\ell\}\times [k]} T_{[k-2]\times [k-2]}
\end{align}
 for all $k, 3 \leq k \leq n+1$, where we set  $[d]:=\{1,2, \ldots, d\}$ for a positive integer $d$.  Before we verify the identity  (\ref{tozhd1}), notice that it also implies 
 \begin{align}
     \left(\frac{T_{\{[k-2], \ell \}\times[k-1]}}{T_{[k-1]\times [k-1]}}\right)'&=\frac{T_{\{[k-2], \ell\}\times \{[k-2], k\}} T_{[k-1]\times [k-1]}-T_{\{[k-2], \ell\}\times [k-1]} T_{[k-1]\times\{[k-2], k\}}}{T^{2}_{[k-1]\times [k-1]}} \nonumber\\
     &= \frac{T_{\{[k-1],\ell\}\times [k]} T_{[k-2]\times [k-2]}}{T^{2}_{[k-1]\times [k-1]}}, \label{tozhd2}
 \end{align}
 for all $k, \ell$ such that $3\leq k \leq n+1$ and $k-1\leq \ell \leq n+1$. Therefore 
\begin{align*}
    \left(\frac{\left( \frac{\left(\frac{\gamma'_{\ell}}{\gamma'_{1}}\right)'}{\left(\frac{\gamma'_{2}}{\gamma'_{1}}\right)'} \right)'}{\left( \frac{\left(\frac{\gamma'_{3}}{\gamma'_{1}}\right)'}{\left(\frac{\gamma'_{2}}{\gamma'_{1}}\right)'} \right)'}\right)' \stackrel{(\ref{ind2})}{=} \left(\frac{T_{\{1,2,\ell\}\times\{1,2,3\}}}{T_{\{1,2,3\}\times\{1,2,3\}}}\right)' \stackrel{(\ref{tozhd2})}{=} \frac{T_{\{[3],\ell\}\times [4]} T_{[2]\times [2]}}{T^{2}_{[3]\times [3]}}.
\end{align*} 
In particular, after step 3,  the entry in (\ref{iteracia}) becomes $\frac{T_{[4]\times [4]} T_{[2]\times [2]}}{T^{2}_{[3]\times [3]}}>0$ because $T_{[k]\times[k]}=\Delta_{k}$. It then follows that after step $k$, the entry in (\ref{iteracia}) takes the form 
\begin{align*}
    \left(\frac{T_{\{[k-1],k+1\}\times[k]}}{T_{[k]\times[k]}}\right)' \stackrel{(\ref{tozhd1})}{=} \frac{T_{[k+1]\times[k+1]} T_{[k-1]\times [k-1]}}{T^{2}_{[k]\times[k]}} = \frac{\Delta_{k+1} \Delta_{k-1}}{\Delta_{k}} >0,
\end{align*}
for all $1\leq k \leq n$. Thus the proof of Lemma~\ref{klasika} is complete provided that the determinant identity (\ref{tozhd1}) is verified. Let $\Delta$ be an invertible $(k-2)\times (k-2)$ matrix,  $p,w,u,q \in \mathbb{R}^{k-2}$, and let $a,b,c,d \in \mathbb{R}$. To verify the identity (\ref{tozhd1}) it suffices to show that 
\begin{align}\label{sila}
    \det \begin{pmatrix}
    \Delta & q^{T} \\
    w & a
    \end{pmatrix} \det \begin{pmatrix}
    \Delta & u^{T} \\
    p & b
    \end{pmatrix} - \det \begin{pmatrix}
    \Delta & u^{T} \\
    w & c
    \end{pmatrix} \det \begin{pmatrix}
    \Delta & q^{T} \\
    p & d
    \end{pmatrix} = \det \begin{pmatrix}
    \Delta & u^{T} & q^{T} \\
    p & b & d \\
    w & c & a
    \end{pmatrix} \det \Delta.
\end{align}
Since $\det \begin{pmatrix}
A & B \\
C & D
\end{pmatrix} = \det A \, \det (D -CA^{-1}B)$ for an invertible $m\times m$ matrix $A$, and arbitrary $n\times n$ matrix $D$, $n\times m$ matrix $B$, and $m\times n$ matrix $C$, we see that (\ref{sila}) simplifies to 
\begin{align*}
    &(\det \Delta)^{2} \left[  (a-w\Delta^{-1} q^{T}) (b-p\Delta^{-1} u^{T})-(c-w \Delta^{-1} u^{T})(d-p\Delta^{-1} q^{T})\right]=\\
    &(\det \Delta)^{2}  \det \left( \begin{pmatrix}
    b & d \\
    c & a
    \end{pmatrix}- \begin{pmatrix}
   p\\
   w
    \end{pmatrix} \Delta^{-1} \begin{pmatrix}
    u^{T} q^{T}
    \end{pmatrix}\right),
\end{align*}
which holds  because 
$\begin{pmatrix}
   p\\
   w
    \end{pmatrix} \Delta^{-1} \begin{pmatrix}
    u^{T} q^{T}
    \end{pmatrix} = \begin{pmatrix}
    p \Delta^{-1}u^{T} & p\Delta^{-1} q^{T}\\
    w \Delta^{-1} u^{T} & w \Delta^{-1}q^{T}
    \end{pmatrix}$. The lemma is  proved. 
\end{proof}

%Next, we will need the following 
%\begin{corollary}
%If the torsion of $\gamma$ is totally positive on $(0,1)$, and $\gamma(0)=0$, then 
%\begin{align}\label{klasika2}
%    \det(\gamma(y_{1}), \ldots, \gamma(y_{n+1}))>0
%\end{align}
%for all $0<y_{1}<y_{2}<\ldots <y_{n}\leq 1$. 
%\end{corollary}
%\begin{proof}
%By Lemma~\ref{klasika} we have 
%\begin{align}\label{klas1}
%    \det(\gamma'(x_{1}), \ldots, \gamma'(x_{n+1})) >0
%\end{align}
%for all $0<x_{1}<\ldots <x_{n+1}<1$. Next, we integrate (\ref{klas1}) with respect to $x_{1}$ on the interval $(0,y_{1})$ we obtain 
%\begin{align}\label{klas2}
%   \det(\gamma(y_{1}), \gamma'(x_{2}), \ldots, \gamma'(x_{n+1}))  
%\end{align}
%holds for all $(x_{2}, \ldots, x_{n+1})$ such that $y_{1}<x_{2}<x_{3}<\ldots<x_{n+1}$. We integrate (\ref{klas2}) with respect to $x_{2}$ on the interval $(y_{1}, y_{2})$ etc. After iterating this procedure we obtain the corollary. 
%\end{proof}

\begin{corollary}\label{klasikac}
Let $a<b$, and let $\beta : [a,b] \to \mathbb{R}^{m}$ be a curve $\beta \in C([a,b])\cap C^{m}((a,b))$ with totally positive torsion. Choose any $a\leq z_{1}<\ldots<z_{m}\leq b$ and $r \in [0,1]\setminus\{z_{1}, \ldots, z_{m}\}$. Then the vectors $\beta(z_{1})-\beta(r), \ldots, \beta(z_{m})-\beta(r)$ are linearly independent in $\mathbb{R}^{m}$. 
\end{corollary}
\begin{proof}
Let $\nu$, $0\leq \nu \leq m$, be chosen in such a way that $r \in [z_{\nu}, z_{\nu+1}]$. Here we set $z_{0}:=a$, and $z_{m+1}:=b$.  We have 
\begin{align*}
    &\det(\beta(z_{1})-\beta(r), \ldots, \beta(z_{m})-\beta(r)) =\\
    &\pm \det(\beta(z_{2})-\beta(z_{1}),\ldots, \beta(r)-\beta(z_{\nu}), \beta(z_{\nu+1})-\beta(r), \ldots, \beta(z_{m})-\beta(z_{m-1}))=\\
    & \pm \int_{z_{m-1}}^{z_{m}}\ldots  \int_{r}^{z_{\nu+1}}\int_{z_{\nu}}^{r}\ldots \int_{z_{1}}^{z_{2}}\det(\beta'(s_{1}), \ldots,\beta'(s_{\nu}), \beta'(s_{\nu+1}),\ldots \beta'(s_{m}))ds_{1}\ldots ds_{m}\neq 0
\end{align*}
by Lemma~\ref{klasika}. 
\end{proof}

Certain parts of the proof of Theorem~\ref{mth010} will require induction on the dimension $n+1$. In particular, we will need to verify the base cases when $n=1$ (the odd case) and $n=2$ (even case). 
%As $n=1$ is trivial it suffices to consider only the case $2+1$. However, to better provide the intuition in high dimensions  we decided to include the complete proof of the theorem in dimensions $2+1$ and $3+1$ separately.  

In what follows without loss of generality  we assume $[a,b]=[0,1]$, and $\gamma(0)=0$. 
\subsection{The proof of Theorem~\ref{mth010} in dimension 1+1}\label{kkl1}
\begin{figure}[t]
\begin{center}
\includegraphics[width=.6\textwidth]{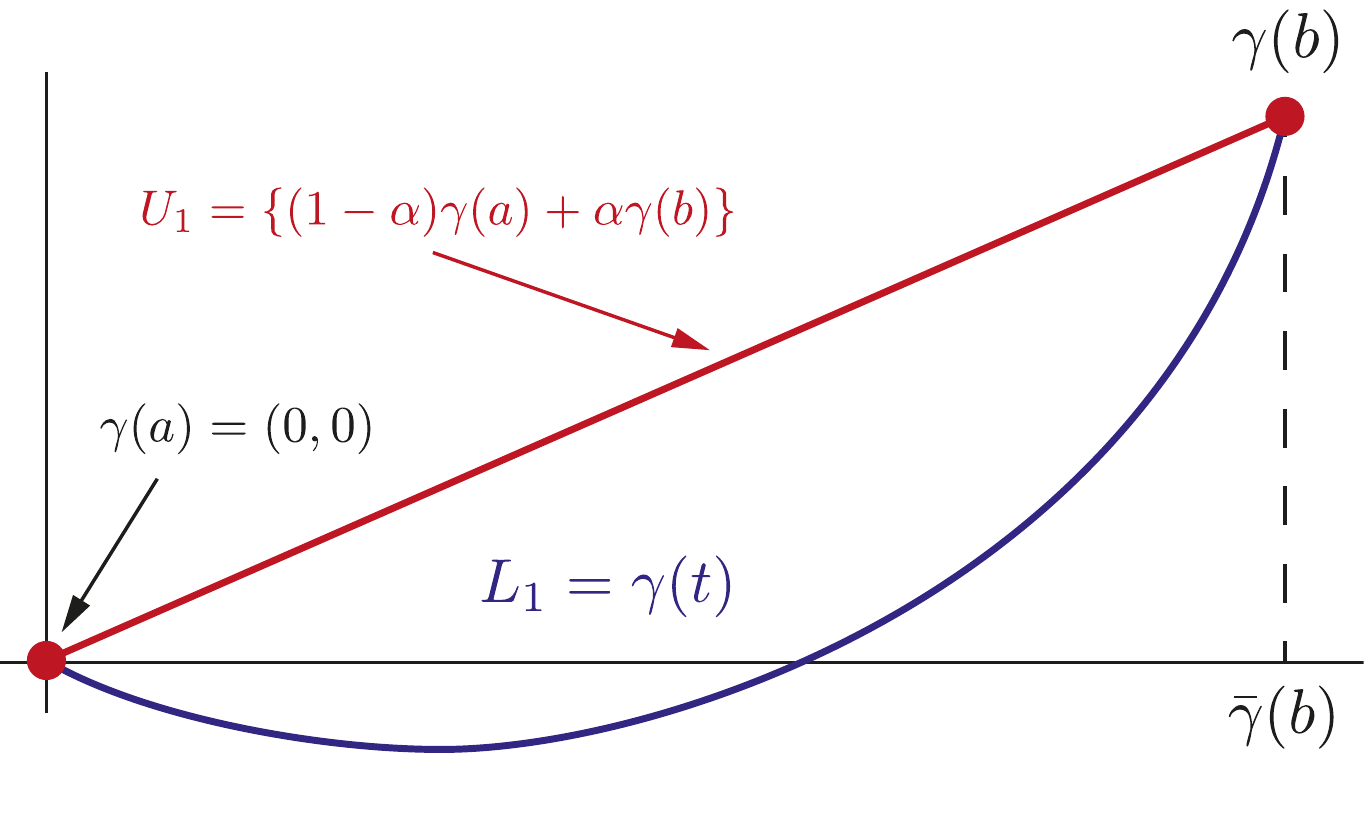}
\end{center}
\caption{ Proof of Theorem \ref{mth010} for dimension $n+1 = 1+1$.}
\label{fig:2d}
\end{figure}

This case is trivial and Theorem~\ref{mth010} essentially follows by looking at Fig.~\ref{fig:2d}.

If we reparametrize the curve $\gamma$ as  $\tilde{\gamma}(t):= \gamma(\gamma_{1}^{-1}(t))$, $t \in  (0,\gamma_{1}(1))$, then $\tilde{\gamma}$ has totally positive torsion. So  $\tilde{\gamma}(t) = (t, g(t)), t \in (0, \gamma_{1}(1))$ where $g(0)=0$, and $\frac{d^{2}}{dt^{2}}g(t)>0$ for  all $t\in (0,\gamma_{1}(1))$.  We have $U_{1}(\beta_{1})=\beta_{1}\gamma(1), \beta_{1} \in [0,1]$,  is the line joining the endpoints of $\tilde{\gamma}$. Also  $L_{1}(x_{1})=\gamma(x_{1}), x_{1} \in [0,1]$, is the curve coinciding with $\tilde{\gamma}$. It is easy to see that in this case Theorem~\ref{mth010} holds true.

\subsection{The proof of Theorem~\ref{mth010} in dimension 2+1}

\subsubsection{The lower hull}\label{3low}
Recall that
\begin{align*}
   \overline{L}_{2} :\Delta_{c}^{1}\times \Delta_{*}^{1} =   [0,1]^{2} \ni (\alpha, x) \mapsto \alpha \overline{\gamma}(x).
\end{align*}
We claim 
\begin{align}
    &\overline{L}_{2}(\partial ([0,1]^{2})) =\partial ( \mathrm{conv}(\overline{\gamma}([0,1]))); \label{3db}\\
   &\overline{L}_{2} :\mathrm{int}([0,1]^{2}) \mapsto \mathrm{int}(\mathrm{conv}(\overline{\gamma}([0,1]))) \quad \text{is diffeomorphism.}  \label{3ddiff}
\end{align}
To verify (\ref{3db}) it suffices  to show that $\overline{\gamma}$ is the convex curve in $\mathbb{R}^{2}$. Convexity of $\overline{\gamma}$ can be verified in a similar way as in Section~\ref{kkl1}. However, here we present one more proof which later will be adapted to higher dimensions too.   Assume contrary, i.e., there exists $0\leq a <b<c \leq 1$ such that $\overline{\gamma}(a), \overline{\gamma}(b), \overline{\gamma}(c)$ lie on the same line, i.e., 
\begin{align}\label{ura1}
 0=\det(\overline{\gamma}(b)-\overline{\gamma}(a), \overline{\gamma}(c)-\overline{\gamma}(b)) =  \int_{a}^{b} \int_{b}^{c} \det(\overline{\gamma}'(y_{1}), \overline{\gamma}'(y_{2}))dy_{1}dy_{2}.
\end{align}
The equation (\ref{ura1}) is in contradiction with Lemma~\ref{klasika} applied to $\overline{\gamma}$. 

\begin{figure}[t]
    \centering
    \includegraphics[width=.49\textwidth]{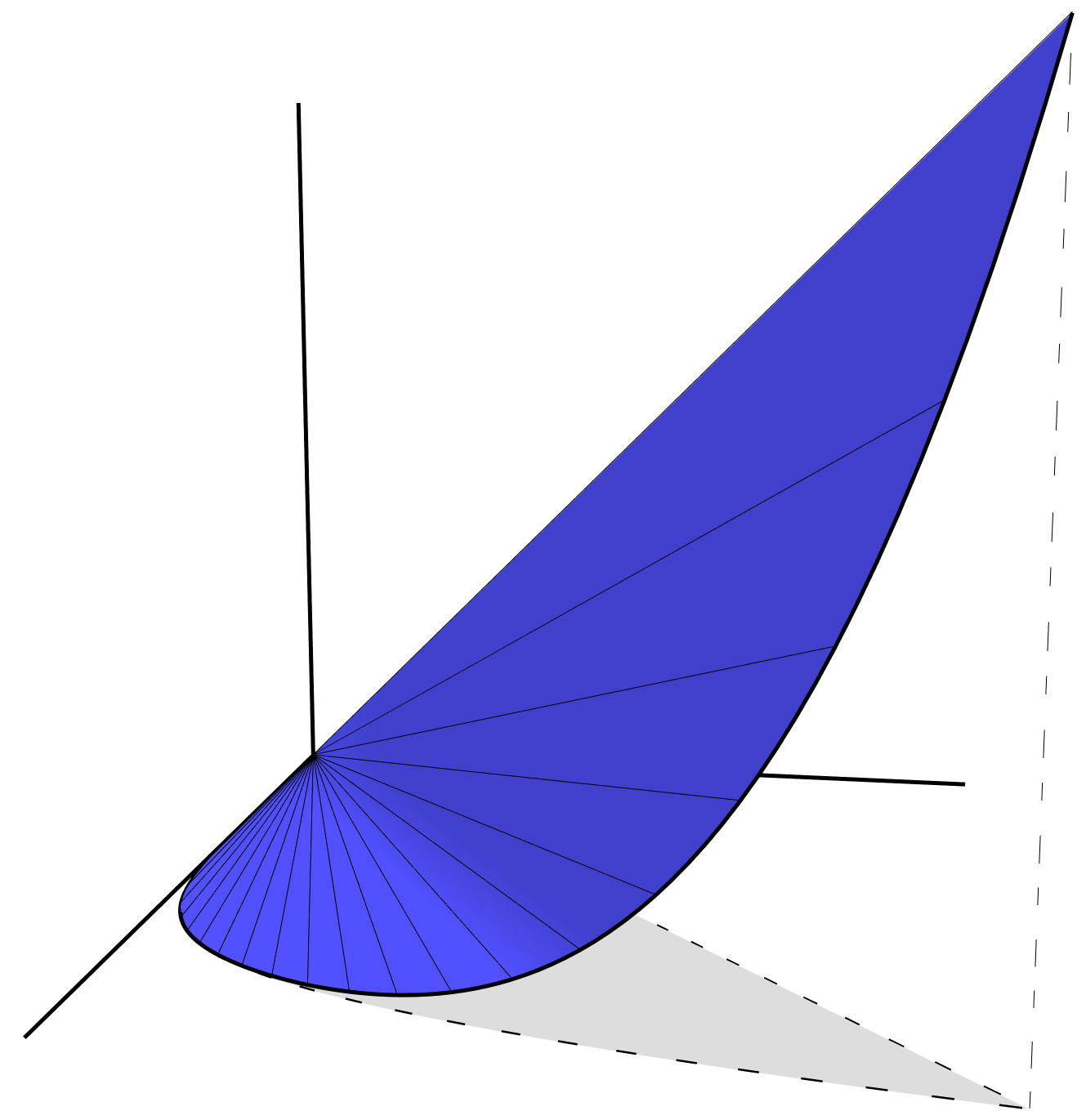}
    \includegraphics[width=.49\textwidth]{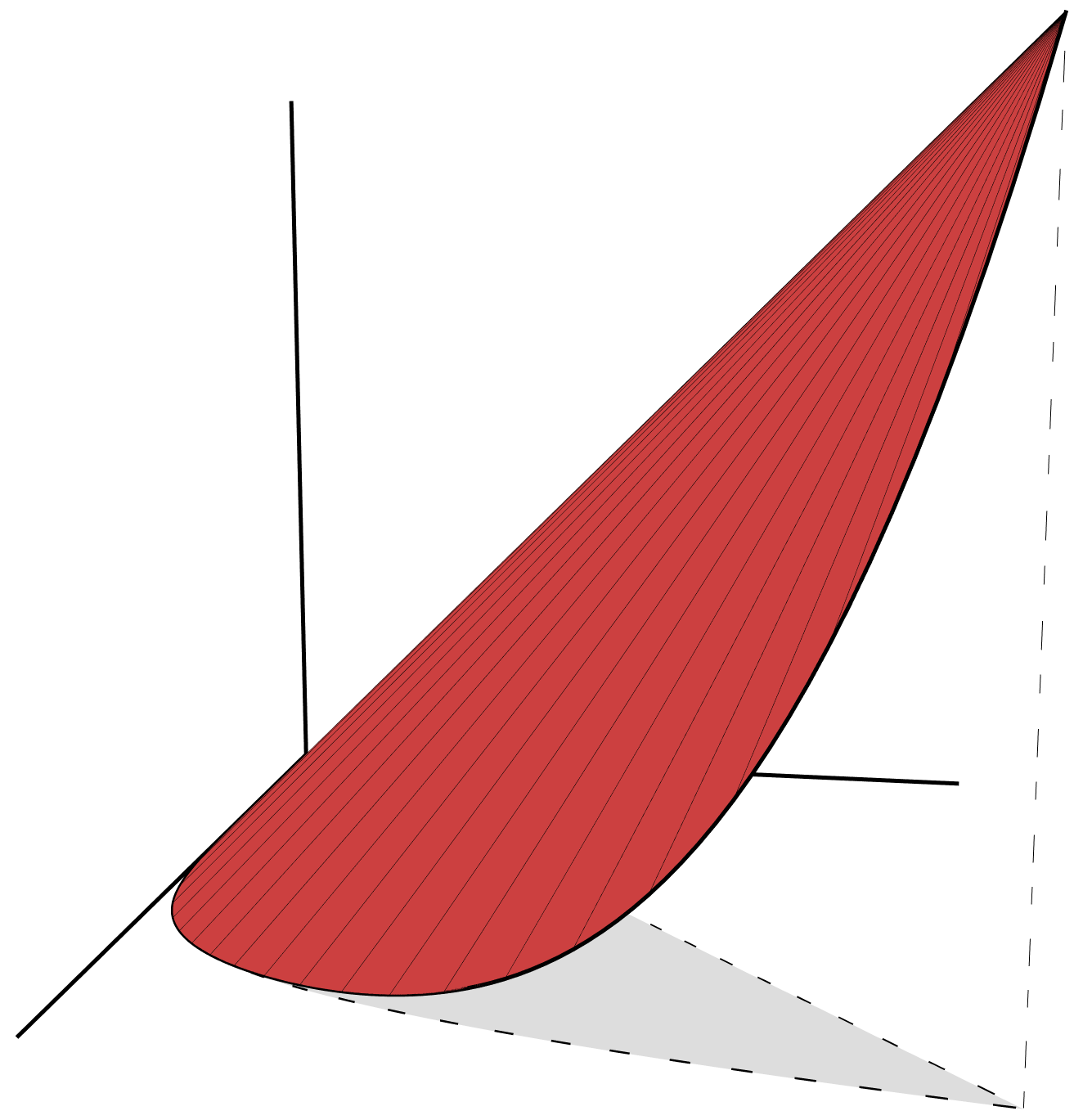}
    \caption{Two pieces of the boundary of the convex hull of $\gamma$: the lower hull $L_{2}$ (left) and the upper hull $U_{2}$}
    \label{fig:3d}
\end{figure}

To verify (\ref{3ddiff}), by Hadamard-Caccioppoli theorem it suffices to check  that the differential of $\overline{L}:=\overline{L}_{2}$ at the interior of $[0,1]^{2}$ has full rank, and the map $\overline{L}_{2}$ is injection. The injectivity will be verified later in all dimensions simultaneously (see the section on proofs of (\ref{diff2l-1u}), (\ref{diff2l-1l}), (\ref{diff2lu}), and  (\ref{diff2ll})).    To verify the full rank property we have  $D \overline{L} = (\overline{L}_{\alpha}, \overline{L}_{x}) =\alpha \det(\overline{\gamma}(x),\overline{\gamma}'(x))$. On the other hand 
\begin{align}\label{3dtozhd1}
    \det (\overline{\gamma}(x),\overline{\gamma}'(x)) = \int_{0}^{x} \det(\overline{\gamma}'(y_{1}), \overline{\gamma}'(x))dy_{1} \stackrel{\text{Lemma}~\ref{klasika}}{>} 0.
\end{align}
Thus, see Fig.~\ref{fig:3d},

$$
L_{2} : \Delta_{c}^{1}\times \Delta_{*}^{1}=[0,1]^{2} \ni (\alpha, x) \mapsto \alpha \gamma(x)
$$
parametrizes a surface in $\mathbb{R}^{3}$ which is a graph of a function $B^{\mathrm{inf}}$ defined on $\mathrm{conv}(\overline{\gamma}([0,1]))$ as follows
\begin{align*}
    B^{\mathrm{inf}}(\alpha \overline{\gamma}(x)) = \alpha \gamma_{3}(x), \quad \text{for all} \quad (\alpha, x) \in [0,1]^{2}. 
\end{align*}
 Let us check that $B^{\mathrm{inf}}$ is convex. Indeed, at any point $(\alpha_{0}, x_{0}) \in \mathrm{int}([0,1]^{2})$ the set of points $\xi \in \mathbb{R}^{3}$ belonging to the tangent  plane at point $L_{2}(\alpha_{0}, x_{0})$ is found as the solution of the equation  
\begin{align}\label{trieq}
\det(L_{\alpha}(\alpha_{0}, x_{0}),L_{x}(\alpha_{0}, x_{0}), \xi-L(\alpha_{0}, x_{0})) = \alpha_{0} \det(\gamma(x_{0}), \gamma'(x_{0}), \xi)=0.
\end{align}
For $\xi=e_{3}$, where $e_{3}=(0,0,1)$ we have 
\begin{align*}
    &\det(\gamma(x_{0}), \gamma'(x_{0}), e_{3}) = \det(\overline{\gamma}(x_{0}), \overline{\gamma}'(x_{0}))\stackrel{(\ref{3dtozhd1})}{>}0.
\end{align*}
Therefore, to verify the convexity of $B^{\mathrm{inf}}$, i.e., the surface $L([0,1]^{2})$ lies above the tangent plane at point $L(\alpha_{0}, x_{0}),$ it suffices to show that 
$$
\det(\gamma(x_{0}), \gamma'(x_{0}), L(\alpha, x)) = \alpha \det(\gamma(x_{0}), \gamma'(x_{0}), \gamma(x))  \geq 0.
$$
If $x=x_{0}$ there is nothing to prove. If $x>x_{0}$ then 
\begin{align*}
    \det(\gamma(x_{0}), \gamma'(x_{0}), \gamma(x)) = \int_{0}^{x_{0}}\int_{x_{0}}^{x}\det(\gamma'(y_{1}),\gamma'(x_{0}), \gamma'(y_{3}))dy_{1}dy_{3} \stackrel{\text{Lemma}~\ref{klasika}}{>}0.
\end{align*}
Similarly, if $x<x_{0}$, by Lemma~\ref{klasika} we have 
\begin{align*}
    \det(\gamma(x_{0}), \gamma'(x_{0}), \gamma(x)) = \int_{x}^{x_{0}} \int_{0}^{x}  \det(\gamma'(y_{1}),\gamma'(x_{0}), \gamma'(y_{3}))dy_{1} dy_{3} >0.
\end{align*}

To verify that $B^{\mathrm{inf}}$ is the maximal convex function defined on $\mathrm{conv}(\overline{\gamma}([0,1]))$ such that $B(\overline{\gamma}(s)) = \gamma_{3}(s)$, notice that since every point $(\xi,B^{\mathrm{inf}}(\xi))$, where $\xi \in \mathrm{conv}(\overline{\gamma}([0,1]))$, is the convex combination of some points of $\gamma$, it follows that any other candidate $\tilde{B}$ would be smaller than $B$ by convexity.

\subsubsection{The upper hull} \label{3up}
Consider the map 
\begin{align*}
    \overline{U}_{2} : \Delta_{c}^{1}\times \Delta_{*}^{1} = [0,1]^{2} \ni (\alpha, x ) \mapsto \alpha \overline{\gamma}(x)+(1-\alpha)\overline{\gamma}(1).
\end{align*}
Similarly as before $\Phi$ satisfies (\ref{3db}) and (\ref{3ddiff}). The property (\ref{3db}) follows from from the convexity of $\overline{\gamma}$. The property (\ref{3ddiff}) follows from 
\begin{align*}
    \det(\overline{U}_{\alpha}, \overline{U}_{x}) = \alpha \det(\overline{\gamma}(x)-\overline{\gamma}(1), \overline{\gamma}'(x))= \int_{x}^{1}\det(\overline{\gamma}'(x), \overline{\gamma}'(y_{2}))dy_{2} \neq 0
\end{align*}
for all $(\alpha, x) \in \mathrm{int}([0,1]^{2})$ by Lemma~\ref{klasika} applied to $\overline{\gamma}$. 

Next, we show that 
\begin{align*}
    B^{\mathrm{sup}}(\alpha \overline{\gamma}(x)+(1-\alpha)\overline{\gamma}(1)) = \alpha \gamma_{3}(x)+(1-\alpha)\gamma_{3}(1)
\end{align*}
defines a minimal concave function on $\mathrm{conv}(\overline{\gamma}([0,1]))$ with the property $B^{\mathrm{sup}}(\overline{\gamma})= \gamma_{3}$, see Fig 3. Let $U(\alpha, x) = \alpha \gamma(x)+(1-\alpha)\gamma(1)$. The equation of the tangent plane at point $U(\alpha_{0}, x_{0})$, where $(\alpha_{0}, x_{0}) \in \mathrm{int}([0,1]^{2})$, is given by 
\begin{align*}
&\det(U_{\alpha}(\alpha_{0}, x_{0}), U_{x}(\alpha_{0}, x_{0}), \xi -U(\alpha_{0}, x_{0})) =\alpha_{0}\det(\gamma(x_{0})-\gamma(1),\gamma'(x_{0}), \xi-\alpha_{0}(\gamma(x_{0})-\gamma(1))-\gamma(1))\\
&=\alpha_{0}\det(\gamma(x_{0})-\gamma(1),\gamma'(x_{0}), \xi-\gamma(1))=0.
\end{align*}
For $\xi=\lambda e_{3}$ with $\lambda \to +\infty$ we have 
\begin{align*}
    &\mathrm{sign}[\det(\gamma(x_{0})-\gamma(1),\gamma'(x_{0}), \lambda e_{3}-\gamma(1))] = \mathrm{sign}[ \det (\overline{\gamma}(x_{0})-\overline{\gamma}(1),\overline{\gamma}'(x_{0})]\\
    &=\mathrm{sign}\left[\int_{x_{0}}^{1}\det (\overline{\gamma}'(x_{0}) ,\overline{\gamma}(y_{2})) dy_{2}\right] >0
\end{align*}
by Lemma~\ref{klasika} applied to $\overline{\gamma}$. Therefore, the concavity of $B^{\mathrm{sup}}$ would follow from the following inequality 
\begin{align*}
   \det(\gamma(x_{0})-\gamma(1),\gamma'(x_{0}), U(\alpha, x)-\gamma(1))=\alpha  \det(\gamma(x_{0})-\gamma(1),\gamma'(x_{0}), \gamma(x)-\gamma(1)) \leq 0
\end{align*}
for all $x_{0}, \alpha, x \in [0,1]$. If $x=x_{0}$ there is nothing to prove. Consider $x>x_{0}$ (the case $x<x_{0}$ is similar). Then 
\begin{align*}
&\det(\gamma(x_{0})-\gamma(1),\gamma'(x_{0}), \gamma(x)-\gamma(1))=\det(\gamma(x_{0})-\gamma(1),\gamma'(x_{0}), \gamma(x)-\gamma(x_{0}))=\\
&-\det(\gamma(x_{0})-\gamma(x), \gamma'(x_{0}), \gamma(1)-\gamma(x_{0})) = -\int_{x}^{x_{0}}\int_{x_{0}}^{1}\det(\gamma'(y_{1}), \gamma'(x_{0}), \gamma'(y_{2}))dy_{2}dy_{1}<0
\end{align*}
by Lemma~\ref{klasika}. 

\vskip1cm

The properties (\ref{giff}) and (\ref{union}) will be verified in sections \ref{giffsub} and \ref{unionsub}.

%Let us also observe that it follows from the construction 
%\begin{align*}
%B^{\mathrm{sup}}(x)=B^{\mathrm{inf}}(x) \quad \text{for all}\quad x \in  \partial\, \mathrm{conv}(\overline{\gamma}([0,1]))
%\end{align*}
%and 
%\begin{align}\label{granica21}
% \partial\,  \mathrm{conv}(\gamma([0,1]))=\{(x,B^{\mathrm{sup}}(x)), x \in \mathrm{conv}(\bar{\gamma}([0,1]))\} \cup \{(x,B^{\mathrm{inf}}(x)), x \in \mathrm{conv}(\bar{\gamma}([0,1]))\}.
%\end{align}

\subsection{The proof of Theorem~\ref{mth010} in an arbitrary dimension $n+1$}

\begin{proof}
Since Theorem~\ref{mth010} contains several statements the whole proof will be split into several parts. 

{\em The proof of claims (\ref{b2l-1}) and (\ref{b2l}).}

The proof will be by induction on $n$. We have checked the statement for $n=1,2$. First we consider the case when $n=2\ell-1$. We shall verify the  claim (\ref{b2l-1}) by showing that $\overline{U}_{2\ell-1}|_{\partial(\Delta_{c}^{\ell} \times \Delta_{*}^{\ell-1})}$, i.e.,  the restriction of $\overline{U}_{2\ell-1}$ on $\partial(\Delta_{c}^{\ell} \times \Delta_{*}^{\ell-1})$, coincides with  maps $U_{2\ell-2}$ and $L_{2\ell-2}$ (similarly for $\overline{L}_{2\ell-1}|_{\partial(\Delta_{c}^{\ell-1} \times \Delta_{*}^{\ell})}$). Since by the induction  the union of the images of $U_{2\ell-2}$ and $L_{2\ell-2}$ coincides with the boundary of the convex hull of $\overline{\gamma}([0,1])$, see (\ref{union}), we obtain the claim.

Recall that 
\begin{align*}
\overline{U}_{2\ell-1} :\Delta_{c}^{\ell} \times \Delta_{*}^{\ell-1} \ni  (\beta_{1}, \ldots, \beta_{\ell}, y_{2},\ldots, y_{\ell}) \mapsto  \beta_{1} \overline{\gamma}(1)+\sum_{j=2}^{\ell} \beta_{j} \overline{\gamma}(y_{j}),
\end{align*}
and 
\begin{align*}
&U_{2\ell-2} :   \Delta_{c}^{\ell-1} \times \Delta_{*}^{\ell-1} \ni (\lambda_{1}, \ldots, \lambda_{\ell-1}, x_{1}, \ldots, x_{\ell-1}) \mapsto \sum_{j=1}^{\ell-1} \lambda_{j} \overline{\gamma}(x_{j})  + (1-\sum_{j=1}^{\ell-1}\lambda_{j}) \overline{\gamma}(1),\\
&L_{2\ell-2} :\Delta_{c}^{\ell-1} \times \Delta_{*}^{\ell-1} \ni  (\lambda_{1}, \ldots, \lambda_{\ell-1}, z_{1}, \ldots, z_{\ell-1}) \mapsto \sum_{j=1}^{\ell-1} \lambda_{j} \overline{\gamma}(z_{j}).
\end{align*}

If $\beta_{1}=0$ then $\overline{U}_{2\ell-1}$ coincides with $L_{2\ell-2}$. If $\sum_{j=1}^{n}\beta_{j}=1$, i.e.,  $\beta_{1}=1-\sum_{j=2}^{\ell}\beta_{j}$, then $\overline{U}_{2\ell-1}$ coincides with $U_{2\ell-2}$. Thus, we have
\begin{align*}
   \partial\,  \mathrm{conv}(\overline{\gamma}([0,1])) \stackrel{\mathrm{induction}}{=} U_{2\ell-2}(\Delta_{c}^{\ell-1} \times \Delta_{*}^{\ell-1}) \cup L_{2\ell-2}(\Delta_{c}^{\ell-1} \times \Delta_{*}^{\ell-1}) \subset  \overline{U}_{2\ell-1}(\partial\,  (\Delta_{c}^{\ell} \times \Delta_{*}^{\ell-1})). 
\end{align*}
On the other hand, if $\beta_{p}=0$ for some $p \in \{2, \ldots, \ell\}$, then $\overline{U}_{2\ell-1}$ coincides with  $L_{2\ell-2}$ restricted to  $z_{1}=1$. If at least one of the following conditions hold: a) $y_{2}=0$; b) $y_{s}=y_{s+1}$ for some $s \in \{2, \ldots, \ell-1\}$; c) $y_{\ell}=1$,  then $\overline{U}_{2\ell-1}$ coincides with $U_{2\ell-2}$ restricted to $x_{1}=0$. Thus we obtain $\partial\,  \mathrm{conv}(\overline{\gamma}([0,1])) =  \overline{U}_{2\ell-1}(\partial\,  (\Delta_{c}^{\ell} \times \Delta_{*}^{\ell-1}))$. 

Next, we verify that $\partial\,  \mathrm{conv}(\overline{\gamma}([0,1])) =  \overline{L}_{2\ell-1}(\partial\,  (\Delta_{c}^{\ell-1} \times \Delta_{*}^{\ell}))$. We recall
\begin{align*}
\overline{L}_{2\ell-1} :\Delta_{c}^{\ell-1} \times \Delta_{*}^{\ell} \ni (\beta_{2}, \ldots, \beta_{\ell}, y_{1},\ldots, y_{\ell}) \mapsto  \sum_{j=2}^{\ell} \beta_{j} \overline{\gamma}(y_{j})+(1-\sum_{j=2}^{\ell} \beta_{j})\overline{\gamma}(y_{1}).
\end{align*}
If $y_{\ell}=1$ then $\overline{L}_{2\ell-1}$ coincides with $U_{2\ell-2}$. If $y_{1}=0$ then  $\overline{L}_{2\ell-1}$ coincides with $L_{2\ell-2}$. Thus, by induction $\partial\,  \mathrm{conv}(\overline{\gamma}([0,1])) \subset   \overline{L}_{2\ell-1}(\partial\,  (\Delta_{c}^{\ell-1} \times \Delta_{*}^{\ell}))$. 

Next, if $y_{s}=y_{s+1}$ for some $s \in \{1, \ldots, \ell-1\}$ then $\overline{L}_{2\ell-1}$ coincides with  $L_{2\ell-2}$ restricted to  $\lambda_{1}=1-\sum_{j=2}^{\ell-1}\lambda_{j}$. Also, if $\sum_{j=2}^{\ell}\beta_{j}=1$ then $\overline{L}_{2\ell-1}$ coincides with $L_{2\ell-2}$. Finally, if $\beta_{s}=0$ for some $s \in \{2, \ldots, \ell\}$ then $\overline{L}_{2\ell-1}$ coincides with $L_{2\ell-2}$ restricted to $\sum_{j=1}^{\ell-1}\lambda_{j}=1$. Thus we obtain 
$\partial\,  \mathrm{conv}(\overline{\gamma}([0,1])) =   \overline{L}_{2\ell-1}(\partial\,  (\Delta_{c}^{\ell-1} \times \Delta_{*}^{\ell}))$.

Next, we assume $n=2\ell$. First we verify (\ref{b2l}). As before we  claim that the restriction of $\overline{U}_{2\ell}$ on $\partial(\Delta_{c}^{\ell} \times \Delta_{*}^{\ell})$ coincides with  maps $U_{2\ell-1}$ and $L_{2\ell-1}$ (similarly for $\overline{L}_{2\ell}$). Since by the induction  the union of the images of $U_{2\ell-1}$ and $L_{2\ell-1}$ coincide with the boundary of the convex hull of $\overline{\gamma}([0,1])$, see (\ref{union}), we obtain the claim. 

We recall that 
\begin{align*}
\overline{U}_{2\ell} :   \Delta_{c}^{\ell} \times \Delta_{*}^{\ell} \ni (\lambda_{1}, \ldots, \lambda_{\ell}, x_{1}, \ldots, x_{\ell}) \mapsto \sum_{j=1}^{\ell} \lambda_{j} \overline{\gamma}(x_{j})  + (1-\sum_{j=1}^{\ell}\lambda_{j}) \overline{\gamma}(1);
\end{align*}
and 
\begin{align*}
   &U_{2\ell-1} :\Delta_{c}^{\ell} \times \Delta_{*}^{\ell-1} \ni  (\beta_{1}, \ldots, \beta_{\ell}, y_{2},\ldots, y_{\ell}) \mapsto  \beta_{1} \overline{\gamma}(1)+\sum_{j=2}^{\ell} \beta_{j} \overline{\gamma}(y_{j});\\
   &L_{2\ell-1} :\Delta_{c}^{\ell-1} \times \Delta_{*}^{\ell} \ni (\beta_{2}, \ldots, \beta_{\ell}, z_{1},\ldots, z_{\ell}) \mapsto  (1-\sum_{j=2}^{\ell} \beta_{j})\overline{\gamma}(z_{1})+\sum_{j=2}^{\ell} \beta_{j} \overline{\gamma}(z_{j}).
\end{align*}
Notice that if $\sum_{j=1}^{\ell} \lambda_{j}=1$ then  $\overline{U}_{2\ell}$ coincides with $L_{2\ell-1}$. On the other hand, if $x_{1}=0$ then $\overline{U}_{2\ell}$ coincides with $U_{2\ell-1}$. Thus, by induction we have $\partial\,  \mathrm{conv}(\overline{\gamma}([0,1])) \subset \overline{U}_{2\ell}(\partial\,  (\Delta_{c}^{\ell} \times \Delta_{*}^{\ell}))$. Also notice that if $\lambda_{p}=0$ for some $p \in \{1, \ldots, \ell\}$ (or $x_{s}=x_{s+1}$ for some $s \in \{1, \ldots, \ell-1\}$, or $x_{\ell}=1$) then $\overline{U}_{2\ell}$ coincides with $U_{2\ell-1}$ restricted to the boundary of  $\Delta_{c}^{\ell-1} \times \Delta_{*}^{\ell}$  (if $\lambda_{p}=0$ or $x_{\ell}=1$ take $\beta_{1} = 1-\sum_{j=2}^{\ell}\beta_{j}$, if $x_{s}=x_{s+1}$ take $\beta_{1} = 1-\sum_{j=2}^{\ell}\beta_{j}$). Thus we obtain $\partial\,  \mathrm{conv}(\overline{\gamma}([0,1]))  = \overline{U}_{2\ell}(\partial\,  (\Delta_{c}^{\ell} \times \Delta_{*}^{\ell}))$. 

Next, we verify the claim $\overline{L}_{2\ell}(\partial\,  (\Delta_{c}^{\ell} \times \Delta_{*}^{\ell})) = \partial\,  \mathrm{conv}(\overline{\gamma}([0,1]))$.
We recall that 
\begin{align*}
\overline{L}_{2\ell} :\Delta_{c}^{\ell} \times \Delta_{*}^{\ell} \ni  (\lambda_{1}, \ldots, \lambda_{\ell}, x_{1}, \ldots, x_{\ell}) \mapsto \sum_{j=1}^{\ell} \lambda_{j} \overline{\gamma}(x_{j}),
\end{align*}
If $\sum_{j=1}^{\ell}\lambda_{j}=1$ then $\overline{L}_{2\ell}$ coincides with $L_{2\ell-1}$. If $x_{\ell}=1$ then $\overline{L}_{2\ell}$ coincides with $U_{2\ell-1}$. Thus by induction we have $\partial\,  \mathrm{conv}(\overline{\gamma}([0,1]))  \subset \overline{L}_{2\ell}(\partial\,  (\Delta_{c}^{\ell} \times \Delta_{*}^{\ell}))$

If $\lambda_{p}=0$ for some $p \in \{1, \ldots, \ell\}$, or $x_{1}=0$,  then $\overline{L}_{2\ell}$ coincides with $U_{2\ell-1}$ if we choose $\beta_{1}=0$. Finally, if $x_{s}=x_{s+1}$ for some $s \in \{1, \ldots, \ell-1\}$, then $\overline{L}_{2\ell}$ coincides with $U_{2\ell-1}$ if we choose $\beta_{1}=0$, and $\beta_{s+1}=\lambda_{s}+\lambda_{s+1}$. Therefore,  we have $\overline{L}_{2\ell}(\partial\,  (\Delta_{c}^{\ell} \times \Delta_{*}^{\ell})) \subset \partial\,  \mathrm{conv}(\overline{\gamma}([0,1]))$, and the claim (\ref{b2l}) is verified.

{\em The proof of claims (\ref{diff2l-1u}), (\ref{diff2l-1l}), (\ref{diff2lu}) and (\ref{diff2ll}).}

We start by showing that the Jacobian of the map $\overline{U}_{n}$ has full rank at the interior points of its domain.  Hence the map is local diffeomoerphism by inverse function theorem. Therefore, the map is surjective, otherwise the image of its domain would have a boundary in the interior of the codomain (boundary goes to boundary by (\ref{b2l}) and (\ref{b2l-1})) and this would contradict the local diffeomoerphism. Next, we  show that the map $\overline{U}_{n}$ is injective, and hence proper. So we conclude that $\overline{U}_{n}$ is  diffeomorphism. The similar reasoning will be done for $\overline{L}_{n}$.

First we verify that the Jacobian matrices  $\nabla \overline{U}_{n}$ and $\nabla \overline{L}_{n}$ have full rank at the interior points of their domains.  

Assume $n=2\ell-1$. We have 
\begin{align*}
    &\det(\nabla \overline{U}_{2\ell-1}) =\det(\overline{\gamma}(1), \overline{\gamma}(x_{2}), \ldots, \overline{\gamma}(x_{\ell}), \beta_{2} \overline{\gamma}'(x_{2}), \ldots, \beta_{\ell}\overline{\gamma}'(x_{\ell}))\\
%    &=\pm \det(\overline{\gamma}(x_{2}), \overline{\gamma}'(x_{2}), \overline{\gamma}(x_{3}), \overline{\gamma}'(x_{3}), \ldots, \overline{\gamma}(x_{\ell}), \overline{\gamma}'(x_{\ell}), \overline{\gamma}(1))\prod_{j=2}^{\ell} \beta_{j} \\
    &=\pm  \det(\overline{\gamma}(x_{2}), \overline{\gamma}'(x_{2}), \overline{\gamma}(x_{3}), \overline{\gamma}'(x_{3}), \ldots, \overline{\gamma}(x_{\ell}), \overline{\gamma}'(x_{\ell}), \overline{\gamma}(1)) \prod_{j=2}^{\ell} \beta_{j}\\
    & = \pm  \det(\overline{\gamma}(x_{2})-\overline{\gamma}(0), \overline{\gamma}'(x_{2}), \overline{\gamma}(x_{3})-\overline{\gamma}(x_{2}), \overline{\gamma}'(x_{3}), \ldots, \overline{\gamma}(x_{\ell})-\overline{\gamma}(x_{\ell-1}), \overline{\gamma}'(x_{\ell}), \overline{\gamma}(1)-\overline{\gamma}(x_{\ell})) \prod_{j=2}^{\ell} \beta_{j}\\
    &=\pm  \prod_{j=2}^{\ell}\beta_{j}\,  \int_{x_{\ell}}^{1}\ldots \int_{x_{2}}^{x_{3}} \int_{0}^{x_{2}}  \det(\overline{\gamma}'(s_{1}), \overline{\gamma}'(x_{2}),\overline{\gamma}'(s_{2}), \ldots, \overline{\gamma}'(x_{\ell}), \overline{\gamma}'(s_{\ell}))ds_{1} ds_{2}\ldots ds_{\ell}.
\end{align*}
Thus $\det(\nabla \overline{U}_{2\ell-1})$ is nonzero by Lemma~\ref{klasika}. 

Next, we verify that $\det(\nabla \overline{L}_{2\ell-1})\neq 0$, Indeed, 
\begin{align*}
&\det(\nabla \overline{L}_{2\ell-1}) = \\
&\det( \overline{\gamma}(x_{2})-\overline{\gamma}(x_{1}), \overline{\gamma}(x_{3})-\overline{\gamma}(x_{1}), \ldots, \overline{\gamma}(x_{\ell})-\overline{\gamma}(x_{1}), \overline{\gamma}'(x_{1}), \ldots, \overline{\gamma}'(x_{\ell}) ) (1-\sum_{j=2}^{\ell}\beta_{j})\prod_{j=2}^{\ell}\beta_{j}=\\
& \det( \overline{\gamma}(x_{2})-\overline{\gamma}(x_{1}), \overline{\gamma}(x_{3})-\overline{\gamma}(x_{2}), \ldots, \overline{\gamma}(x_{\ell})-\overline{\gamma}(x_{\ell-1}), \overline{\gamma}'(x_{1}), \ldots, \overline{\gamma}'(x_{\ell}) ) (1-\sum_{j=2}^{\ell}\beta_{j})\prod_{j=2}^{\ell}\beta_{j}=\\
&\pm \det(\overline{\gamma}'(x_{1}),  \overline{\gamma}(x_{2})-\overline{\gamma}(x_{1}), \overline{\gamma}'(x_{2}), \overline{\gamma}(x_{3})-\overline{\gamma}(x_{2}), \ldots, \overline{\gamma}(x_{\ell})-\overline{\gamma}(x_{\ell-1}), \overline{\gamma}'(x_{\ell}) ) (1-\sum_{j=2}^{\ell}\beta_{j})\prod_{j=2}^{\ell}\beta_{j}=\\
&\pm (1-\sum_{j=2}^{\ell}\beta_{j})\prod_{j=2}^{\ell}\beta_{j} \times \\
&\int_{x_{\ell-1}}^{x_{\ell}}\ldots \int_{x_{2}}^{x_{3}}\int_{x_{1}}^{x_{2}}\det( \overline{\gamma}'(x_{1}), \overline{\gamma}'(s_{1}), \overline{\gamma}'(x_{2}), \overline{\gamma}'(s_{2}), \ldots, \overline{\gamma}'(s_{\ell-1}), \overline{\gamma}'(x_{\ell})) ds_{1} ds_{2}\ldots ds_{\ell-1} \neq 0
\end{align*}
by Lemma~\ref{klasika}.

Assume $n=2\ell$. We have 
\begin{align*}
    &\det(\nabla \overline{U}_{2\ell}) = \det(\overline{\gamma}(x_{1})-\overline{\gamma}(1), \ldots, \overline{\gamma}(x_{\ell})-\overline{\gamma}(1), \overline{\gamma}'(x_{1}), \ldots, \overline{\gamma}'(x_{\ell})) \prod_{j=1}^{\ell} \lambda_{j}=\\
    &\pm \det(\overline{\gamma}'(x_{1}), \overline{\gamma}(x_{1})-\overline{\gamma}(x_{2}), \overline{\gamma}'(x_{2}), \overline{\gamma}(x_{2})-\overline{\gamma}(x_{3}), \ldots, \overline{\gamma}'(x_{\ell}), \overline{\gamma}(x_{\ell})-\overline{\gamma}(1))\prod_{j=1}^{\ell} \lambda_{j} = \\
    & \pm \int_{x_{\ell}}^{1}\ldots \int_{x_{2}}^{x_{3}}\int_{x_{1}}^{x_{2}}\det(\overline{\gamma}'(x_{1}), \overline{\gamma}'(s_{1}), \overline{\gamma}'(x_{2}), \overline{\gamma}'(s_{2}), \ldots, \overline{\gamma}'(x_{\ell}), \overline{\gamma}'(s_{\ell}))ds_{1}ds_{2}\ldots ds_{\ell}    \prod_{j=1}^{\ell}\lambda_{j}
\end{align*}
which is nonzero by Lemma~\ref{klasika}. 

Finally, we verify $\det(\nabla \overline{L}_{2\ell}) \neq 0$. We have 
\begin{align*}
&\det(\nabla \overline{L}_{2\ell}) =\det(\overline{\gamma}(x_{1}), \ldots, \overline{\gamma}(x_{\ell}), \overline{\gamma}'(x_{1}), \ldots, \overline{\gamma}'(x_{\ell})) \prod_{j=1}^{\ell}\lambda_{j}=\\
&\pm \det(\overline{\gamma}(x_{1})-\overline{\gamma}(0),\overline{\gamma}'(x_{1}), \overline{\gamma}(x_{2})-\overline{\gamma}(x_{1}), \overline{\gamma}'(x_{2}), \ldots,\overline{\gamma}(x_{\ell})-\overline{\gamma}(x_{\ell-1}),\overline{\gamma}'(x_{\ell})) \prod_{j=1}^{\ell}\lambda_{j} = \\
&\pm  \int_{x_{\ell-1}}^{x_{\ell}} \ldots \int_{x_{1}}^{x_{2}}\int_{0}^{x_{1}}\det(\overline{\gamma}'(s_{1}), \overline{\gamma}'(x_{1}), \overline{\gamma}'(s_{2}), \overline{\gamma}'(x_{2}), \ldots, \overline{\gamma}'(s_{\ell}), \overline{\gamma}'(x_{\ell})) ds_{1} ds_{2}\ldots ds_{\ell}\prod_{j=1}^{\ell}\lambda_{j}.
\end{align*}
Thus $\det(\nabla \overline{L}_{2\ell}) \neq 0$ by Lemma~\ref{klasika}. 

Next, we show that the map $\overline{U}_{n}$  is injective in the interior of its domain. Assume $n=2\ell$. Let $(\lambda_{1}, \ldots, \lambda_{\ell}, x_{1}, \ldots, x_{\ell})$ and $(\beta_{1}, \ldots, \beta_{\ell}, y_{1}, \ldots, y_{\ell})$ be two different points in $\mathrm{int}(\Delta_{c}^{\ell}\times \Delta_{*}^{\ell})$ such that $\overline{U}_{\ell}$ takes the same values on these points. Then 
\begin{align}\label{linin}
    \sum_{j=1}^{\ell}\lambda_{j}(\overline{\gamma}(x_{j})-\overline{\gamma}(1)) - \sum_{k=1}^{\ell} \beta_{k} (\overline{\gamma}(y_{k})-\overline{\gamma}(1))=0. 
\end{align}
We claim that (\ref{linin})  holds if and only if $x_{j}=y_{j}$ and $\lambda_{j}=\beta_{j}$ for all $j=1, \ldots, \ell$. Indeed, we need the following 
\begin{lemma}\label{ltorsion}
For any numbers $z_{j}$, $1\leq j \leq 2\ell$, such that  $0<z_{1}<z_{2}<\ldots<z_{2\ell}\leq  1$, and any $r \in [0,1]\setminus\{z_{1}, \ldots, z_{2\ell}\}$, the vectors $\overline{\gamma}(z_{1})-\overline{\gamma}(r), \ldots, \overline{\gamma}(z_{2\ell})-\overline{\gamma}(r)$ are linearly independent in $\mathbb{R}^{2\ell}$. 
\end{lemma} 
\begin{proof}
The lemma follows from Corollary~\ref{klasikac} applied to $\beta=\overline{\gamma}$.
\end{proof}
Let $N$ be the cardinality of the set $Q=\{x_{1}, \ldots, x_{\ell}\} \cap \{y_{1}, \ldots, y_{\ell}\}$. If $N=\ell$ then necessarily $x_{j}=y_{j}$ for all $j=1, \ldots, \ell$, and the equation (\ref{linin}) combined with Lemma~\ref{ltorsion} implies that $\lambda_{j}=\beta_{j}$ for all $j=1, \ldots, \ell$. Therefore, assume $N<\ell$. Then we can split the sum (\ref{linin}) into the sum of 3 terms: the sum of $\lambda_{j} (\overline{\gamma}(x_{j})-\overline{\gamma}(1))$ where $x_{j} \notin Q$; the sum $(\lambda_{j}-\beta_{i_{j}})(\overline{\gamma}(x_{j})-\overline{\gamma}(1))$ where $x_{j} \in Q$; and the sum $\beta_{j} (\overline{\gamma}(y_{j})-\overline{\gamma}(1))$ where $y_{j}\notin Q$. Since $\beta_{j}$ and $\lambda_{j}$ cannot be zero, then applying Lemma~\ref{ltorsion} with $r=1$ we get a contradiction.

Next, we verify the injectivity of $\overline{L}_{2\ell}$ on the interior of its domain. Let $(\lambda_{1}, \ldots, \lambda_{\ell}, x_{1}, \ldots, x_{\ell})$ and $(\beta_{1}, \ldots, \beta_{\ell}, y_{1}, \ldots, y_{\ell})$ belong to $\mathrm{int}(\Delta_{c}^{\ell}\times \Delta_{*}^{\ell})$ and satisfy 
\begin{align*}
    \sum_{j=1}^{\ell} \lambda_{j} \overline{\gamma}(x_{j}) -\sum_{k=1}^{\ell}\beta_{k} \overline{\gamma}(y_{k})=0. 
\end{align*}
By applying Lemma~\ref{ltorsion} with $r=0$ and invoking the set $Q$ as before  we obtain   $x_{j}=y_{j}$, $\lambda_{j}=\beta_{j}$ for all $j=1, \ldots, \ell$. 

Assume $n=2\ell-1$. To verify the injectivity of $\overline{U}_{2\ell-1}$ on the interior of $\Delta_{c}^{\ell}\times \Delta_{*}^{\ell-1}$ we pick points $(\lambda_{1}, \ldots, \lambda_{\ell}, x_{2}, \ldots, x_{\ell})$ and $(\beta_{1}, \ldots, \beta_{\ell}, y_{2}, \ldots, y_{\ell})$ from $\mathrm{int}(\Delta_{c}^{\ell}\times \Delta_{*}^{\ell-1})$, and we assume 
\begin{align}\label{kidev1}
    (\lambda_{1}-\beta_{1})\overline{\gamma}(1)+\sum_{j=2}^{\ell}\lambda_{j} \overline{\gamma}(x_{j})-\sum_{j=2}^{\ell}\beta_{j} \overline{\gamma}(y_{j})=0.
\end{align}
\begin{lemma}\label{lltorsion}
For any numbers $0<z_{1}<\ldots<z_{2\ell-2}<1$ the vectors $\overline{\gamma}(z_{1}), \ldots, \overline{\gamma}(z_{2\ell-2}), \overline{\gamma}(1)$ are linearly independent in $\mathbb{R}^{2\ell-1}$.
\end{lemma}
\begin{proof}
The lemma follows from Corollary~\ref{klasikac} applied to $\beta=\overline{\gamma}$, $z_{2\ell-1}=1$, and $r=0$. 
\end{proof}
Invoking the set $Q$, and repeating the same reasoning as in the case of injectivity of $\overline{U}_{2\ell}$, we see that the equality (\ref{kidev1}) combined with Lemma~\ref{lltorsion} implies $x_{j}=y_{j}$ for all $j=2, \ldots, \ell$, and $\lambda_{j}=\beta_{j}$ for all $j=1, \ldots, \ell$. 

To verify the injectivity of $\overline{L}_{2\ell-1}$ on the interior of $\Delta_{c}^{\ell-1}\times \Delta_{*}^{\ell}$ we pick points $(\lambda_{2}, \ldots, \lambda_{\ell}, x_{1}, \ldots, x_{\ell})$ and $(\beta_{2}, \ldots, \beta_{\ell}, y_{1}, \ldots, y_{\ell})$ from $\mathrm{int}(\Delta_{c}^{\ell-1}\times \Delta_{*}^{\ell})$, and we assume 
\begin{align}\label{eq0011}
    (1-\sum_{j=2}^{\ell}\lambda_{j})\overline{\gamma}(x_{1})+\sum_{j=2}^{\ell}\lambda_{j}\overline{\gamma}(x_{j})=(1-\sum_{j=2}^{\ell}\beta_{j})\overline{\gamma}(y_{1})+\sum_{j=2}^{\ell}\beta_{j}\overline{\gamma}(y_{j}).
\end{align}
Without loss of generality assume $y_{1}\leq x_{1}$. We rewrite (\ref{eq0011}) as follows 
\begin{align}\label{sum01}
    (1-\sum_{j=2}^{\ell}\lambda_{j})(\overline{\gamma}(x_{1})-\overline{\gamma}(y_{1}))+\sum_{j=2}^{\ell}\lambda_{j}(\overline{\gamma}(x_{j})-\overline{\gamma}(y_{1}))-\sum_{j=2}^{\ell}\beta_{j}(\overline{\gamma}(y_{j})-\overline{\gamma}(y_{1}))=0.
\end{align}
Notice that if the points $x_{1}, \ldots, x_{\ell}, y_{1}, \ldots, x_{\ell}$ are different from each other,  and they belong to the interval $(0,1)$, then the vectors $\overline{\gamma}(x_{1})-\overline{\gamma}(y_{1}), \ldots, \overline{\gamma}(x_{\ell})-\overline{\gamma}(y_{1}), \overline{\gamma}(y_{2})-\overline{\gamma}(y_{1}), \ldots, \overline{\gamma}(y_{\ell})-\overline{\gamma}(y_{1})$ are linearly independent. The proof of the linear independence proceeds absolutely in the same way as the proof of Lemma~\ref{ltorsion} therefore we omit the proof to avoid the repetitions. Let $Q = \{x_{2}, \ldots, x_{\ell}\}\cap \{y_{2}, \ldots, y_{\ell}\}$, $X=\{x_{2}, \ldots, x_{\ell}\}$ and $Y=\{y_{2}, \ldots, y_{\ell}\}$. Then (\ref{sum01}) takes the form 
\begin{align} 
     &(1-\sum_{j=2}^{\ell}\lambda_{j})(\overline{\gamma}(x_{1})-\overline{\gamma}(y_{1}))+\sum_{j\, :\, x_{j} \in X\setminus Q}\lambda_{j}(\overline{\gamma}(x_{j})-\overline{\gamma}(y_{1}))+\nonumber\\
     &\sum_{j\, :\, x_{j} \in Q}(\lambda_{j}-\beta_{k_{j}})(\overline{\gamma}(x_{j})-\overline{\gamma}(y_{1}))
     -\sum_{j\, :\, y_{j} \in Y\setminus Q}\beta_{j}(\overline{\gamma}(y_{j})-\overline{\gamma}(y_{1}))=0. \label{pirx}
\end{align}

If $y_{1}<x_{1}$ then from the linear independence we obtain that $x_{j}=y_{j}$ for all $j=1,\ldots, \ell$, and $\lambda_{j}=\beta_{j}$ for all $j=2,\ldots, \ell$. In what follows we assume $y_{1}<x_{1}$. 

Notice that if for any $y \in Y \setminus Q$ we have $y\neq x_{1}$ then (\ref{pirx}) contradicts to the linear independence.  On the other hand if for some $y_{j^{*}}\in Y\setminus Q$  we have $y_{j^{*}}=x_{1}$ (we remark that there can be only one such $y_{j^{*}}$ in $Y\setminus Q$, moreover, $y_{j^{*}}\notin Q$) then (\ref{pirx}) we can rewrite as 
\begin{align} 
     &(1-\beta_{j}^{*}-\sum_{j=2}^{\ell}\lambda_{j})(\overline{\gamma}(x_{1})-\overline{\gamma}(y_{1}))+\sum_{j\, :\, x_{j} \in X\setminus Q}\lambda_{j}(\overline{\gamma}(x_{j})-\overline{\gamma}(y_{1}))+\nonumber\\
     &\sum_{j\, :\, x_{j} \in Q}(\lambda_{j}-\beta_{k_{j}})(\overline{\gamma}(x_{j})-\overline{\gamma}(y_{1}))
     -\sum_{j\, :\, y_{j} \in Y\setminus Q, \, y_{j}\neq y_{j^{*}}}\beta_{j}(\overline{\gamma}(y_{j})-\overline{\gamma}(y_{1}))=0. \label{pirx1}
\end{align}
Invoking the linear independence we must have $1-\beta_{j}^{*}-\sum_{j=2}^{\ell}\lambda_{j}=0$. Since $\lambda_{j}, \beta_{j} >0$ we have $X\setminus Q$ and  $Y\setminus (Q \cup\{y_{j^{*}}\})$ are empty. Then $Q$ has cardinality $\ell-1$ and $Q$ does not contain $y_{j^{*}}$ which is a contradiction.

\subsubsection{The proof of (\ref{welld})}

Assume $n=2\ell$. Since $\overline{U}_{n}$ and $\overline{L}_{n}$ are diffeomorphisms between $\mathrm{int}(\Delta_{c}^{\ell}\times \Delta_{*}^{\ell})$ 
and $\mathrm{int}(\mathrm{conv}(\overline{\gamma}([0,1])))$ we see that the equations 
\begin{align}
    &B^{\sup}(\overline{U}(t))=U^{z}(t), \label{be1}\\
    &B^{\inf}(\overline{L}(t))=L^{z}(t) \label{be2}
\end{align}
for all $t \in \mathrm{int}(\Delta_{c}^{\ell}\times \Delta_{*}^{\ell})$ define functions $B^{\sup}$ and $B^{\inf}$ uniquely on $\mathrm{int}(\mathrm{conv}(\overline{\gamma}([0,1])))$. We would like to extend the definitions of $B^{\sup}$ and $B^{\inf}$ to the boundary of $\mathrm{conv}(\overline{\gamma}([0,1]))$ just by taking $t \in \partial (\Delta_{c}^{\ell}\times \Delta_{*}^{\ell})$ in (\ref{be1}) and  (\ref{be2}). To make sure that the choice $t \in \partial (\Delta_{c}^{\ell}\times \Delta_{*}^{\ell})$  in (\ref{be1}) defines $B^{\sup}$ (and $B^{\inf}$) uniquely and continuously on  $\mathrm{conv}(\overline{\gamma}([0,1]))$ we shall verify the following
\begin{lemma}\label{gran1}
If $\overline{U}(t_{1})=\overline{U}(t_{2})$ for some $t_{1}, t_{2} \in  \Delta_{c}^{\ell}\times \Delta_{*}^{\ell}$, then $U^{z}(t_{1})=U^{z}(t_{2})$. Similarly, if $\overline{L}(t_{1})=\overline{L}(t_{2})$ for some $t_{1}, t_{2} \in  \Delta_{c}^{\ell}\times \Delta_{*}^{\ell}$, then $L^{z}(t_{1})=L^{z}(t_{2})$. 
\end{lemma}
\begin{proof}
Without loss of generality we can assume that $t_{1}, t_{2} \in \partial (\Delta_{c}^{\ell}\times \Delta_{*}^{\ell})$ otherwise the lemma follows from (\ref{b2l}), (\ref{diff2lu}), and (\ref{diff2ll}). 

First we show  $\overline{L}(t_{1})=\overline{L}(t_{2})$ for some $t_{1}, t_{2} \in  \partial (\Delta_{c}^{\ell}\times \Delta_{*}^{\ell})$ implies $L^{z}(t_{1})=L^{z}(t_{2})$. If $t_{1}=t_{2}$ there is nothing to prove, therefore, we assume $t_{1}\neq t_{2}$. For $t_{1} = (\lambda_{1}, \ldots, \lambda_{\ell}, x_{1}, \ldots, x_{\ell}) \in \partial (\Delta_{c}^{\ell}\times \Delta_{*}^{\ell})$ we have
\begin{align*}
    \overline{L}_{2\ell}(t_{1})  = \sum_{j=1}^{\ell} \lambda_{j}\overline{\gamma}(x_{j}).
\end{align*}
Among $\lambda_{1}, \ldots, \lambda_{\ell}$ many of them can be zero so we reduce the sum into $\sum_{j=1}^{\ell_{1}} \lambda_{q_{j}} \overline{\gamma}(x_{q_{j}})$ where $\lambda_{q_{j}}>0$, $\ell_{1}\leq \ell$, and $0\leq x_{q_{1}}\leq \ldots \leq x_{q_{\ell_{1}}}\leq 1$. Next, among $x_{q_{1}}, \ldots, x_{q_{\ell_{1}}}$  many can be equal to each other. Those $x_{q_{j}}$ who are equal to each other we group them together, and those $x_{j}$'s which are zero we remove from the sum by reducing the sum if necessary. This brings as to the following expression
\begin{align*}
\overline{L}_{2\ell}(t_{1})=\sum_{k=1}^{m}\lambda_{I_{k}} \overline{\gamma}(x_{I_{k}})
\end{align*}
where $ I_{k} \subset \{1, \ldots, \ell\}$,  the sets $I_{k}$ are disjoint for all $k=1, \ldots, m$. Here, $0<x_{I_{1}}<\ldots<x_{I_{m}}\leq 1$; for any $k, 1\leq k \leq m$ we have  $x_{j}=x_{I_{k}}$ for all $j \in I_{k}$; for any $k$, $1\leq k \leq m$ we set $0<\lambda_{I_{k}} := \sum_{j \in I_{k}} \lambda_{j}$. We remark that if $I_{k}=\emptyset$ then the term $\lambda_{I_{k}} \overline{\gamma}(x_{I_{k}})$ is zero by definition. 

Similarly, for $t_{2} = (\beta_{1}, \ldots, \beta_{\ell}, y_{1}, \ldots, y_{\ell}) \in \partial (\Delta_{c}^{\ell}\times \Delta_{*}^{\ell})$ we can write 
\begin{align*}
    \overline{L}_{2\ell}(t_{2})=\sum_{k=1}^{v}\beta_{J_{k}} \overline{\gamma}(y_{J_{k}})
\end{align*}
with $v \leq \ell$. 

As in the previous section, from linear independence of the vectors $\overline{\gamma}(z_{1}), \ldots, \overline{\gamma}(z_{2\ell})$, where $0<z_{1}<\ldots<z_{2\ell}\leq 1$, it follows that $ \overline{L}_{2\ell}(t_{1})= \overline{L}_{2\ell}(t_{2})$ holds if and only if $v=m$, $x_{I_{k}}=y_{J_{k}}$, and $\lambda_{I_{k}}=\beta_{J_{k}}$ for all $k=1, \ldots, m$. Hence $L^{z}_{2\ell}(t_{1})=L^{z}_{2\ell}(t_{2})$.

The proof for the map $\overline{U}_{2\ell}$ proceeds in the same way as for $\overline{L}_{2\ell}$. Indeed, the equality $\overline{U}_{2\ell}(t_{1})=\overline{U}_{2\ell}(t_{2})$ implies 
$\sum_{j=1}^{\ell} \lambda_{j}(\overline{\gamma}(x_{j})-\overline{\gamma}(1)) =   \sum_{j=1}^{\ell} \beta_{j}(\overline{\gamma}(y_{j})-\overline{\gamma}(1)).$
By removing zero terms, and grouping the similar terms inside the sums as before we obtain the equation 
\begin{align*}
\sum_{k=1}^{m} \lambda_{I_{k}}(\overline{\gamma}(x_{I_{k}})-\overline{\gamma}(1))=\sum_{k=1}^{v} \beta_{J_{k}}(\overline{\gamma}(y_{J_{k}})-\overline{\gamma}(1)),
\end{align*}
where we also removed the terms containing those $x_{j}$ and $y_{i}$ which are equal to $1$.
Applying Lemma~\ref{ltorsion} with $r=1$  
  we obtain that $v=m$ and $x_{I_{k}}=y_{J_{k}}$ for all $k=1,\ldots, m$, and $\lambda_{I_{k}}=\beta_{J_{k}}$. Hence $U^{z}_{2\ell}(t_{1})=U^{z}_{2\ell}(t_{2})$
\end{proof}

Next, we prove the analog of Lemma~\ref{gran1} for $n=2\ell-1$. 
\begin{lemma}\label{gran2}
If $\overline{U}(t_{1})=\overline{U}(t_{2})$ for some $t_{1}, t_{2} \in  \Delta_{c}^{\ell}\times \Delta_{*}^{\ell-1}$, then $U^{z}(t_{1})=U^{z}(t_{2})$. Similarly, if $\overline{L}(t_{1})=\overline{L}(t_{2})$ for some $t_{1}, t_{2} \in  \Delta_{c}^{\ell-1}\times \Delta_{*}^{\ell}$, then $L^{z}(t_{1})=L^{z}(t_{2})$. 
\end{lemma}
\begin{proof}
Without loss of generality we can assume that $t_{1}, t_{2} \in  \partial (\Delta_{c}^{\ell}\times \Delta_{*}^{\ell-1})$ (similarly,  $t_{1}, t_{2} \in  \partial (\Delta_{c}^{\ell-1}\times \Delta_{*}^{\ell})$ in the second claim of the lemma) otherwise the lemma follows from (\ref{b2l-1}), (\ref{diff2l-1u}), and (\ref{diff2l-1l}).

We show that the equality $\overline{U}(t_{1})=\overline{U}(t_{2})$ for some 
$t_{1}=(\lambda_{1}, \ldots, \lambda_{\ell}, x_{2}, \ldots, x_{\ell})$, and  $t_{2}=(\beta_{1}, \ldots, \beta_{\ell}, y_{2}, \ldots, y_{\ell})$ in $\partial (\Delta_{c}^{\ell}\times \Delta_{*}^{\ell-1})$ implies $L^{z}(t_{1})=L^{z}(t_{2})$. We can further assume  $t_{1}\neq t_{2}$ otherwise there is nothing to prove. We have
\begin{align}\label{jau1}
    \lambda_{1} \overline{\gamma}(1)+\sum_{j=2}^{\ell}\lambda_{j} \overline{\gamma}(x_{j}) = \beta_{1} \overline{\gamma}(1)+\sum_{j=2}^{\ell}\beta_{j}\overline{\gamma}(y_{j}).
\end{align}
As in the previous lemma, in the left hand side of (\ref{jau1}) we reduce the sum by removing those $\lambda_{j}$'s which are equal to zero. We further reduce the sum by considering only  positive $x_{j}$'s. Next, among the numbers $0\leq x_{2}\leq \ldots \leq x_{\ell}\leq 1$, those who are equal to each other we group them together, and those $x_{j}$'s which are equal to $1$ we group with $\lambda_{1} \overline{\gamma}(1)$. Eventually, the left hand side of (\ref{jau1}) takes the form $\lambda_{I_{0}}\overline{\gamma}(1)+\sum_{j=1}^{m}\lambda_{I_{j}}\overline{\gamma}(x_{j})$,  where $m\leq \ell-1$, $0<x_{I_{1}}<\ldots<x_{I_{m}}<1$, and $\lambda_{I_{j}} = \sum_{j \in I_{j}}\lambda_{j}$ with $\lambda_{I_{0}}\geq 0$ and $\lambda_{I_{j}}>0$ for all $j=1, \ldots, m$. Making a similar reduction in the right hand side of (\ref{jau1}), we see that (\ref{jau1}) takes the form 
\begin{align}\label{jau2}
    (\lambda_{I_{0}}-\beta_{J_{0}})\overline{\gamma}(1)+\sum_{j=1}^{m}\lambda_{I_{j}}\overline{\gamma}(x_{I_{j}}) - \sum_{j=1}^{v}\beta_{J_{j}}\overline{\gamma}(y_{J_{j}})=0.
\end{align}
Since $1+m+v\leq 2\ell-1$ it follows from Lemma~\ref{lltorsion} that (\ref{jau2}) holds if and only if $m=v$, $\lambda_{I_{j}}=\beta_{J_{j}}$ for all $j=2, \ldots, m$, and $x_{I_{j}}=y_{J_{j}}$ for all $j=1,\ldots, m$. It then follows that $U^{z}(t_{1})=U^{z}(t_{2})$. 

Next, we show that the equality $\overline{L}(t_{1})=\overline{L}(t_{2})$ for some 
$t_{1}=(\lambda_{2}, \ldots, \lambda_{\ell}, x_{1}, \ldots, x_{\ell})$, and  $t_{2}=(\beta_{2}, \ldots, \beta_{\ell}, y_{1}, \ldots, y_{\ell})$ in $\partial (\Delta_{c}^{\ell-1}\times \Delta_{*}^{\ell})$ implies $L^{z}(t_{1})=L^{z}(t_{2})$. Without loss of generality assume $t_{1} \neq t_{2}$ and $y_{1}\leq x_{1}$. The equality $\overline{L}(t_{1})=\overline{L}(t_{2})$ implies
\begin{align*}
    (1-\sum_{j=2}^{\ell}\lambda_{j})\overline{\gamma}(x_{1})+\sum_{j=2}^{\ell}\lambda_{j} \overline{\gamma}(x_{j})= (1-\sum_{j=2}^{\ell}\beta_{j})\overline{\gamma}(y_{1})+\sum_{j=2}^{\ell}\beta_{j} \overline{\gamma}(y_{j}),
\end{align*}
which we can rewrite as 
\begin{align}\label{sum19}
    (1-\sum_{j=2}^{\ell}\lambda_{j})(\overline{\gamma}(x_{1})-\overline{\gamma}(y_{1}))+\sum_{j=2}^{\ell}\lambda_{j}(\overline{\gamma}(x_{j})-\overline{\gamma}(y_{1}))-\sum_{j=2}^{\ell}\beta_{j}(\overline{\gamma}(y_{j})-\overline{\gamma}(y_{1}))=0.
\end{align}
We would like to show $L^{z}(t_{1})-L^{z}(t_{2})=0$. Notice that 
\begin{align} 
     &L^{z}(t_{1})-L^{z}(t_{2}) =(1-\sum_{j=2}^{\ell}\lambda_{j})(\gamma_{n+1}(x_{1})-\gamma_{n+1}(y_{1}))+\sum_{j=2}^{\ell}\lambda_{j}(\gamma_{n+1}(x_{j})-\gamma_{n+1}(y_{1}))-\label{sum13}\\
     &\sum_{j=2}^{\ell}\beta_{j}(\gamma_{n+1}(y_{j})-\gamma_{n+1}(y_{1})).  \nonumber
\end{align}
Rearranging and grouping equal terms in (\ref{sum19}) as in the previous arguments we can rewrite (\ref{sum19}) as  
\begin{align}
    &(1-\sum_{j=1}^{m_{1}} \lambda_{I^{1}_{j}} -\beta_{I_{0}})(\overline{\gamma}(x_{1})-\overline{\gamma}(y_{1})) + \sum_{j=1}^{m_{2}}\lambda_{I_{j}^{2}}(\overline{\gamma}(x_{I_{j}^{2}})-\overline{\gamma}(y_{1})) \nonumber\\
    &+\sum_{j=1}^{m_{3}}(\lambda_{I_{j}^{3}} - \beta_{J_{j}^{1}})(\overline{\gamma}(x_{I_{j}^{3}})-\overline{\gamma}(y_{1})) - \sum_{j=1}^{m_{4}} \beta_{J_{j}^{2}}(\overline{\gamma}(y_{J_{j}^{2}})-\overline{\gamma}(y_{1}))=0, \label{sum18}
\end{align}
where $m_{1}, m_{2}, m_{4}$ are non-negative integers with $1+m_{2}+m_{3}+m_{4}\leq 2\ell-1$ (if $m_{k}=0$ then the corresponding sum is set to be zero), $I_{j}^{1}, I_{j}^{2}, I_{j}^{3}, J_{j}^{1}, J_{j}^{2}$ are subsets of $\{2, \ldots, \ell\}$, $\beta_{I_{0}}\geq 0$, $\lambda_{I_{j}^{k}}=\sum_{j \in I_{j}^{k}} \lambda_{j}>0$, $\beta_{J_{j}^{k}}=\sum_{j \in J_{j}^{k}} \beta_{j}>0$, $\lambda_{I_{j}^{3}}\neq \beta_{J_{j}^{1}}$, and the points $x_{1}, \{x_{I_{j}^{2}}\}_{j=1}^{m_{2}},\{x_{I_{j}^{3}}\}_{j=1}^{m_{3}},\{y_{J_{j}^{2}}\}_{j=1}^{m_{4}}$ are different from each other, none of them (except of $x_{1}$) coincides with $y_{1}$, and all of them (except of $x_{1}$)  belong to $(0,1]$. We remark that $x_{1}$ can be equal to $y_{1}$.  In a similar way we can rewrite  (\ref{sum13}) as (\ref{sum18}), i.e., 
\begin{align*}
    &L^{z}(t_{1})-L^{z}(t_{2})=\\
    &(1-\sum_{j=1}^{m_{1}} \lambda_{I^{1}_{j}} -\beta_{I_{0}})(\gamma_{n+1}(x_{1})-\gamma_{n+1}(y_{1})) + \sum_{j=1}^{m_{2}}\lambda_{I_{j}^{2}}(\gamma_{n+1}(x_{I_{j}^{2}})-\gamma_{n+1}(y_{1})) \\
    &+\sum_{j=1}^{m_{3}}(\lambda_{I_{j}^{3}} - \beta_{J_{j}^{1}})(\gamma_{n+1}(x_{I_{j}^{3}})-\gamma_{n+1}(y_{1})) - \sum_{j=1}^{m_{4}} \beta_{J_{j}^{2}}(\gamma_{n+1}(y_{J_{j}^{2}})-\gamma_{n+1}(y_{1})). 
\end{align*}
 The next lemma  follows from Corollary~\ref{klasikac}.
\begin{lemma}\label{lgtorsion}
For any numbers $z_{j}$, $1\leq j \leq 2\ell-1$, such that  $0<z_{1}<z_{2}<\ldots<z_{2\ell}\leq  1$, and any $r \in [0,1]\setminus\{z_{1}, \ldots, z_{2\ell}\}$, the vectors $\overline{\gamma}(z_{1})-\overline{\gamma}(r), \ldots, \overline{\gamma}(z_{2\ell-1})-\overline{\gamma}(r)$ are linearly independent in $\mathbb{R}^{2\ell-1}$. 
\end{lemma} 
If $y_{1}=x_{1}$ then $L^{z}(t_{1})-L^{z}(t_{2})=0$ follows from (\ref{sum18}) and Lemma~\ref{lgtorsion}. If $y_{1}<x_{1}$, then applying Lemma~\ref{lgtorsion} to (\ref{sum18}) we see that $1-\sum_{j=1}^{m_{1}} \lambda_{I^{1}_{j}} -\beta_{I_{0}}=0$ and $m_{2}=m_{3}=m_{4}=0$, which implies that $L^{z}(t_{1})-L^{z}(t_{2})=0$.  Lemma~\ref{gran2} is proved.
\end{proof}

\subsubsection{The proof of (\ref{mincon1}) and (\ref{maxcon2})}
We start with (\ref{mincon1}). Assume $n=2\ell-1$. First we show that $B^{\sup}(\overline{\gamma})=\gamma_{n+1}$.  We remind that 
\begin{align*}
    B^{\sup}(\beta_{1} \overline{\gamma}(1)+\sum_{j=2}^{\ell}\beta_{j} \overline{\gamma}(x_{j}))=\beta_{1} \gamma_{n+1}(1)+\sum_{j=2}^{\ell}\beta_{j} \gamma_{n+1}(x_{j}), 
\end{align*}
holds for all $(\beta_{1}, \ldots, \beta_{\ell}, x_{2}, \ldots, x_{\ell}) \in \Delta_{c}^{\ell}\times \Delta_{*}^{\ell-1}$. We claim that if $\beta_{1} \overline{\gamma}(1)+\sum_{j=2}^{\ell}\beta_{j} \overline{\gamma}(x_{j}) = \overline{\gamma}(y)$ for some $y \in [0,1]$ then $\beta_{1} \gamma_{n+1}(1)+\sum_{j=2}^{\ell}\beta_{j} \gamma_{n+1}(x_{j})=\gamma_{n+1}(y)$. Indeed, $\overline{\gamma}(y)=\overline{U}(t_{2})$ with $t_{2} = (1,0, \ldots, 0, y) \in \Delta_{c}^{\ell}\times \Delta_{*}^{\ell-1}$, and $\beta_{1} \gamma_{n+1}(1)+\sum_{j=2}^{\ell}\beta_{j} \gamma_{n+1}(x_{j})=\overline{U}(t_{1})$ with $t_{1}=(\beta_{1}, \ldots, \beta_{\ell}, x_{1}, \ldots, x_{\ell}) \in \Delta_{c}^{\ell}\times \Delta_{*}^{\ell-1}$. Thus the claim follows from Lemma~\ref{gran2}. 

Next, we show that $B^{\sup}$ is concave on $\mathrm{conv}(\overline{\gamma}([0,1]))$. As the surface parametrized by $U_{n}(t)$, $t \in \Delta_{c}^{\ell}\times \Delta_{*}^{\ell-1}$, coincides with the graph $\{ (x,B^{\sup}(x)), x \in \mathrm{conv}(\overline{\gamma}([0,1]))\}$, and $B^{\sup} \in C(\mathrm{conv}(\overline{\gamma}([0,1])))$, it suffices to show that the tangent plane $T$ at   $U_{n}(s)$, for any $s=(\lambda_{1}, \ldots, \lambda_{\ell}, y_{2}, \ldots, y_{\ell})\in \mathrm{int}(\Delta_{c}^{\ell}\times \Delta_{*}^{\ell-1})$,  lies {\em above} the surface $U_{n}$. The equation of the tangent plane $T$  at $U(s):=U_{n}(s)$ is given as 
\begin{align*}
    T(x):=\det(U_{\beta_{1}}(s), \ldots, U_{\beta_{\ell}}(s), U_{x_{2}}(s), \ldots, U_{x_{\ell}}(s), x-U(s))=0, \quad x \in \mathbb{R}^{n+1}. 
\end{align*}
We have 
\begin{align*}
    T(x) = \lambda_{1} \cdots \lambda_{\ell}\det(\gamma(1), \gamma(y_{2}), \ldots, \gamma(y_{\ell}), \gamma'(y_{2}), \ldots, \gamma'(y_{\ell}), x).
\end{align*}
To show that the tangent plane $T$ lies {\em above} the surface, first we should find the sign of  $T(\lambda e_{n+1})$  as $\lambda \to \infty$, where $e_{n+1}=(\underbrace{0, \ldots, 0, 1}_{n+1})$. For sufficiently large positive $\lambda$ we have
\begin{align*}
    \mathrm{sign}(T(\lambda e_{n+1})) = \mathrm{sign}( \det(\overline{\gamma}(1), \overline{\gamma}(y_{2}), \ldots, \overline{\gamma}(y_{\ell}), \overline{\gamma}'(y_{2}), \ldots, \overline{\gamma}'(y_{\ell}))).
\end{align*}
On the other hand we have 
\begin{align*}
    &\det(\overline{\gamma}(1), \overline{\gamma}(y_{2}), \ldots, \overline{\gamma}(y_{\ell}), \overline{\gamma}'(y_{2}), \ldots, \overline{\gamma}'(y_{\ell})) = \\
    &(-1)^{\frac{(\ell-1)(\ell-2)}{2}} \det(\overline{\gamma}(y_{2}), \overline{\gamma}'(y_{2}), \ldots, \overline{\gamma}(y_{\ell}), \overline{\gamma}'(y_{\ell}), \overline{\gamma}(1))
     =\\
     &(-1)^{\frac{(\ell-1)(\ell-2)}{2}}\det(\overline{\gamma}(y_{2})-\overline{\gamma}(0), \overline{\gamma}'(y_{2}), \ldots, \overline{\gamma}(y_{\ell})-\overline{\gamma}(y_{\ell-1}), \overline{\gamma}'(y_{\ell}), \overline{\gamma}(1)-\overline{\gamma}(y_{\ell}))=\\
     &(-1)^{\frac{(\ell-1)(\ell-2)}{2}}  \int_{y_{\ell}}^{1} \int_{y_{\ell-1}}^{y_{\ell}}\ldots \int_{0}^{y_{2}}\det (\overline{\gamma}'(v_{2}), \overline{\gamma}'(y_{2}), \ldots,\overline{\gamma}'(v_{\ell}), \overline{\gamma}'(y_{\ell}), \overline{\gamma}'(v_{
     \ell+1}))dv_{2} \ldots dv_{\ell} dv_{\ell+1}.
\end{align*}
Thus, Lemma~\ref{klasika} applied to $\overline{\gamma}$ shows that  $\mathrm{sign}(T(\lambda e_{n+1}))$, for sufficiently large $\lambda$, coincides with $(-1)^{\frac{(\ell-1)(\ell-2)}{2}}$. Therefore, the surface $U(t)$ being {\em below} the tangent plane $T$ simply means that $(-1)^{\frac{(\ell-1)(\ell-2)}{2}}T(U(t))\leq 0$ for all $t=(\beta_{1}, \ldots, \beta_{\ell}, x_{2}, \ldots, x_{\ell}) \in \Delta_{c}^{\ell}\times \Delta_{*}^{\ell-1}$.  We have 
\begin{align*}
     T(U(t)) =  \sum_{j=2}^{\ell} \beta_{j} \det(\gamma(1), \gamma(y_{2}), \ldots, \gamma(y_{\ell}), \gamma'(y_{2}), \ldots, \gamma'(y_{\ell}), \gamma(x_{j}))\, \prod_{k=1}^{\ell}\lambda_{k}. 
\end{align*}
It suffices to verify that 
\begin{align}\label{nacili1}
    (-1)^{\frac{(\ell-1)(\ell-2)}{2}}\det(\gamma(1), \gamma(y_{2}), \ldots, \gamma(y_{\ell}), \gamma'(y_{2}), \ldots, \gamma'(y_{\ell}), \gamma(u)) \leq 0
\end{align}
for all $u \in [0,1]$. We have
\begin{align}
    &(-1)^{\frac{(\ell-1)(\ell-2)}{2}}\det(\gamma(1), \gamma(y_{2}), \ldots, \gamma(y_{\ell}), \gamma'(y_{2}), \ldots, \gamma'(y_{\ell}), \gamma(u)) \nonumber \\
    &=\det(\gamma(y_{2}), \gamma'(y_{2}), \ldots, \gamma(y_{\ell}), \gamma'(y_{\ell}), \gamma(1), \gamma(u)). \label{perm1}
\end{align}
If $u \in [y_{\ell},1]$, then 
\begin{align*}
  &\det(\gamma(y_{2}), \gamma'(y_{2}), \ldots, \gamma(y_{\ell}), \gamma'(y_{\ell}), \gamma(1), \gamma(u))=\\
  &-\det(\gamma(y_{2}), \gamma'(y_{2}), \ldots, \gamma(y_{\ell}), \gamma'(y_{\ell}),\gamma(u),\gamma(1))=\\
  &-\det(\gamma(y_{2})-\gamma(0), \gamma'(y_{2}), \ldots, \gamma(y_{\ell})-\gamma(y_{\ell-1}), \gamma'(y_{\ell}),\gamma(u)-\gamma(y_{\ell}),\gamma(1)-\gamma(u))=\\
  &-\int_{u}^{1}\int_{y_{\ell}}^{x_{j}}\int_{y_{\ell-1}}^{y_{\ell}}\ldots \int_{0}^{y_{2}}\det(\gamma'(v_{2}), \gamma'(y_{2}), \ldots, \gamma'(v_{\ell}), \gamma'(y_{\ell}),\gamma'(v_{\ell+1}),\gamma'(v_{\ell+2}))dv_{2}\ldots dv_{\ell}dv_{\ell+1}dv_{\ell+2}
\end{align*}
is non-positive  by Lemma~\ref{klasika}. 

If $u \in [0,y_{2}]$ we again use (\ref{perm1}). Next, we move the column $\gamma(u)$ to the left of the column $\gamma(y_{2})$. Notice that we will acquire the negative sign because passing the couples $\gamma(y_{i}), \gamma'(y_{i})$ does not change the sign of the determinant, the negative sign arises by passing $\gamma(1)$. Using the similar integral representation as before together with Lemma~\ref{klasika} we see that the inequality (\ref{nacili1}) holds true in the case $u \in [0, \gamma(y_{2})]$. The case $u \in [y_{i}, y_{i+1}]$ for some $i \in \{2, \ldots, \ell-1\}$, is similar to the previous case. Indeed, first we apply (\ref{perm1}), then we place the column $\gamma(u)$ between the columns $\gamma'(y_{i}), \gamma(y_{i+1})$ (thus we acquire the negative sign), we use the similar integral representation as before together with Lemma~\ref{klasika} to conclude that (\ref{nacili1}) holds true in this case too. This finishes the proof of concavity of $B^{\sup}$ on $\mathrm{conv}(\overline{\gamma}([0,1]))$. 

Next, we show that $B^{\sup}$ is the minimal concave function in a family of concave functions  $G$ on $\mathrm{conv}(\overline{\gamma}([0,1]))$ with the obstacle condition $G(\overline{\gamma}(s)) \geq \gamma_{n+1}(s)$ for all $s \in [0,1]$. Indeed, pick an arbitrary point $x \in \mathrm{conv}(\overline{\gamma}([0,1]))$. We would like to show $G(x) \geq B^{\sup}(x)$. There exists $(\lambda_{1}, \ldots, \lambda_{\ell}, y_{2}, \ldots, y_{\ell}) \in \Delta_{c}^{\ell}\times \Delta_{*}^{\ell-1}$ such that $x = \lambda_{1} \overline{\gamma}(1)+\sum_{j=2}^{\ell}\lambda_{j} \overline{\gamma}(y_{j}).$
Therefore 
\begin{align*}
    B^{\sup}(x) = \lambda_{1} \gamma_{n+1}(1)+\sum_{j=2}^{\ell}\lambda_{j}\gamma_{n+1}(y_{j})\leq \lambda_{1} G(\overline{\gamma}(1))+\sum_{j=2}^{\ell}\lambda_{j}G(\overline{\gamma}(y_{j})) \leq G(x).
\end{align*}

Next we consider $B^{\sup}$ when $n=2\ell$. We only check the concavity of $B^{\sup}$ because the remaining properties (minimality and the obstacle condition $B^{\sup}(\overline{\gamma})=\gamma_{n+1}$) are verified similarly as in the dimension $n=2\ell-1$.  The equation of the tangent plane $T$ at point 
\begin{align*}
    U(s) :=U_{n}(s) = \sum_{j=1}^{\ell} \beta_{j} \gamma(y_{j})+(1-\sum_{j=1}^{\ell}\beta_{j})\gamma(1),
\end{align*}
where  $s=(\beta_{1}, \ldots, \beta_{\ell}, y_{1}, \ldots, y_{\ell}) \in \mathrm{int}(\Delta_{c}^{\ell}\times \Delta_{*}^{\ell})$, is given as
\begin{align*}
    T(x) := \det(U_{\beta_{1}}, \ldots, U_{\beta_{\ell}}, U_{y_{1}}, \ldots, U_{y_{\ell}}, x-U(s))=0, \quad x \in \mathbb{R}^{n+1}.
\end{align*}
We have 
\begin{align*}
    \mathrm{sign}(T(x)) = \mathrm{sign}(\det(\gamma(y_{1})-\gamma(1), \ldots, \gamma(y_{\ell})-\gamma(1),\gamma'(y_{1}), \ldots, \gamma'(y_{\ell}), x-\gamma(1))). 
\end{align*}
Next, 
\begin{align*}
  \mathrm{sign}(T(\lambda e_{n+1})) =   \mathrm{sign}(\det(\overline{\gamma}(y_{1})-\overline{\gamma}(1), \ldots, \overline{\gamma}(y_{\ell})-\overline{\gamma}(1),\overline{\gamma}'(y_{1}), \ldots, \overline{\gamma}'(y_{\ell})))
\end{align*}
as $\lambda \to +\infty$. On the other hand we have 
\begin{align*}
    &\det(\overline{\gamma}(y_{1})-\overline{\gamma}(1), \ldots, \overline{\gamma}(y_{\ell})-\overline{\gamma}(1),\overline{\gamma}'(y_{1}), \ldots, \overline{\gamma}'(y_{\ell})) = \\
    &(-1)^{\ell}\det(\overline{\gamma}(y_{2})-\overline{\gamma}(y_{1}), \ldots, \overline{\gamma}(y_{\ell})-\overline{\gamma}(y_{\ell-1}),\overline{\gamma}(1)-\overline{\gamma}(y_{\ell}),\overline{\gamma}'(y_{1}), \ldots, \overline{\gamma}'(y_{\ell}))=\\
    &(-1)^{\frac{\ell(\ell-1)}{2}}\det(\overline{\gamma}'(y_{1}),\overline{\gamma}(y_{2})-\overline{\gamma}(y_{1}), \ldots, \overline{\gamma}'(y_{\ell}),\overline{\gamma}(1)-\overline{\gamma}(y_{\ell}))=\\
    &(-1)^{\frac{\ell(\ell-1)}{2}} \int_{y_{\ell}}^{1}\ldots \int_{y_{1}}^{y_{2}}  \det(\overline{\gamma}'(y_{1}),\overline{\gamma}'(x_{1}), \ldots, \overline{\gamma}'(y_{\ell}),\overline{\gamma}'(x_{\ell}))dx_{1} \ldots dx_{\ell}. 
\end{align*}
Thus, it follows from Lemma~\ref{klasika} that $\mathrm{sign}(T(\lambda e_{n+1})) = (-1)^{\frac{\ell(\ell-1)}{2}}$ as $\lambda \to \infty$. Therefore, to verify concavity of $B^{\sup}$ it suffices to show $(-1)^{\frac{\ell(\ell-1)}{2}} T(U(t)) \leq 0$ for all $t = (\lambda_{1}, \ldots, \lambda_{\ell}, x_{1}, \ldots, x_{\ell}) \in \Delta_{c}^{\ell}\times \Delta_{*}^{\ell}$. We have 
\begin{align*}
    T(U(t))=  \sum_{j=1}^{\ell} \lambda_{j} \det(\gamma(y_{1})-\gamma(1), \ldots, \gamma(y_{\ell})-\gamma(1),\gamma'(y_{1}), \ldots, \gamma'(y_{\ell}), \gamma(x_{j})-\gamma(1))\prod_{i=1}^{\ell}\beta_{i}.
\end{align*}
It suffices to show that $(-1)^{\frac{\ell(\ell-1)}{2}} \det(\gamma(y_{1})-\gamma(1), \ldots, \gamma(y_{\ell})-\gamma(1),\gamma'(y_{1}), \ldots, \gamma'(y_{\ell}), \gamma(u)-\gamma(1)) \leq 0$ for all $u \in [0,1]$. Assume $u \in [y_{i}, y_{i+1}]$ for some $i\in \{1,\ldots, \ell-1\}$. We have 
\begin{align*}
    &(-1)^{\frac{\ell(\ell-1)}{2}} \det(\gamma(y_{1})-\gamma(1), \ldots, \gamma(y_{\ell})-\gamma(1),\gamma'(y_{1}), \ldots, \gamma'(y_{\ell}), \gamma(u)-\gamma(1))=\\
    &-\det(\gamma'(y_{1}), \gamma(1)-\gamma(y_{1}), \ldots, \gamma'(y_{\ell}), \gamma(1)-\gamma(y_{\ell}), \gamma(1)-\gamma(u))=-\\
    &\det(\gamma'(y_{1}), \gamma(1)-\gamma(y_{1}), \ldots, \gamma'(y_{i}), \gamma(1)-\gamma(y_{i}),\gamma(1)-\gamma(u), \gamma'(y_{i+1}), \gamma(1)-\gamma(y_{i+1}), \ldots, \gamma(1)-\gamma(y_{\ell}))=\\
    &-\det(\gamma'(y_{1}), \gamma(y_{2})-\gamma(y_{1}), \ldots, \gamma'(y_{i}), \gamma(u)-\gamma(y_{i}),\gamma(y_{i+1})-\gamma(u), \gamma'(y_{i+1}), \gamma(y_{i+2})-\gamma(y_{i+1}), \ldots, \gamma(1)-\gamma(y_{\ell}))\\
    &=-  \int_{y_{\ell}}^{1}\ldots \int_{u}^{y_{i+1}} \int_{y_{i}}^{u} \ldots \int_{y_{1}}^{y_{2}}\det(\gamma'(y_{1}), \gamma'(v_{1}), \ldots, \gamma'(y_{i}), \gamma'(w), \gamma'(v_{i}), \gamma'(y_{i+1}), \ldots, \gamma'(v_{\ell}))  dv_{1}\ldots  dw dv_{i} \ldots dv_{\ell}
\end{align*}
which has a nonpositive sign by Lemma~\ref{klasika} (here $y_{i+2}$ for $i=\ell-1$ is set to be $1$). The cases $u \in [0,y_{1}]$, and $u \in [y_{\ell},1]$ are treated similarly.

Next, we verify (\ref{maxcon2}). The obstacle condition $B^{\inf}(\overline{\gamma})=\gamma_{n+1}$ and the minimality (assuming $B^{\inf}$ is convex) are verified similarly as in the case $B^{\sup}$. So, in what follows we only verify convexity of $B^{\inf}$.

Assume $n=2\ell-1$. The equation of the tangent plane $T$ at point
\begin{align*}
    L(s):=L_{n}(s) = (1-\sum_{j=2}^{\ell}\beta_{j})\gamma(y_{1})+\sum_{j=2}^{\ell}\beta_{j} \gamma(y_{j}),
\end{align*}
where $s=(\beta_{2}, \ldots, \beta_{\ell}, y_{1}, \ldots, y_{\ell}) \in \mathrm{int}(\Delta_{c}^{\ell-1}\times \Delta_{*}^{\ell})$ is given by 
\begin{align*}
    &T(x):=\det(L_{\beta_{2}}, \ldots, L_{\beta_{\ell}}, L_{y_{1}}, \ldots, L_{y_{\ell}}, x-L(s)) =\\
    &\det(\gamma(y_{2})-\gamma(y_{1}), \ldots, \gamma(y_{\ell})-\gamma(y_{1}), \gamma'(y_{1}), \ldots, \gamma'(y_{\ell}), x-\gamma(y_{1}))\, (1-\sum_{j=2}^{\ell}\beta_{j}) \prod_{j=2}^{\ell}\beta_{j}.
\end{align*}
We have 
\begin{align*}
    \mathrm{sign}(T(\lambda e_{n+1})) = \mathrm{sign}(\det(\overline{\gamma}(y_{2})-\overline{\gamma}(y_{1}), \ldots, \overline{\gamma}(y_{\ell})-\overline{\gamma}(y_{1}), \overline{\gamma}'(y_{1}), \ldots, \overline{\gamma}'(y_{\ell})))
\end{align*}
as $\lambda \to +\infty$. On the other hand
\begin{align*}
    &\det(\overline{\gamma}(y_{2})-\overline{\gamma}(y_{1}), \ldots, \overline{\gamma}(y_{\ell})-\overline{\gamma}(y_{1}), \overline{\gamma}'(y_{1}), \ldots, \overline{\gamma}'(y_{\ell})) =\\
    &(-1)^{\frac{\ell(\ell-1)}{2}}\det(\overline{\gamma}'(y_{1}),\overline{\gamma}(y_{2})-\overline{\gamma}(y_{1}), \ldots,\overline{\gamma}'(y_{\ell-1}), \overline{\gamma}(y_{\ell})-\overline{\gamma}(y_{\ell-1}),\overline{\gamma}'(y_{\ell}))=\\
    &(-1)^{\frac{\ell(\ell-1)}{2}} \int_{y_{\ell-1}}^{y_{\ell}}\ldots \int_{y_{1}}^{y_{2}} \det(\overline{\gamma}'(y_{1}),\overline{\gamma}'(v_{2}), \ldots,\overline{\gamma}'(y_{\ell-1}), \overline{\gamma}'(v_{\ell}),\overline{\gamma}'(y_{\ell})) dv_{2} \ldots dv_{\ell}.
\end{align*}
Thus $\mathrm{sign}(T(\lambda e_{n+1})) = (-1)^{\frac{\ell(\ell-1)}{2}}$ by Lemma~\ref{klasika} as $\lambda \to +\infty$. Therefore, $B^{\inf}$ is convex if $(-1)^{\frac{\ell(\ell-1)}{2}} T(L(t))\geq 0$ for all $t = (\lambda_{2}, \ldots, \lambda_{\ell}, x_{1}, \ldots, x_{\ell}) \in \Delta_{c}^{\ell-1}\times \Delta_{*}^{\ell}$. We have 
\begin{align*}
    &T(L(t)) = \det(\gamma(y_{2})-\gamma(y_{1}), \ldots, \gamma(y_{\ell})-\gamma(y_{1}), \gamma'(y_{1}), \ldots, \gamma'(y_{\ell}), L(t)-\gamma(y_{1}))\, (1-\sum_{j=2}^{\ell}\beta_{j}) \prod_{j=2}^{\ell}\beta_{j}=\\
    &(1-\sum_{k=2}^{\ell}\lambda_{k})\det(\gamma(y_{2})-\gamma(y_{1}), \ldots, \gamma(y_{\ell})-\gamma(y_{1}), \gamma'(y_{1}), \ldots, \gamma'(y_{\ell}), \gamma(x_{1})-\gamma(y_{1}))\, (1-\sum_{j=2}^{\ell}\beta_{j}) \prod_{j=2}^{\ell}\beta_{j}\\
    &+\sum_{k=2}^{\ell}\lambda_{k} \det(\gamma(y_{2})-\gamma(y_{1}), \ldots, \gamma(y_{\ell})-\gamma(y_{1}), \gamma'(y_{1}), \ldots, \gamma'(y_{\ell}), \gamma(x_{k})-\gamma(y_{1}))\, (1-\sum_{j=2}^{\ell}\beta_{j}) \prod_{j=2}^{\ell}\beta_{j}.
\end{align*}
Thus, to verify convexity of $B^{\inf}$, it suffices to show
\begin{align*}
    (-1)^{\frac{\ell(\ell-1)}{2}}\det(\gamma(y_{2})-\gamma(y_{1}), \ldots, \gamma(y_{\ell})-\gamma(y_{1}), \gamma'(y_{1}), \ldots, \gamma'(y_{\ell}), \gamma(u)-\gamma(y_{1}))\geq 0
\end{align*}
for all $u \in [0,1]$. Notice that 
\begin{align*}
     &(-1)^{\frac{\ell(\ell-1)}{2}}\det(\gamma(y_{2})-\gamma(y_{1}), \ldots, \gamma(y_{\ell})-\gamma(y_{1}), \gamma'(y_{1}), \ldots, \gamma'(y_{\ell}), \gamma(u)-\gamma(y_{1})) =\\
     & \det(\gamma'(y_{1}),\gamma(y_{2})-\gamma(y_{1}), \ldots, \gamma'(y_{\ell-1}),\gamma(y_{\ell})-\gamma(y_{1}), \gamma'(y_{\ell}), \gamma(u)-\gamma(y_{1})). 
\end{align*}
Next, assume $u \in [y_{i}, y_{i+1}]$ for some $i \in \{1, \ldots, \ell-1\}$ (the cases $u \in [0, y_{1}]$ and $u\in [y_{\ell},1]$ are considered similarly). We have 
\begin{align*}
&\det(\gamma'(y_{1}),\gamma(y_{2})-\gamma(y_{1}), \ldots, \gamma'(y_{\ell-1}),\gamma(y_{\ell})-\gamma(y_{1}), \gamma'(y_{\ell}), \gamma(u)-\gamma(y_{1}))=\\
&\det(\gamma'(y_{1}),\gamma(y_{2})-\gamma(y_{1}), \ldots, \gamma'(y_{i}), \gamma(u)-\gamma(y_{1}),\gamma(y_{i+1})-\gamma(y_{1}), \gamma'(y_{i+1}), \ldots, \gamma(y_{\ell})-\gamma(y_{1}), \gamma'(y_{\ell}))=\\
&\det(\gamma'(y_{1}),\gamma(y_{2})-\gamma(y_{1}), \ldots, \gamma'(y_{i}), \gamma(u)-\gamma(y_{i}),\gamma(y_{i+1})-\gamma(u), \gamma'(y_{i+1}), \ldots, \gamma(y_{\ell})-\gamma(y_{\ell-1}), \gamma'(y_{\ell}))\\
&=\int_{y_{\ell-1}}^{y_{\ell}}\ldots \int_{u}^{y_{i+1}}\int_{y_{i}}^{u}\ldots \int_{y_{1}}^{y_{2}}\\
&\det(\gamma'(y_{1}),\gamma'(v_{1}), \ldots, \gamma'(y_{i}), \gamma'(w),\gamma'(v_{i}), \gamma'(y_{i+1}), \ldots, \gamma'(v_{\ell-1}), \gamma'(y_{\ell})) dv_{1}\ldots dv_{i}dw\ldots dv_{\ell-1}.
\end{align*}
Thus $T(L(t))\geq 0$ by Lemma~\ref{klasika}. 

Next, we consider  $B^{\inf}$ when $n=2\ell$. As in the previous cases we only verify convexity of $B^{\inf}$ (minimality and the obstacle condition $B^{\inf}(\overline{\gamma})=\gamma_{n+1}$ are verified easily).   

The equation of the tangent plane $T$ at point 
\begin{align*}
    L(s):=L_{n}(s) = \sum_{j=1}^{\ell}\beta_{j} \gamma(y_{j}),
\end{align*}
where $s=(\beta_{1}, \ldots, \beta_{\ell}, y_{1}, \ldots, y_{\ell}) \in \mathrm{int}(\Delta_{c}^{\ell}\times \Delta_{*}^{\ell})$ is given by 
\begin{align*}
    &T(x):=\det(L_{\beta_{1}}, \ldots, L_{\beta_{\ell}}, L_{y_{1}}, \ldots, L_{y_{\ell}}, x-L(s)) =
    \det(\gamma(y_{1}), \ldots, \gamma(y_{\ell}), \gamma'(y_{1}), \ldots, \gamma'(y_{\ell}), x)\,  \prod_{j=1}^{\ell}\beta_{j}.
\end{align*}
We have 
\begin{align*}
    \mathrm{sign}(T(\lambda e_{n+1})) = \mathrm{sign}(\det(\overline{\gamma}(y_{1}), \ldots, \overline{\gamma}(y_{\ell}), \overline{\gamma}'(y_{1}), \ldots, \overline{\gamma}'(y_{\ell})))
\end{align*}
as $\lambda \to +\infty$. On the other hand
\begin{align}
    &\det(\overline{\gamma}(y_{1}), \ldots, \overline{\gamma}(y_{\ell}), \overline{\gamma}'(y_{1}), \ldots, \overline{\gamma}'(y_{\ell})) =\label{nishani1}\\
    &(-1)^{\frac{\ell(\ell-1)}{2}}\det(\overline{\gamma}(y_{1})-\overline{\gamma}(0),\overline{\gamma}'(y_{1}), \ldots, \overline{\gamma}(y_{\ell})-\overline{\gamma}(y_{\ell-1}),\overline{\gamma}'(y_{\ell}))= \nonumber\\
    &(-1)^{\frac{\ell(\ell-1)}{2}} \int_{y_{\ell-1}}^{y_{\ell}}\ldots \int_{0}^{y_{1}} \det(\overline{\gamma}'(v_{1}),\overline{\gamma}'(y_{1}), \ldots,\overline{\gamma}'(v_{\ell}), \overline{\gamma}'(y_{\ell})) dv_{1} \ldots dv_{\ell}. \nonumber
\end{align}
Thus $\mathrm{sign}(T(\lambda e_{n+1})) = (-1)^{\frac{\ell(\ell-1)}{2}}$ by Lemma~\ref{klasika} as $\lambda \to +\infty$. Therefore, $B^{\inf}$ is convex if $(-1)^{\frac{\ell(\ell-1)}{2}} T(L(t))\geq 0$ for all $t = (\lambda_{1}, \ldots, \lambda_{\ell}, x_{1}, \ldots, x_{\ell}) \in \Delta_{c}^{\ell}\times \Delta_{*}^{\ell}$. We have 
\begin{align*}
    &T(L(t)) = \sum_{k=1}^{\ell}\lambda_{k}\det(\gamma(y_{1}), \ldots, \gamma(y_{\ell}), \gamma'(y_{1}), \ldots, \gamma'(y_{\ell}), \gamma(x_{k}))\,  \prod_{j=1}^{\ell}\beta_{j}.
\end{align*}
Thus, to verify convexity of $B^{\inf}$, it suffices to show
\begin{align*}
    (-1)^{\frac{\ell(\ell-1)}{2}}\det(\gamma(y_{1}), \ldots, \gamma(y_{\ell}), \gamma'(y_{1}), \ldots, \gamma'(y_{\ell}), \gamma(u))\geq 0 \quad \text{for all} \quad u \in [0,1].
\end{align*}
Notice that 
\begin{align*}
(-1)^{\frac{\ell(\ell-1)}{2}}\det(\gamma(y_{1}), \ldots, \gamma(y_{\ell}), \gamma'(y_{1}), \ldots, \gamma'(y_{\ell}), \gamma(u)) = \det(\gamma(y_{1}),\gamma'(y_{1}), \ldots, \gamma(y_{\ell}), \gamma'(y_{\ell}), \gamma(u)). 
\end{align*}
Next, assume $u\in [y_{i}, y_{i+1}]$ for some $i \in \{1, \ldots, \ell-1\}$ (the cases $u\in [0,y_{1}]$ or $u \in [y_{\ell},1]$ are considered similarly). Set $y_{0}=0$.  We have 
\begin{align*}
&\det(\gamma(y_{1}),\gamma'(y_{1}), \ldots, \gamma(y_{\ell}), \gamma'(y_{\ell}), \gamma(u)) =\\
&\det(\gamma(y_{1}),\gamma'(y_{1}), \ldots, \gamma(y_{i}), \gamma'(y_{i}), \gamma(u), \gamma(y_{i+1}), \gamma'(y_{i+1}), \ldots)=\\
&\det(\gamma(y_{1})-\gamma(0),\gamma'(y_{1}), \ldots, \gamma(y_{i}) - \gamma(y_{i-1}), \gamma'(y_{i}), \gamma(u)-\gamma(y_{i}), \gamma(y_{i+1})-\gamma(u), \gamma'(y_{i+1}), \ldots)=\\
&\int_{y_{\ell-1}}^{y_{\ell}}\ldots \int_{u}^{y_{i-1}}\int_{y_{i}}^{u}\int_{y_{i-1}}^{y_{i}}\ldots \int_{0}^{y_{1}}\\
&\det(\gamma'(v_{1}),\gamma'(y_{1}), \ldots, \gamma'(v_{i}) , \gamma'(y_{i}), \gamma'(w), \gamma'(v_{i+1}), \gamma'(y_{i+1}), \ldots, \gamma'(v_{\ell}), \gamma'(y_{\ell})) dv_{1} \ldots dv_{i} dw dv_{i+1}\ldots dv_{\ell}.
\end{align*}
Thus $T(L(t))\geq 0$ by Lemma~\ref{klasika}.

\subsubsection{The proof of (\ref{giff})}\label{giffsub}
First we show the implication $B^{\sup}(u)=B^{\inf}(u) \Rightarrow u \in \partial\, \mathrm{conv}(\overline{\gamma}([0,1]))$. Consider the case $n=2\ell$.  Assume contrary, i.e., $u \in \mathrm{int}(\mathrm{conv}(\overline{\gamma}([0,1])))$.  Then using (\ref{diff2lu}), (\ref{diff2ll}) we can find  $t = (\lambda_{1}, \ldots, \lambda_{\ell}, x_{1}, \ldots, x_{\ell})$ and $s=(\beta_{1}, \ldots, \beta_{\ell}, y_{1}, \ldots, y_{\ell})$, both in $\mathrm{int}(\Delta_{c}^{\ell}\times \Delta_{*}^{\ell})$, such that 
\begin{align*}
    u=\sum_{j=1}^{\ell}\lambda_{j} \overline{\gamma}(x_{j})+(1-\sum_{j=1}^{\ell}\lambda_{j})\overline{\gamma}(1)=\sum_{j=1}^{\ell}\beta_{j}\overline{\gamma}(y_{j}).
\end{align*}
The equality $B^{\sup}(u)=B^{\inf}(u)$ implies (see (\ref{be1}), (\ref{be2}))
\begin{align*}
    \sum_{j=1}^{\ell}\lambda_{j} \gamma(x_{j})+(1-\sum_{j=1}^{\ell}\lambda_{j})\gamma(1)=\sum_{j=1}^{\ell}\beta_{j}\gamma(y_{j}).
\end{align*}
We see that $\gamma(1)$ is a linear combination of  $2\ell$ vectors $\gamma(x_{j}), \gamma(y_{j}), j=1,\ldots, \ell$ which leads us to a contradiction with Corollary~\ref{klasikac}. Thus $u \in \partial\, \mathrm{conv}(\overline{\gamma}([0,1]))$. 

Next, consider the case $n=2\ell-1$ and assume the contrary, i.e., $u \in \mathrm{int}(\mathrm{conv}(\overline{\gamma}([0,1])))$. Similarly as before we have 
\begin{align}\label{ukanaskneli}
    \lambda_{1} \gamma(1)+\sum_{j=2}^{\ell}\lambda_{j} \gamma(x_{j})=(1-\sum_{j=2}^{\ell}\beta_{j})\gamma(y_{1})+\sum_{j=2}^{\ell}\beta_{j} \gamma(y_{j})
\end{align}
for some  $t=(\lambda_{1}, \ldots, \lambda_{\ell}, x_{2}, \ldots, x_{\ell}) \in  \mathrm{int}(\Delta_{c}^{\ell}\times \Delta_{*}^{\ell-1})$ and $s = (\beta_{2}, \ldots, \beta_{\ell}, y_{1}, \ldots, y_{\ell}) \in  \mathrm{int}(\Delta_{c}^{\ell}\times \Delta_{*}^{\ell-1})$. The equality (\ref{ukanaskneli}) shows that $\gamma(1)$ is a linear combination of $2\ell-1$ vectors $\{\gamma(x_{j})\}_{j=2}^{\ell}$, $\{\gamma(y_{j})\}_{j=1}^{\ell}$ which contradicts to Corollary~\ref{klasikac}. 

Next we show the implication  $u \in \partial\, \mathrm{conv}(\overline{\gamma}([0,1])) \Rightarrow B^{\sup}(u)=B^{\inf}(u)$. Consider $n=2\ell$. Suppose 
\begin{align*}
    \overline{U}(t) \stackrel{\mathrm{def}}{=}\sum_{j=1}^{\ell}\lambda_{j} \overline{\gamma}(x_{j})+(1-\sum_{j=1}^{\ell}\lambda_{j})\overline{\gamma}(1)=\sum_{j=1}^{\ell}\beta_{j}\overline{\gamma}(y_{j}) \stackrel{\mathrm{def}}{=} \overline{L}(s)   
\end{align*}
for some $t = (\lambda, \ldots, \lambda_{\ell}, x_{1}, \ldots, x_{\ell})$ and $s=(\beta_{1}, \ldots, \beta_{\ell}, y_{1}, \ldots, y_{\ell})$, both in $\partial (\Delta_{c}^{\ell}\times \Delta_{*}^{\ell}).$ The goal is to show that 
\begin{align}\label{zgvari1}
    U^{z}(t) \stackrel{\mathrm{def}}{=} \sum_{j=1}^{\ell}\lambda_{j} \gamma_{n+1}(x_{j})+(1-\sum_{j=1}^{\ell}\lambda_{j})\gamma_{n+1}(1)=\sum_{j=1}^{\ell}\beta_{j}\gamma_{n+1}(y_{j})\stackrel{\mathrm{def}}{=} L^{z}(s).    
\end{align}
We claim that (\ref{zgvari1}) follows from the second part of Lemma~\ref{gran1}. For this it suffices to show that any point $\overline{U}(t)$, $t \in \partial (\Delta_{c}^{\ell}\times \Delta_{*}^{\ell})$, can be written as $\overline{L}(s_{1})$  for some $s_{1}=(\beta'_{1}, \ldots, \beta'_{\ell}, y'_{1}, \ldots, y'_{\ell}) \in \partial (\Delta_{c}^{\ell}\times \Delta_{*}^{\ell})$.  Indeed, as $t \in \partial (\Delta_{c}^{\ell}\times \Delta_{*}^{\ell})$ several cases can happen.  1) If $\sum_{j=1}^{\ell}\lambda_{j}=1$, then choose $\beta'_{j}=\lambda_{j}$, $j=1, \ldots, \ell-1$, $\beta'_{\ell}=1-\sum_{j=1}^{\ell-1}\lambda_{j}$, and $y'_{j}=x_{j}$, $j=1, \ldots, \ell$. Then 
\begin{align}\label{shemtxveva1}
L^{z}(s) = B^{\inf}(\overline{L}(s)) \stackrel{\mathrm{Lemma}~\ref{gran1}}{=}  B^{\inf}(\overline{L}(s_{1})) =\sum_{j=1}^{\ell}\beta'_{j}\gamma_{n+1}(y'_{j}) = U^{z}(t).
\end{align}
Next, 2) if at least one $\lambda_{j}=0$, say $\lambda_{p}=0$ for some $p \in \{1, \ldots, \ell\}$, then take $\beta'_{1}=\lambda_{1}, \ldots, \beta'_{p-1}=\lambda_{p-1}, \beta_{p}=\lambda_{p+1}, \ldots, \beta'_{\ell-1}=\lambda_{\ell}, \beta_{\ell}=1-\sum_{j=1}^{\ell} \lambda_{j}$, and $y'_{1}=x_{1}, \ldots, y'_{p-1}=x_{p-1}, y'_{p}=x_{p+1}, \ldots, y'_{\ell-1}=x_{\ell}, y'_{\ell}=1$ and repeat (\ref{shemtxveva1}). Next 3) if $x_{\ell}=1$, choose $(\beta'_{j}, y'_{j})=(\lambda_{j}, x_{j})$ for $j=1, \ldots, \ell-1$, and $(\beta'_{\ell}, y'_{\ell})=(\sum_{j=1}^{\ell-1}\lambda_{j},1)$ and repeat (\ref{shemtxveva1}). 4) If $x_{p}=x_{p+1}$ for some $p\in \{1, \ldots, \ell-1\}$ then take $y'_{j}=x_{j}$ for $j=1, \ldots, p$; $y'_{j}=x_{j+1}$ for $j=p+1, \ldots, \ell-1$; $y'_{\ell}=1$; $\beta'_{1}=\lambda_{1}$, \ldots, $\beta'_{p}=\lambda_{p}+\lambda_{p+1}$,   $\beta'_{p+1}=\lambda_{p+2}, \ldots, \beta'_{\ell-1}=\lambda_{\ell}$, $\beta'_{\ell}=1-\sum_{j=1}^{\ell}\lambda_{j}$ and repeat (\ref{shemtxveva1}). Finally, 5) if $x_{1}=0$ choose $\beta'_{j}=\lambda_{j+1}$, $j=1,\ldots, \ell-1$; $\beta'_{\ell}=1-\sum_{j=1}^{\ell}\lambda_{j}$; $y'_{j}=x_{j+1}$, $j=1, \ldots, \ell-1$; $y'_{\ell}=1$, and apply (\ref{shemtxveva1}). 

Next, consider $n=2\ell-1$.  Suppose 
\begin{align*}
   \overline{U}(t) \stackrel{\mathrm{def}}{=}   \sum_{j=2}^{\ell}\beta_{j} \overline{\gamma}(x_{j}) + \beta_{1} \overline{\gamma}(1)=(1-\sum_{j=2}^{\ell}\lambda_{j})\overline{\gamma}(y_{1})+\sum_{j=2}^{\ell}\lambda_{j} \overline{\gamma}(y_{j})  \stackrel{\mathrm{def}}{=} \overline{L}(s)
\end{align*}
for some $t=(\beta_{1}, \ldots, \beta_{\ell},x_{2}, \ldots, x_{\ell}) \in \partial (\Delta_{c}^{\ell}\times \Delta_{*}^{\ell-1})$ and $s=(\lambda_{2}, \ldots, \lambda _{\ell}, y_{1}, \ldots, y_{\ell}) \in \partial (\Delta_{c}^{\ell-1}\times \Delta_{*}^{\ell})$. We would like to show 
\begin{align}\label{shemxtveva2}
       U^{z}(t) \stackrel{\mathrm{def}}{=}  \sum_{j=2}^{\ell}\beta_{j} \gamma_{n+1}(x_{j})+ \beta_{1} \gamma_{n+1}(1)=(1-\sum_{j=2}^{\ell}\lambda_{j})\gamma_{n+1}(y_{1})+\sum_{j=2}^{\ell}\lambda_{j} \gamma_{n+1}(y_{j})  \stackrel{\mathrm{def}}{=} L^{z}(s).
\end{align}
As in the case $n=2\ell-1$ we claim that (\ref{shemxtveva2}) follows from Lemma~\ref{gran2}. It suffices to show that for any point $\overline{U}(t)$, $t \in \partial (\Delta_{c}^{\ell}\times \Delta_{*}^{\ell-1})$, there exists a point $s_{1} =(\lambda'_{2}, \ldots, \lambda'_{\ell}, y'_{1}, \ldots, \lambda'_{\ell}) \in \partial(\Delta_{c}^{\ell-1}\times \Delta_{*}^{\ell})$ such that $\overline{U}(t)=\overline{L}(s_{1})$. Several instances may happen.  1) if $\sum_{j=1}^{\ell}\beta_{j}=1$. Let 
\begin{align*}
(\lambda'_{2}, \ldots, \lambda'_{\ell-1}, \lambda'_{\ell}, y'_{1}, \ldots, y'_{\ell-1}, y'_{\ell})=(\beta_{3}, \ldots, \beta_{\ell}, \beta_{1}, x_{2}, \ldots, x_{\ell}, 1).
\end{align*}
Notice that $1-\sum_{j=2}^{\ell}\lambda'_{j}=\beta_{2}$. 2) if $\beta_{p}=0$ for some $p \in \{1,\ldots, \ell-1\}$ then let 
\begin{align*}
&(\lambda'_{2}, \ldots, \lambda'_{p-1}, \lambda'_{p}, \ldots, \lambda'_{\ell-1}, \lambda'_{\ell}, y'_{1}, y'_{2}, \ldots, y'_{p-1}, y'_{p}, \ldots, y'_{\ell-1}, y'_{\ell}) =\\
&(\beta_{2}, \ldots, \beta_{p-1}, \beta_{p+1}, \ldots, \beta_{\ell}, \beta_{1}, 0, x_{2}, \ldots, x_{p-1}, x_{p+1}, \ldots, x_{\ell}, 1).
\end{align*}
3) if $\beta_{1}=0$ then we choose $y'_{1}=0$ and 
\begin{align*}
(\lambda'_{2}, \ldots, \lambda'_{\ell}, y'_{2}, \ldots, y'_{\ell})=(\beta_{2}, \ldots, \beta_{\ell}, x_{2}, \ldots, x_{\ell}).
\end{align*}
4) if $x_{2}=0$ then we choose $y_{1}=0$ and 
\begin{align*}
(\lambda'_{2}, \ldots, \lambda'_{\ell-1}, \lambda'_{\ell}, y'_{2}, \ldots, y'_{\ell-1}, y'_{\ell})=(\beta_{3}, \ldots, \beta_{\ell}, \beta_{1}, x_{3}, \ldots, x_{\ell}, 1).
\end{align*}
5) if $x_{\ell}=1$ then let  $y_{1}=0$ and 
\begin{align*}
(\lambda'_{2}, \ldots, \lambda'_{\ell-1}, \lambda'_{\ell}, y'_{2}, \ldots, y'_{\ell-1}, y'_{\ell})=(\beta_{2}, \ldots, \beta_{\ell-1}, \beta_{\ell}+\beta_{1}, x_{2}, \ldots, x_{\ell-1}, 1).
\end{align*}
Finally, 6) if $x_{p}=x_{p+1}$ for some $p \in \{2, \ldots, \ell-1\}$ take $y_{1}=0$ and 
\begin{align*}
&(\lambda'_{2}, \ldots, \lambda'_{p-1},\lambda'_{p}, \lambda'_{p+1}, \ldots,  \lambda'_{\ell-1}, \lambda'_{\ell}, y'_{2}, \ldots, y'_{p}, y'_{p+1}, \ldots,  y'_{\ell-1}, y'_{\ell})=\\
&(\beta_{2}, \ldots, \beta_{p-1},\beta_{p}+\beta_{p+1}, \beta_{p+2},\ldots, \beta_{\ell}, \beta_{1}, x_{2}, \ldots, x_{p}, x_{p+2}, \ldots x_{\ell}, 1).
\end{align*}
Under such choices we have 
\begin{align*}
L^{z}(s) = B^{\inf}(\overline{L}(s)) \stackrel{\mathrm{Lemma}~\ref{gran2}}{=}  B^{\inf}(\overline{L}(s_{1})) =(1-\sum_{j=2}^{\ell}\lambda'_{j})\gamma_{n+1}(y'_{1})+\sum_{j=2}^{\ell}\lambda'_{j} \gamma_{n+1}(y'_{j})  = U^{z}(t).
\end{align*}
This finishes the proof of (\ref{giff}).

\subsubsection{The proof of (\ref{union})} \label{unionsub}
The inclusion 
\begin{align*}
 \{(x,B^{\mathrm{sup}}(x)), x \in \mathrm{conv}(\bar{\gamma}([0,1]))\} \cup \{(x,B^{\mathrm{inf}}(x)), x \in \mathrm{conv}(\bar{\gamma}([0,1]))\} \subset \partial\,  \mathrm{conv}(\gamma([0,1])) 
\end{align*}
is trivial. Indeed, it follows from (\ref{be1}) that the point $(x,B^{\sup}(x))$ is a convex combination of some points of $\gamma([0,1])$, therefore, $(x,B^{\sup}(x)) \in  \mathrm{conv}(\gamma([0,1]))$. On the other hand, no point of the form $(x,s)$, where $s>B^{\sup}(x)$ belongs to $\mathrm{conv}(\gamma([0,1]))$. Indeed, otherwise  $(x,s) = \sum_{j=1}^{m}\lambda_{j} \gamma(t_{j})$ for some $t_{j} \in [0,1]$ and nonnegative $\lambda_{j}$ such that $\sum_{j=1}^{m} \lambda_{j}=1$. Then  
\begin{align*}
    B^{\sup}(x)=B^{\sup}\left( \sum \lambda_{j} \overline{\gamma}(t_{j})\right)\stackrel{(\ref{mincon1})}{\geq} \sum \lambda_{j} B^{\sup}(\overline{\gamma}(t_{j})) \stackrel{(\ref{mincon1})}{=} \sum \lambda_{j} \gamma_{n+1}(t_{j})=s
\end{align*}
gives a contradiction. Thus $(x,B^{\sup}(x)) \in  \partial\,  \mathrm{conv}(\gamma([0,1]))$. In a similar way we have $(x,B^{\inf}(x)) \in  \partial\,  \mathrm{conv}(\gamma([0,1]))$ for $x \in \mathrm{conv}(\bar{\gamma}([0,1]))$. 

To verify the inclusion
\begin{align*}
    \partial\,  \mathrm{conv}(\gamma([0,1])) \subset \{(x,B^{\mathrm{sup}}(x)), x \in \mathrm{conv}(\bar{\gamma}([0,1]))\} \cup \{(x,B^{\mathrm{inf}}(x)), x \in \mathrm{conv}(\bar{\gamma}([0,1]))\}
\end{align*}
we pick a point  $(x,t) \in \partial\,  \mathrm{conv}(\gamma([0,1]))$ where $x \in \mathbb{R}^{n}$, i.e., $x \in \mathrm{conv}(\bar{\gamma}([0,1]))$. Clearly $B^{\inf}(x) \leq t \leq B^{\sup}(x)$. Assume contrary that $B^{\inf}(x) < t < B^{\sup}(x)$.  If $x \in \partial\, \mathrm{conv}(\bar{\gamma}([0,1]))$ then by (\ref{giff}) we have $B^{\inf}(x) = B^{\sup}(x)$, therefore, we get a contradiction. If $x \in \mathrm{int}(\mathrm{conv}(\bar{\gamma}([0,1])))$ then  (\ref{giff}) and continuity of $B^{\sup}$ and  $B^{\inf}$ imply that there exists a ball $U_{\varepsilon}(x)$ or radius $\varepsilon>0$ centered at point $x$ such that $U_{\varepsilon}(x) \subset \mathrm{int}(\mathrm{conv}(\bar{\gamma}([0,1])))$ and $B^{\inf}(s)<t-\delta < t+\delta<B^{\sup}(s)$ for all $s \in U_{\varepsilon}(x)$ and some $\delta>0$. Then 
\begin{align*}
    &(x,t) \in U_{\min\{\varepsilon, \delta\}}((x,t))\subset \{(s,y)\, :\, B^{\inf}(s)\leq y\leq B^{\sup}(s), s \in U_{\min\{\varepsilon, \delta\}}(x)\}=\\
    &\mathrm{conv}(\{(s,B^{\inf}(s)), \, s \in U_{\min\{\varepsilon, \delta\}}(x)\} \cup \{(s,B^{\sup}(s)), \, s \in U_{\min\{\varepsilon, \delta\}}(x)\}) \subset \mathrm{conv} (\gamma([0,1])),
\end{align*}
 where $U_{\min\{\varepsilon, \delta\}}((x,t))$ is the ball in $\mathbb{R}^{n+1}$ centered at $(x,t)$ with radius $\min\{\varepsilon, \delta\}$. We obtain a contradiction with the assumption that $(x,t)$ belongs to the boundary of $\mathrm{conv} (\gamma([0,1]))$.

The proof of Theorem~\ref{mth010} is complete. 
\end{proof}
\subsection{The proof of Proposition~\ref{sensitive}}

Take $\gamma(t) = (t,t^{4}, -t^{3})$ on $[-1,1]$. We have 
\begin{align*}
    (\gamma', \gamma'', \gamma''') = \begin{pmatrix}
    1 & 0 & 0 \\
    4t^{3} & 12t^{2} & 24 t\\
    -3t^{2} & -6t & -6
    \end{pmatrix}.
\end{align*}
All the leading principal minors of the matrix $(\gamma', \gamma'', \gamma''')$ are positive on $[-1,1]\setminus \{0\}$, and we notice that $2\times 2$ and $3\times 3$ the leading principal minors vanish at $t=0$. Assume contrary to Proposition~\ref{sensitive} that the map $B^{\sup}(x,y)$ defined on $\mathrm{conv}(\overline{\gamma}([-1,1]))$  by (\ref{vog}) is concave. We have 
\begin{align}
    B(\lambda (a,a^{4})+(1-\lambda)(1,1)) = -\lambda a^{3} -(1-\lambda), \lambda \in [0,1], a \in (-1,1).
\end{align}
In particular, $g(y):=B(0,y), y \in [0,1],$ must be concave. The restriction $\lambda a + (1-\lambda)=0$ implies $\lambda = \frac{1}{1-a}$. Therefore 
\begin{align*}
    \lambda a^{4} + (1-\lambda) = a^{3}+a^{2}+a \quad \text{and} \quad -\lambda a^{3} -(1-\lambda) = a^{2}+a.
\end{align*}
Since $-a^{3}-a^{2}-a = y \in [0,1]$ we must have $a \in [-1,0]$. Thus $g(-a^{3}-a^{2}-a) = a^{2}+a$ for $a \in [-1,0]$. differentiating both sides in $a$ two times we obtain 
\begin{align*}
    &g'(-a^{3}-a^{2}-a) = -\frac{2a+1}{3a^{2}+2a+1},\\
    &g''(-a^{3}-a^{2}-a) = \frac{-6a(a+1)}{(3a^{2}+2a+1)^{3}} > 0 \quad \text{for} \quad a \in [-1,0).
\end{align*}
Thus $g''>0$ gives a contradiction. 

\subsection{The proof of Theorem~\ref{mth1}}
We verify (\ref{extr01}). The verification of (\ref{extr02}) is similar. Denote 
\begin{align*}
    M^{\sup}(x) :=\sup_{a\leq Y\leq b} \{ \mathbb{E} \gamma_{n+1}(Y)\, :\, \mathbb{E} \overline{\gamma}(Y)=x\}, \quad x \in \mathrm{conv}(\overline{\gamma}([a,b])).
\end{align*}

First we show the inequality $M^{\sup} \leq B^{\sup}$ on $\mathrm{conv}(\overline{\gamma}([a,b]))$. Indeed, let $x \in \mathrm{conv}(\overline{\gamma}([a,b]))$.  Pick an arbitrary random variable $Y$ with values in $[a,b]$, such that $\mathbb{E} \overline{\gamma}(Y)=x$.  Then 
\begin{align*}
    \mathbb{E} \gamma_{n+1}(Y) \stackrel{(\ref{mincon1})}{=} \mathbb{E}B^{\sup}(\overline{\gamma}(Y)) \stackrel{(\ref{mincon1})+\mathrm{Jensen}}{\leq} B^{\sup}(\mathbb{E}\overline{\gamma}(Y))=B^{\sup}(x).
\end{align*}
Taking the supremum over all $Y$, $a\leq Y\leq b$, such that $\mathbb{E} \overline{\gamma}(Y)=x$, gives the inequality $M^{\sup}(x) \leq B^{\sup}(x)$. 

To verify the reverse inequality $M^{\sup}(x) \geq B^{\sup}(x)$ it suffices to construct at least one random variable $Y=Y(x)$, $a\leq Y\leq b$,  such that $\mathbb{E} \overline{\gamma}(Y)=x$ and $\mathbb{E} \gamma_{n+1}(Y)=B^{\sup}(x)$. Notice that $Y=\zeta(x)$, where $\zeta(x)$ is defined in Theorem~\ref{mth1}, satisfies $a\leq \zeta(x) \leq b$, $\mathbb{E} \overline{\gamma}(\zeta(x))=x$. It also follows from (\ref{vog})  that   $\mathbb{E} \gamma_{n+1}(\zeta(x))=B^{\sup}(x)$. 
 
 \subsection{The proof of Corollary~\ref{nobel2}}. The moment curve $\gamma$ has totally positive torsion on $[0,1]$, hence, Theorem~\ref{mth010} applies. 
 
 First we work with $B^{\sup}(x)=x_{n+1}$. Consider the case $n=2\ell$. By Theorem~\ref{mth010} there exists a unique point $(\lambda_{1}, \ldots, \lambda_{\ell}, y_{1}, \ldots, y_{\ell}) \in \mathrm{int}(\Delta_{c}^{\ell}\times \Delta_{*}^{\ell})$ such that $\sum_{j=1}^{\ell}\lambda_{j} \overline{\gamma}(y_{j})+(1-\sum_{j=1}^{\ell}\lambda_{j})\overline{\gamma}(1)=x$ then the value $x_{n+1}:=B^{\sup}(x)$ equals to $\sum_{j=1}^{\ell}\lambda_{j} y_{j}^{2\ell+1}+(1-\sum_{j=1}^{\ell}\lambda_{j})$. We would like to show that the linear equation  
 \begin{align}\label{gant01}
     \det 
 \begin{pmatrix}a_{0} & a_{1} & \ldots & a_{\ell}\\
 \vdots & & & \\
 a_{\ell} & a_{\ell+1} & \ldots & a_{2\ell}\end{pmatrix}=0,
 \end{align}
 where $a_{k}:=x_{k}-x_{k+1}$, $k=0, \ldots, 2\ell$, $x_{0}:=1$, has a unique solution in $x_{n+1}$ which equals to $\sum_{j=1}^{\ell}\lambda_{j} y_{j}^{2\ell+1}+(1-\sum_{j=1}^{\ell}\lambda_{j})$. First we check why  $x_{n+1}=\sum_{j=1}^{\ell}\lambda_{j} y_{j}^{2\ell+1}+(1-\sum_{j=1}^{\ell}\lambda_{j})$ solves (\ref{gant01}). Notice that $a_{k} = \langle y^{k}, \beta \rangle$, where $y^{k} := (y_{1}^{k}, \ldots, y_{\ell}^{k})$, and $\beta := (\lambda_{1}(1-y_{1}), \ldots, \lambda_{\ell}(1-y_{\ell}))$. The $j$'th column of the matrix in (\ref{gant01}), call it $w_{j}$, $j=0, \ldots, \ell$, we can write as $w_{j} = AD^{j}\beta^{T}$, where $A$ is $(\ell+1)\times \ell$ matrix with $m$'th column $(1, y_{m}, \ldots, y_{m}^{\ell})^{T}$, and $D$ is  $\ell\times\ell$  diagonal matrix with  diagonal entries $y_{1}, \ldots, y_{\ell}$. Since there exists a nonzero vector $(z_{0}, \ldots, z_{\ell})\in \mathbb{R}^{\ell+1}$ such that $z_{0}D^{0}+\ldots+z_{\ell}D^{\ell}=0$ (the number of variables $z_{j}$ is greater than the number of equations, i.e., $\ell$), it follows that the vectors $\{w_{0}, \ldots, w_{\ell}\}$ are linearly dependent, so (\ref{gant01}) holds true. 
 
 To show the uniqueness of the solution $x_{n+1}$ it suffices to show that the leading  $\ell\times \ell$ principal minor $R$ of the matrix in $(\ref{gant01})$ has nonzero determinant. Notice that  $R=\det(\tilde{w}_{0}, \ldots, \tilde{w}_{\ell-1})$, where $\tilde{w}_{j}=\tilde{A}D^{j}\beta^{T}$  and $\tilde{A}$ is obtained from $A$ by removing the last row. Assume contrary that $R=0$. Then there exists a nonzero vector $(z_{0}, \ldots, z_{\ell-1})\in \mathbb{R}^{\ell}$ such that $\tilde{A}(z_{0}D^{0}+\ldots+z_{\ell-1}D^{\ell-1})\beta^{T}=0$. As $\det(\tilde{A})\neq 0$ (Vandermonde matrix) we have $(z_{0}D^{0}+\ldots+z_{\ell-1}D^{\ell-1})\beta^{T}=0$. Since the entries of $\beta^{T}$ are nonzero and the matrix $(z_{0}D^{0}+\ldots+z_{\ell-1}D^{\ell-1})$ is diagonal we must have $z_{0}D^{0}+\ldots+z_{\ell-1}D^{\ell-1}=0$. The last equation rewrites as $\tilde{A}^{T}z^{T}=0$ where $z=(z_{0}, \ldots, z_{\ell-1})\neq 0$, which is a contradiction. 
 
 Next, consider $n=2\ell-1$. In this case $x = (1-\sum_{j=1}^{\ell}\lambda_{j})\overline{\gamma}(0)+\sum_{j=2}^{\ell}\lambda_{j}\overline{\gamma}(y_{j})+\lambda_{1}\overline{\gamma}(1)$ for a unique $(\lambda_{1}, \ldots, \lambda_{\ell}, y_{2}, \ldots, y_{\ell}) \in \mathrm{int}(\mathrm{conv}(\overline{\gamma}([0,1])))$, and the value $x_{n+1}:=B^{\sup}(x)$ is $(1-\sum_{j=1}^{\ell}\lambda_{j})\gamma_{n+1}(0)+\sum_{j=2}^{\ell}\lambda_{j}\gamma_{n+1}(y_{j})+\lambda_{1}\gamma_{n+1}(1)$. Set $b_{k}:=x_{k}-x_{k+1}$, $k=1, \ldots, 2\ell-1$. As before we would like to show that the linear equation 
 \begin{align}\label{gant02}
     \det 
 \begin{pmatrix}b_{1} & b_{2} & \ldots & b_{\ell}\\
 \vdots & & & \\
 b_{\ell} & b_{\ell+1} & \ldots & b_{2\ell-1}\end{pmatrix}=0,
 \end{align}
 has a unique solution in $x_{n+1}$ which equals to $(1-\sum_{j=1}^{\ell}\lambda_{j})\gamma_{n+1}(0)+\sum_{j=2}^{\ell}\lambda_{j}\gamma_{n+1}(y_{j})+\lambda_{1}\gamma_{n+1}(1)$. To check that such a choice for $x_{n+1}$ solves (\ref{gant02}), notice that $b_{k} =\langle y^{k}, \beta \rangle$, where $y^{k} = (y_{2}^{k}, \ldots, y_{\ell}^{k})$ and $\beta = (\lambda_{2}(1-y_{2}), \ldots, \lambda_{\ell}(1-y_{\ell}))$. 
 
 The $j$'th column of the matrix in (\ref{gant02}), call it $w_{j}$, $j=1, \ldots, \ell$, we can write as $w_{j} = AD^{j}\beta^{T}$, where $A$ is $\ell\times(\ell-1)$ matrix with $m$'th column $(y_{m}, \ldots, y_{m}^{\ell})^{T}$, $m=2, \ldots, \ell$, and $D$ is  $(\ell-1)\times(\ell-1)$  diagonal matrix with  diagonal entries $y_{2}, \ldots, y_{\ell}$. Since there exists a nonzero vector $(z_{1}, \ldots, z_{\ell})\in \mathbb{R}^{\ell}$ such that $z_{1}D+\ldots+z_{\ell}D^{\ell}=0$ (the number of variables $z_{j}$ is greater than the number of equations, i.e., $\ell-1$), it follows that the vectors $\{w_{1}, \ldots, w_{\ell}\}$ are linearly dependent, so (\ref{gant02}) holds true. 
 
 To show the uniqueness of the solution $x_{n+1}$ it suffices to show that the leading  $(\ell-1)\times (\ell-1)$  principal minor $R$ of the matrix in $(\ref{gant02})$ has nonzero determinant. Notice that  $R=\det(\tilde{w}_{1}, \ldots, \tilde{w}_{\ell-1})$, where $\tilde{w}_{j}=\tilde{A}D^{j}\beta^{T}$,  and $\tilde{A}$ is obtained from $A$ by removing the last row. Assume contrary that $R=0$. Then there exists nonzero vector $(z_{1}, \ldots, z_{\ell-1})\in \mathbb{R}^{\ell-1}$ such that $\tilde{A}(z_{1}D+\ldots+z_{\ell-1}D^{\ell-1})\beta^{T}=0$. As $\det(\tilde{A})\neq 0$ (Vandermonde matrix) we have $(z_{1}D+\ldots+z_{\ell-1}D^{\ell-1})\beta^{T}=0$. Since the entries of $\beta^{T}$ are nonzero and the matrix $(z_{1}D+\ldots+z_{\ell-1}D^{\ell-1})$ is diagonal we must have $z_{1}D+\ldots+z_{\ell-1}D^{\ell-1}=0$. The last equation rewrites as $\tilde{A}^{T}z^{T}=0$ where $z=(z_{1}, \ldots, z_{\ell-1})\neq 0$, which is a contradiction. 
 
 Next we work with $B^{\inf}(x)$. Consider $n=2\ell$. There is a unique point $(\lambda_{1}, \ldots, \lambda_{\ell}, y_{1}, \ldots, y_{\ell}) \in \mathrm{int}(\Delta_{c}^{\ell}\times\Delta_{*}^{\ell})$ such that $\sum_{j=1}^{\ell}\lambda_{j}\overline{\gamma}(y_{j})=x$. It suffices to show that the linear equation
  \begin{align}\label{gant03}
     \det 
 \begin{pmatrix}x_{1} & x_{2} & \ldots & x_{\ell+1}\\
 \vdots & & & \\
 x_{\ell+1} & x_{\ell+2} & \ldots & x_{2\ell+1}\end{pmatrix}=0,
 \end{align}
 has a unique solution $x_{2\ell+1}=\sum_{j=1}^{\ell}\lambda_{j}\gamma_{n+1}(y_{j})$.  The $j$'th column of the matrix in (\ref{gant03}),  call it $w_{j}$, $j=1, \ldots, \ell+1$, we can write as $w_{j} = AD^{j}\lambda^{T}$, where $A$ is $(\ell+1)\times\ell$ matrix with $m$'th column $(y_{m}, \ldots, y_{m}^{\ell+1})^{T}$, $m=1, \ldots, \ell$, $D$ is  $\ell\times\ell$  diagonal matrix with  diagonal entries $y_{1}, \ldots, y_{\ell}$, and $\lambda=(\lambda_{1}, \ldots, \lambda_{\ell})$. The rest of the reasoning (including the uniqueness of the solution $x_{n+1}$) is similar to the one we just discussed for $B^{\sup}$ and $n=2\ell$. 
 
 Finally, consider $n=2\ell-1$. There exists a unique point $(\beta_{2}, \ldots, \beta_{\ell}, y_{1}, \ldots, y_{\ell})\in \mathrm{int}(\Delta_{c}^{\ell-1}\times \Delta_{*}^{\ell})$ such that $\sum_{j=1}^{\ell}\beta_{j} \gamma(y_{j})=x$, where $\beta_{1}:=1-\sum_{j=2}^{\ell}\beta_{j}$. It suffices to show that the linear equation 
  \begin{align}\label{gant04}
     \det 
 \begin{pmatrix}1 & x_{1} & \ldots & x_{\ell}\\
 \vdots & & & \\
 x_{\ell} & x_{\ell+1} & \ldots & x_{2\ell}\end{pmatrix}=0,
 \end{align}
 has a unique solution $x_{2\ell}=\sum_{j=1}^{\ell}\beta_{j}\gamma_{n+1}(y_{j})$. The $j$'th column of the matrix in (\ref{gant03}),  call it $w_{j}$, $j=1, \ldots, \ell+1$, we can write as $w_{j} = AD^{j-1}\beta^{T}$, where $A$ is $(\ell+1)\times\ell$ matrix with $m$'th column $(1, y_{m}, \ldots, y_{m}^{\ell})^{T}$, $m=1, \ldots, \ell$, $D$ is  $\ell\times\ell$  diagonal matrix with  diagonal entries $y_{1}, \ldots, y_{\ell}$, and $\beta=(\beta_{1}, \ldots, \beta_{\ell})$. The rest of the reasoning (including the uniqueness of the solution $x_{n+1}$) is similar to the one we just discussed for $B^{\sup}$ and $n=2\ell$. 
 
 \subsection{The proof of Corollary~\ref{karatecor}}

 Assume contrary that there exist $n+1$ points, $\gamma(t_{1}), \ldots, \gamma(t_{n+1})$, where $a\leq t_{1}<\ldots <t_{n+1}\leq b$, which lie  in a single affine hyperplane. In particular, we have 
 \begin{align}\label{lies}
     \det(\gamma(t_{2})-\gamma(t_{1}), \gamma(t_{3})-\gamma(t_{1}), \ldots, \gamma(t_{n+1})-\gamma(t_{1}))=0.
 \end{align}
 On the other hand, we have 
 \begin{align*}
    &\det(\gamma(t_{2})-\gamma(t_{1}), \gamma(t_{3})-\gamma(t_{1}), \ldots, \gamma(t_{n+1})-\gamma(t_{1})) = \\
    &\det(\gamma(t_{2})-\gamma(t_{1}), \gamma(t_{3})-\gamma(t_{2}), \ldots, \gamma(t_{n+1})-\gamma(t_{n})) = \\
    &\int_{t_{n}}^{t_{n+1}}\ldots\int_{t_{2}}^{t_{3}}\int_{t_{1}}^{t_{2}} \det(\gamma'(s_{1}),\gamma'(s_{2}) \ldots, \gamma'(s_{n})) ds_{1}ds_{2}\ldots ds_{n} >0
 \end{align*}
 by Lemma~\ref{klasika}. Thus we have a contradiction with (\ref{lies}).

 \subsection{The proof of Corollary~\ref{provolume}}
 To prove the formulas for the volume we apply Theorem~\ref{mth010}, where $\gamma$ in Corollary~\ref{provolume} will be used as $\overline{\gamma}$ in Theorem~\ref{mth010}. Let $n=2\ell$. To verify
 \begin{align}
 \mathrm{Vol}(\mathrm{conv}(\gamma([a,b])))& \label{moculoba}\\
    =\frac{(-1)^{\frac{\ell(\ell-1)}{2}}}{(2\ell)!}& \int_{a\leq x_{1}\leq \ldots \leq x_{\ell} \leq b} \mathrm{det}(\gamma(x_{1})-\gamma(a), \ldots, \gamma(x_{\ell})-\gamma(a), \gamma'(x_{1}), \ldots, \gamma'(x_{\ell})) dx, \nonumber
\end{align}
notice that according to Theorem~\ref{mth010} the map $U:=U_{2\ell}$, where
\begin{align*}
U_{2\ell} : \Delta_{c}^{\ell}\times \Delta_{*}^{\ell} \ni (\lambda_{1}, \ldots, \lambda_{\ell}, x_{1}, \ldots, x_{\ell}) \mapsto (1-\sum_{j=1}^{\ell}\lambda_{j})\gamma(a)+\sum_{j=1}^{\ell}\lambda_{j} \gamma(x_{j}), 
\end{align*}
 is diffeomorphism between $\mathrm{int}(\Delta_{c}^{\ell}\times \Delta_{*}^{\ell})$ and $\mathrm{int}(\mathrm{conv}(\gamma([a,b])))$. In particular, by change of variables formula, we have 
 \begin{align*}
&\mathrm{Vol}(\mathrm{conv}(\gamma([a,b]))) = \int_{\Delta_{c}^{\ell}} \int_{\Delta_{*}^{\ell}} |\det(U_{\lambda_{1}}, \ldots, U_{\lambda_{\ell}}, U_{x_{1}}, \ldots, U_{x_{\ell}})| d\lambda\,  dx=\\
&\int_{\Delta_{c}^{\ell}}\lambda_{1}\ldots \lambda_{\ell} d\lambda \, \int_{\Delta_{*}^{\ell}} |\mathrm{det}(\gamma(x_{1})-\gamma(a), \ldots, \gamma(x_{\ell})-\gamma(a), \gamma'(x_{1}), \ldots, \gamma'(x_{\ell}))|dx.
 \end{align*}
 Next, using the identity 
 \begin{align}\label{distr}
     \int_{\Delta_{c}^{\ell}}\lambda^{p_{1}-1}_{1}\ldots \lambda_{\ell}^{p_{\ell-1}}(1-\sum_{j=1}^{\ell}\lambda_{j})^{p_{0}-1} d\lambda = \frac{\prod_{j=0}^{\ell} \Gamma(p_{j})}{\Gamma(\sum_{j=0}^{\ell} p_{j})}
 \end{align}
 valid for all $p_{0}, \ldots, p_{\ell}>0$ (see Dirichlet distribution in \cite{book1}), and the property 
 \begin{align*}
    &|\mathrm{det}(\gamma(x_{1})-\gamma(a), \ldots, \gamma(x_{\ell})-\gamma(a), \gamma'(x_{1}), \ldots, \gamma'(x_{\ell}))| \\ &=(-1)^{\frac{\ell(\ell-1)}{2}}\mathrm{det}(\gamma(x_{1})-\gamma(a), \ldots, \gamma(x_{\ell})-\gamma(a), \gamma'(x_{1}), \ldots, \gamma'(x_{\ell}))
 \end{align*}
 whenever $a<x_{1}<\ldots <x_{\ell}<b$, see (\ref{nishani1}), we recover (\ref{moculoba}). The other three identities in Corollary ~\ref{provolume} are obtained in the same way by repeating the  computations with $L_{2\ell}$, and in the case of odd dimensions with  $U_{2\ell-1}$ and $L_{2\ell-1}$.

 \subsection{The proof of Corollary~\ref{area1}}
 Let $n=2\ell$ (the case $n=2\ell-1$ is similar and will be omitted), and let us verify the identity 
\begin{align*}
\mathrm{Area}(\partial \; \mathrm{conv}(\gamma([a,b]))) = \frac{1}{n!} \int_{a\leq x_{1}\leq \ldots \leq x_{\ell}\leq b} \left( \sqrt{\det S_{a}^{\mathrm{Tr}}S_{a}} +\sqrt{\det S_{b}^{\mathrm{Tr}}S_{b}} \right) dx, 
\end{align*}
where $S_{r} = (\gamma(x_{1})-\gamma(r), \ldots, \gamma(x_{\ell})-\gamma(r), \gamma'(x_{1}), \ldots, \gamma'(x_{\ell}))$. By (\ref{union}) we have 
\begin{align*}
\partial\,  \mathrm{conv}(\gamma([a,b]))=\{(x,B^{\mathrm{sup}}(x)), x \in \mathrm{conv}(\bar{\gamma}([a,b]))\} \cup \{(x,B^{\mathrm{inf}}(x)), x \in \mathrm{conv}(\bar{\gamma}([a,b]))\}.
\end{align*}
On the other hand, by (\ref{giff}) and (\ref{b2l}) the set $\{(x,B^{\mathrm{sup}}(x)), x \in \mathrm{conv}(\bar{\gamma}([a,b]))\} \cap \{(x,B^{\mathrm{inf}}(x)), x \in \mathrm{conv}(\bar{\gamma}([a,b]))\}$ is contained in the image of $C^{1}$ map of the set $\partial (\Delta^{\ell}_{c}\times \Delta^{\ell}_{*})$ which has zero $n$ dimensional Lebesgue measure. Therefore, it follows from (\ref{diff2lu}) and (\ref{diff2ll}) that 
\begin{align*}
&\mathrm{Area}(\partial \; \mathrm{conv}(\gamma([a,b]))) = \\
&\mathrm{Area}(\{(x,B^{\mathrm{sup}}(x)), x \in \mathrm{conv}(\bar{\gamma}([a,b]))\}) + \mathrm{Area}(\{(x,B^{\mathrm{inf}}(x)), x \in \mathrm{conv}(\bar{\gamma}([a,b]))\}) = \\
&\int_{\Delta_{c}^{\ell}\times \Delta_{*}^{\ell}}\sqrt{ \det{ A^{\mathrm{Tr}} A}} \, dx d\lambda + \int_{\Delta_{c}^{\ell}\times \Delta_{*}^{\ell}}\sqrt{ \det{ C^{\mathrm{Tr}} C}}\,  dx d\lambda,
\end{align*}
where $A = (U_{\lambda_{1}}, \ldots, U_{\lambda_{\ell}}, U_{x_{1}}, \ldots, U_{x_{\ell}})$ with $U:=U_{n}$, and $C=(L_{\lambda_{1}}, \ldots, L_{\lambda_{\ell}}, L_{x_{1}}, \ldots, L_{x_{\ell}})$ with $L:=L_{n}$. Notice that $A^{\mathrm{Tr}}A =RS^{\mathrm{Tr}}_{b}S_{b}R$ where $R$ is $2\ell \times 2\ell$ diagonal matrix with diagonal entries $r_{1}=\ldots=r_{\ell}=1$, and $r_{\ell+1}=\lambda_{1}, \ldots, r_{\ell+\ell}=\lambda_{\ell}$. Similarly $C^{\mathrm{Tr}}C = RS_{a}^{\mathrm{Tr}}S_{a}R$. Therefore, 
\begin{align*}
&\int_{\Delta_{c}^{\ell}\times \Delta_{*}^{\ell}}\sqrt{ \det{ A^{\mathrm{Tr}} A}} \, dx d\lambda + \int_{\Delta_{c}^{\ell}\times \Delta_{*}^{\ell}}\sqrt{ \det{ C^{\mathrm{Tr}} C}}\,  dx d\lambda =\\
&\int_{\Delta_{c}^{\ell}} \lambda_{1}\cdots \lambda_{\ell} d\lambda  \int_{\Delta_{*}^{\ell}}\sqrt{\det S^{\mathrm{Tr}}_{b}S_{b}}dx + \int_{\Delta_{c}^{\ell}} \lambda_{1}\cdots \lambda_{\ell} d\lambda  \int_{\Delta_{*}^{\ell}}\sqrt{\det S^{\mathrm{Tr}}_{a}S_{a}}dx  \stackrel{(\ref{distr})}{=}\\
&\frac{1}{(2\ell)!}\int_{\Delta_{*}^{\ell}} \left(\sqrt{\det S^{\mathrm{Tr}}_{b}S_{b}} +\sqrt{\det S^{\mathrm{Tr}}_{a}S_{a}} \right) \, dx.
\end{align*}

This finishes the proof of Corollary~\ref{area1}.

\end{document}